\documentclass[12pt]{article}

\newtheorem{theorem}{Theorem}

\newtheorem{definition}{Definition}
\newtheorem{remark}{Remark}

\usepackage{amsmath,amssymb,ascmac,color,ulem}

\begin{document}
\title{ 
The Construction Problem of Algebraic Potentials
and Reflection Groups
}

\author{Jiro Sekiguchi\thanks{Department of Mathematics, Tokyo University of Agriculture and Technology}}

\maketitle

\begin{abstract}
This paper has two aims.
The first one is the construction problem of algebraic potentials of Frobenius manifolds.
We show examples of such potentials for the cases of reflection groups of types $H_4,E_6,E_7,E_8$
and also include those which are already known.
The second one is an application of such potentials to singularity theory.
We introduce families of hypersurfaces of ${\bf C}^3$
which are deformations of $E_n$-singularities  $(n=6,7,8)$
but are not the versal families of $E_n$-singularities.
We study the properties of the families.
In particular we show the correspondence between such families and the algebraic potentials constructed
in the first aim.
Moreover we discuss the relationship between the complex reflection groups $ST33$ and $ST34$
and the two families corresponding to the $E_6$-singularity and the $E_7$-singularity.
\end{abstract}

\section{Introduction}
\label{section:sec0}

In the first half of this paper, we shall show examples of algebraic potentials of Frobenius manifolds and 
related topics.
In the second half, we shall discuss applications of the algebraic potentials, mainly to singularity theory.

Not many examples of algebraic  potentials are knwon in spite of its importance.
One of the questions concerning algebraic potentials is
 the conjecture proposed by B. Dubrovin below:

\vspace{5mm}
{\bf The conjecture}:
Massive irreducible algebraic Frobenius manifolds with positive degrees $d_i$ correspond to primitive conjugacy classes in Coxeter groups.

\vspace{5mm}
This conjecture means that if the solution of the WDVV equation=the potential
is algebraic,
it  corresponds to a conjugacy class of the Coxeter group.
This is a generalization of  C. Hertling's result that 
there is a one-to-one correspondence between
the polynomial potentials
and Coxeter groups.
O. Pavlyk \cite{P1}, Y. Dinar \cite{D1}, \cite{D2}, \cite{D3} are
the basic literatures on this conjecture.
Dinar  investigated the relationship between the nilpotent orbits of
simple Lie algebras and algebraic potentials.

The author, together with collaborators  M. Kato and T. Mano,
previously defined the notion of potential vector fields,
which is a generalization of the potential,
in \cite{KMS1}, \cite{KMS2}  and calculated the potential vector fields
corresponding to the algebraic function solutions of
the Painlev\'{e} VI equation in \cite{KMS3}, \cite{KMS4}.
Recently, although limited to the real reflection groups of types $H_4,E_6,E_7,E_8$,
the author calculated some algebraic potentials.
The resulting algebraic potentials are given in
\S\ref{section:sec3} of the main text.
Various calculation methods considered in the stage of preparing
\cite{KMS3}, \cite{KMS4} are helpful in obtaining these results.

It is possible to apply algebraic potentials to other fields of mathematics  if we obtain 
their concrete forms.
Along this line, we shall discuss applications of the algebraic potentials, mainly 
to singularity theory.
We shall construct three families of hypersurfaces which are deformations
of $E_n$-singularities  $(n=6,7,8)$
but are different from the versal families of $E_n$-singularities.
The study of the families has at least three interesting consequences.
The first one is related to nilpotent orbits.
The polynomial potentials of Frobenius manifolds correspond to Coxeter-Killing conjugacy classes
of real reflection groups and in the case of Weyl groups, by the result of B. Kostant, the Coxeter-Killing classes correspond to regular nilpotent orbits of the corresponding semisimple Lie algebras.
On the other hand, polynomial potentials correspond to the versal families of
 simple singularities.
In this sense, the families introduced in this paper are regarded as those corresponding to
subregular nilpotent orbits of simple Lie algebras
of types $E_n\>(n=6,7,8)$.
The second one is the relationship between the family corresponding to $E_6$ and the complex
reflection group $ST33$ and that  between the family corresponding to $E_7$ and the complex
reflection group $ST34$
as will be shown in \S4.
The third one  is a byproduct of the second one and concerns with one of Arnold's problems.
In the book \cite{An}, p.20, it is written that

``1974-5 Find and application of the (Shephard-Todd) complex reflection groups to singularity theory''

It is possible to construct subfamilies of the families  for $E_6,\>E_7$ which are related with
ST33 group and ST34 group.
The study of these families is expected to give a new insight 
on the role of the groups $ST33$ and $ST34$ on the deformation theory of $E_6,\>E_7$-singularities.

We are going to explain the contents of this paper.
In section 2, we begin with reviewing
the WDVV equation and formulate the conjecture of Dubrovin on algebraic potentials.
We close this section with
mentioning four dimensional algebraic potentials having tri-Hamiltonian structure.
In section 3, we collect examples of algebraic potentials.
In \S3.1, we first explain an idea of constructing algebraic potentials 
which are defined by quadratic equations and using the idea we shall show examples
of algebraic potentials of types $D_4(1),\>H_4(1),\>E_n(1)\>(n=6,7,8)$.
Among these, those corresponding to $D_4(1)$ and $E_8(1)$ were constructed by O. Pavlyk \cite{P1}
and Y. Dinar (cf. \cite{DS}) but the remaining ones are new.
By using the data obtained by T. Douvropoulos in \S3.2, we shall construct examples of algebraic potentials
related to 
the reflection group of type $H_4$ in \S3.3.
From \S4, we shall discuss applications of algebraic potentials to singularity theory.
The  purpose of \S4 is to introduce three  families of hypersurfaces in ${\bf C}^3$
which are deformations of the simple singularity of types $E_n$ ($n=6,7,8$)
but are different from the versal deformations of these three singularities.
We discuss the relationship between the algebraic potentials $E_n(1)\>(n=6,7,8)$
 constructed in \S\ref{section:sec3} and these families.
 
 In the case of the reflection group of type $H_3$, there are one polynomial potential
$F_{H_3}$ (cf. (\ref{equation:potential-H3})) and two algebraic potentials $F_{(H_3)'}$, $F_{(H_3)''}$
(cf. (\ref{equation:ex-(H3)'case}),  (\ref{equation:ex-(H3)''case})).
It is known that the flat coordinates of a polynomial potential are regarded as basic invariants of the
corresponding reflection group.
Then it is a basic question whether  the flat coordinates and the algebraic functions
are written in terms of
 the basic invariants of $W(H_3)$ or not.
In section 5, we shall discuss   this question.
Another question is to construct a family of hypersurfaces
which corresponds to algebraic potentials.
We discuss this problem in the cases $F_{H_3}$ and $F_{(H_3)''}$.
As to the former, this problem is already answereded by T. Yano \cite{Ya}.
The questions similar to those explained above for the the reflection group
of type $H_4$ shall be also discussed in this section.
In the Appendix, we shall explain the idea how to obtain the coodinate transformation 
(\ref{equation:trans-e8-case}).

This is an extended version of the paper of the same title  written in Japanese
 appeared in Prceedings of the Symposium on Representation Theory 2019, 151-165 (2019).

\section{The WDVV equation}
\label{section:sec1}

We start this section with formulating the definition of
the WDVV equation.
Secondly we explain the results on the polynomial and algebraic solutions of the WDVV equation in three dimensions.
Thirdly we state the conjecture of Dubrovin concerning the algebraic solutions to the WDVV equation.
In the last part of this section, we mension the algebraic potentials having tri-Hamiltonian structure.

\subsection{Definition of the WDVV equation}
\label{subsection:sec1-1}

The WDVV equation is the system of non-linear partial differential equations  formulated by 
Witten and Dijkgraaf-Verlinde-Verlinde.
We review it briefly.
For the details, refer to
Dubrovin \cite{Du}, C. Sabbah \cite{Sab}, etc.
The formulation below is due to \cite{KMS2}.

Let $F=F_0+F_1$ be a function of
$(x_1,x_2,\cdots,x_n)$.
Here 
$F_0$ is  the polynomial defined by
\begin{equation}
\label{equation:def-F0}
F_0=\left\{\begin{array}{ll}
\displaystyle{\frac{1}{2}x_1x_n^2+\sum_{j=2}^{n/2}x_jx_{n-j+1}x_n}&(n:\>{\rm even})\\
\displaystyle{\frac{1}{2}x_1x_n^2+\sum_{j=2}^{m-1}x_jx_{2m-j}x_n+\frac{1}{2}x_{m}{}^2x_n}&
(n:\>{\rm odd},\>m=(n+1)/2)\\
\end{array}\right.
\end{equation}
and $F_1$ is a function of $(x_1,x_2,\cdots,x_{n-1})$ and is independent of $x_n$.

We assume that $F$ is weighted homogeneous.
Namely, there are non-zero constants
$w,w_1,w_2,\cdots,w_n$ such that
$EF=wF$, where
$
\displaystyle{
E=\frac{1}{w_n}\sum_{j=1}^nw_jx_j\partial_{x_j}}
$
is the Euler vector field.
The constants
$w_1,w_2,\cdots,w_n$ are assumed to be rational numbers and
$
0<w_1\le w_2\le\cdots\le w_n=1.
$
We define the vector-valued function $P$ by
$$
P=(\partial_{x_n}F,\partial_{x_{n-1}}F,\cdots,\partial_{x_1}F)
$$
From the definition,
$P$ is the gradient vector of $F$.
Moreover we put
$$
C={}^t(\partial_{x_1}P,\partial_{x_2}P,\cdots,\partial_{x_n}P),
$$
which is the Hessian of $F$.
We need the partial differentiation of $C$:
$$
\tilde{B}_j=\partial_{x_j}C\quad(j=1,2,\cdots,n)
$$
It is clear that ${\tilde B}_n=I_n$, the identity matrix.

\begin{definition}
The system of partial differential equations
$$
\left\{\begin{array}{l}
EF=wF\\
{}[\tilde{B}_j,\>\tilde{B}_k]=O\quad(j,k=1,2,\cdots,n)\\
\end{array}
\right.
$$
for $F$ is called the WDVV equation.
If $F$ is its solution, $F$ is called a potential.
Moreover the coordinate  $(x_1,\cdots,x_n)$ is called a flat coordinate.
\end{definition}

We introduce the $n\times n$ matrix $T=EC $ and $\Delta_F=\det(T)$.
In this paper,  $\Delta_F$  is called the discriminant of $F$ by the reason that
in the case where $F$ is a polynomial potential,
$\Delta_F$ is regarded as the discriminant
of the reflection group  corresponding to $F$.

\subsection{Three dimensional case}
\label{subsection:sec1-2}
\subsubsection{Polynomial potentials}
\label{subsubsection:sec1-2-1}

It is important but difficult to find the potentials.
B. Dubrovin came up a very interesting conclusion in the course of
 finding  polynomial potentials when $n=3$
to understand what the WDVV equation is.
(Cf. \cite{Du})
The polynomial potential means that $F$ is a polynomial and is a solution of the WDVV equation.
Among  such potentials, Dubrovin extracted only three of them.
These are
\begin{equation}
\label{equation:potential-A3}
F=\displaystyle{\frac{x_1x_3^2+x_2^2x_3}{2}+\frac{x_1^2x_2^2}{4}+\frac{x_1^5}{60},}
\end{equation}
\begin{equation}
\label{equation:potential-B3}
F=\displaystyle{\frac{x_1x_3^2+x_2^2x_3}{2}+\frac{x_1x_2^3}{6}+
\frac{x_1^3x_2^2}{6}+\frac{x_1^7}{210},}
\end{equation}
\begin{equation}
\label{equation:potential-H3}
F=\displaystyle{\frac{x_1x_3^2+x_2^2x_3}{2}+\frac{x_1^2x_2^3}{6}+
\frac{x_1^5x_2^2}{20}+\frac{x_1^{11}}{3960}}.
\end{equation}

We easily see that the weights
$w_1,w_2,w_3$ in these cases are as follows.

\vspace{5mm}
\centerline{TABLE I}
$$
\begin{array}{c|c|c}\hline
w_1&w_2&w_3\\\hline\hline
1/2&3/4&1\\
1/3&2/3&1\\
1/5&3/5&1\\\hline
\end{array}
$$
Then we observe the relationship between these numbers and the degrees of basic
invariants of the reflection groups of types $A_3,B_3,H_3$.
Let
$d_1,d_2,d_3\>(d_1<d_2<d_3)$ be such degrees.
They are given in TABLE II:

\vspace{5mm}
\centerline{TABLE II}
$$
\begin{array}{c||c|c|c}\hline
&d_1&d_2&d_3\\\hline\hline
A_3&2&3&4\\
B_3&2&4&6\\
H_3&2&6&10\\\hline
\end{array}
$$
Comparison of the data for $(w_1,w_2,w_3)$ with those for $(d_1,d_2,d_3)$ shows 
that
the ratios of both coincide with each other.

Use the notation of 
\S\ref{subsection:sec1-1} as is and set $T=EC$.
Then
$\det(T)$ is considered the discriminant
of the corresponding reflection group.
We explain its meaning.
Let $W$ be an irreducible reflection group and let $P_1,\cdots,P_n$ be the basic invariants of
$W$.
We may assume that each $P_j$ is homogeneous and set $d_j=\deg(P_j)$.
Furthermore we may assume that $2=d_1\le d_2\le \cdots\le d_n$.
Then there is a weighted homogeneous polynomial $\delta(x_1,\cdots,x_n)$
such that the discriminant of $W$
coincides with $\delta(P_1,\cdots,P_n)$.
There is one thing to say here.
The  polynomial  $\delta(x_1,\cdots,x_n)$ depends on the choice
of basic invariants $P_1,\cdots,P_n$.
The meaning of ``$\det(T)$ is considered the discriminant
of the corresponding reflection group'' is 
that if the type of $W$ is one of $A_3,\>B_3,\>H_3$,
then  there exist a polynomial $\delta(x_1,x_2,x_3)$
and  a set of basic invariants $\{P_1,P_2,P_3\}$
such that $\delta(P_1,P_2,P_3)$
is the discriminant of $W$
and that $\det(T)$ coincides with $\delta(x_1,x_2,x_3)$ up to a constant factor.
This suggests the existence of a special kind of a generator system of a reflection group
corresponding to the potential $F$ of (\ref{equation:potential-A3}),  (\ref{equation:potential-B3}),
 (\ref{equation:potential-H3}).
Dubrovin found the literature \cite{SYS}
which treated  such a generator system, called a flat generator system.

\begin{remark}
K. Saito already recognized the existence of $F$ at least in the cases $A_3,\>B_3,\>H_3$
before Dubrovin's study (unpublished).
\end{remark}
Dubrovin developed further his investigation and constructed
algebraic solutions to the
Painlev\'{e} VI equation from the three polynomials in 
(\ref{equation:potential-A3}),  (\ref{equation:potential-B3}),
 (\ref{equation:potential-H3}).

The study of Dubrovin leads us to formulate two problems:

Problem 1.
(Dubrovin's conjecture) Is there a one to one correspondence
between the totality of the polynomial potentials and that of real reflection groups?

Problem 2. Is it possible to generalize the existence of flat generator systems
to the case of complex reflection groups?

\begin{remark}
Problem 1 is answered by C. Hertling affirmatively.

As to Problem 2, refer to \cite{KMS1}, \cite{KMS2}.
\end{remark}

\subsubsection{The study by Dubrovin-Mazzocco }
\label{subsubsection:sec1-2-2}

Dubrovin gave a one-to-one correspondence between the WDVV equation
for $n=3$ and the 1-parameter family of the
Pianlev\'{e} VI equations.
This 1-parameter family is displayed as Painlev\'{e} VI$\mu$.
Dubrovin-Mazzocco\cite{DM} studied the solutions contained in
Painlev\'{e} VI$\mu$.
And in addition to the polynomial solutions,
they showed that there are two kinds of algebraic function solutions of the Painlev\'{e} VI equations.
One display of the algebraic function solutions they found would be
10 pages when printed in normal character size.
P. Boalch succeeded in getting a compact display of one line.
These two are related to the reflection group of type $H_3$.
The potential corresponding to the above mentioned
real reflection group of type $H_3$ is represented by $(H_3)$,
and the remaining two are $(H_3)'$, $(H_3)''$.
Their concrete forms are  (cf. \cite{KMS2})
\begin{equation}
\label{equation:ex-(H3)'case}
F_{(H_3)'}=\displaystyle{\frac{x_1x_3^2+x_2^2x_3}{2}-
\frac{x_1^4z}{6}-\frac{7x_1^3z^4}{72}-\frac{17x_1^2z^7}{315}-\frac{2x_1z^{10}}{81}
-\frac{64z^{13}}{15795}}
\end{equation}
(Here $z$ is an algebraic function of $x_1,x_2$ defined by the equation
$
3x_2-3x_1z-z^4=0.
$)
\begin{equation}
\label{equation:ex-(H3)''case}
F_{(H_3)''}=\displaystyle{\frac{x_1x_3^2+x_2^2x_3}{2}+
\frac{4063x_1^7}{1701}+\frac{19x_1^5z^2}{135}-\frac{73x_1^3z^4}{27}+\frac{11x_1z^{6}}{9}
-\frac{16z^{7}}{35}}
\end{equation}
(Here $z$ is an algebraic function of $x_1,x_2$ defined by the equation
$-x_2-x_1^2+z^2=0$.)

\subsection{The conjecture of Dubrovin}
\label{subsection:sec1-3}

Dubrovin  proposed many conjectures.
The Dubrovin conjecture here is:

\vspace{5mm}
{\bf The conjecture}:
Massive irreducible algebraic Frobenius manifolds with positive degrees $d_i$ correspond to primitive conjugacy classes in Coxeter groups.

\vspace{5mm}
I will explain the terminology of this conjecture.
``Frobenius manifold'' can be thought of a synonym for  ``a solution
of the WDVV equation=potential''.
The number ``$d_i$'' is the same as $w_i$.
``massive'' means that the matrix $T$ is generically semisimple.
The Coxeter group here means a finite Coxeter group.
``The primitive conjugacy class'' is a conjugacy class
whose original Coxeter group is generated by
the entire reflections
that appear in the display when the representative element is
represented by the product of reflections.
What is an ``algebraic Frobenius manifold''?
As with the polynomial solution of the WDVV equation,
it can be understood as an algebraic solution.
However, when it comes to what an algebraic solution is, its meaning is subtle.
Apart from that, the background that Dubrovin formulated this conjecture was
Drinfel'd-Sokolov reduction,
Drinfel'd-Sokolov hierarchy(see \cite{DS}).
Dubrovin expanded his perspective to the general case from
the result of the real reflection groups of rank 3.
A typical example of primitive conjugacy class is that of the Coxeter-Killing element.
It is known that there is a close relationship
between the order of the basic invariants and the eigenvalues of 
the Coxeter-Killing element.
In particular, a duality relationship holds.
If we use the weights $w_1,w_2,\cdots,w_n(=1)$,
then 
$w_j+w_{n-j+1}=1+w_1\>(j=1,2,\cdots,n)$ holds.
This duality is important and
holds in the case of  other conjugacy classes.

\subsection{Algebraic potentials having tri-Hamiltonian structure}
\label{subsection:sec1-4}

The algebraic function solution of the WDVV equation is called the algebraic potential.
In fact, $(H_3)',(H_3)''$ are examples of algebraic potentials.
According to the result of Dubrovin-Mazzocco,
the algebraic potentials of the three variables are three kinds of polynomial potentials
and
$(H_3)',(H_3)''$ and there are five kinds in total.
On the other hand, we find the following:

\begin{tabular}{l}
There is one kind of primitive conjugacy class of the reflection group of type $A_3$.\\
There is one kind of primitive conjugacy class of the reflection group of type $B_3$.\\
There are three kinds of primitive conjugacy classes of the reflection group of type $H_3$.\\
\end{tabular}\\

In this way, we can see that the Dubrovin conjecture holds in the three dimensional case.

Higher rank case is an open problem.
One of interesting results in four dimensional case is treated by S. Romano \cite{Rom1},
which we are going to explain briefly.
We consider the condition

\vspace{2mm}
\noindent
(TH) $n=2k$ and $0<2\mu+1=w_1=\cdots=w_k<w_{k+1}=\cdots=w_{2k}=1$.

\vspace{3mm}
If  (TH) holds, then the  Frobenius manifold defined by $F$ has tri-Hamiltonian structure
(cf. \cite{Rom1}).
S. Romano \cite{Rom1} proved a complete description of semisimple tri-Hamiltonian Frobenus manifolds
in the lowest non-trivial dimension $n=4$.
His approach relies on the interpretation of semisimple Frobenius manifolds
as isomonodromic deformation spaces of certain Fuchsian linear systems 
which are equivalent to the Painlev\'{e} VI$\mu$ equation
\begin{equation}
\label{equation:Painleve-DM}
\begin{array}{lll}
\displaystyle{\frac{d^2y}{dt^2}}&=&
\displaystyle{\frac{1}{2}\left(\frac{1}{y}+\frac{1}{y-1}+\frac{1}{y-t}\right)\left(\frac{dy}{dt}\right)^2-
\left(\frac{1}{t}+\frac{1}{t-1}+\frac{1}{y-t}\right)\left(\frac{dy}{dt}\right)}\\[4mm]
&&\displaystyle{\hspace{10mm}+\frac{y(y-1)(y-t)}{2t^2(t-1)^2}\left((2\mu-1)^2+\frac{t(t-1)}{(y-t)^2}\right).}\\
\end{array}
\end{equation}
The above special case of the Painlev\'{e} VI equation was appeared in \cite{DM}.
It is shown in \cite{DM} that the solutions to
(\ref{equation:Painleve-DM}) parametrize semisimple Frobenius manifolds of dimension 3.
This means the existence of a map of  3-dimensional Frobenius manifolds to 4-dimensional tri-Hamiltonian ones,
since both classes are parametrized by the same PVI transcendents.

O. Pavlyk's paper \cite{P1} is one of the fundamental
papers in the study of higher dimensional case.
In particular
he gave the following example for $D_4$:
$$
F=\frac{1}{2}x_1x_4^2 + x_2x_3x_4 
-\frac{19}{15}x_1^5 +\frac{14}{3}x_1^3x_2^2 - x_1x_2^4 - \frac{4}{3}x_1^3x_3 - 
 4x_1x_2^2x_3 - 
  \frac{1}{2}x_1x_3^2 +\frac{8}{15}z^5,
$$
where $z$ is an algebraic function of $x_1,x_2,x_3$ defined by the equation
$$
-x_1^2+x_2^2+x_3+z^2=0.
$$
In this example,
$w_1=w_2=1/2,\>w_3=w_4=1$,
which shows that his example is regarded the first example of semisiple Frobenius manifold
having tri-Hamiltonian structure.

Moreover the papers by Y. Dinar
\cite{D1},\cite{D2},\cite{D3} are also fundamental in this direction.
He calculated an example in the case of type $F_4$ in \cite{D1}:
\begin{equation}
\label{equation:pvf-F4}
\begin{array}{ll}
F=&\displaystyle{\frac{1}{2}x_1x_4^2+
x_2x_3x_4 +
\frac{17314}{35}x_1^7 + 1822x_1^6x_2 + 1992x_1^5x_2^2 + 650x_1^4x_2^3 - 
  290x_1^3x_2^4 - 
    \frac{1374}{5}x_1^2x_2^5}\\[4mm]
&\displaystyle{ -\frac{356}{5}x_1x_2^6
 -\frac{214}{35}x_2^7 + 
 \frac{3}{2}(x_1 - x_2)(x_1 + x_2)^2(131x_1^3 + 195x_1^2x_2 + 93x_1x_2^2 + 
    13x_2^3)z }\\[4mm]
&\displaystyle{   -\frac{9}{2}(x_1 + x_2)^4(5x_1 + 4x_2)z^2}
\displaystyle{+\left(\frac{125}{2}x_1^4 + 20x_1^3x_2 - 
  69x_1^2x_2^2 - 
     52x_1x_2^3 -\frac{ 19}{2}x_2^4\right)z^3}\\[4mm]
&\displaystyle{-3(x_1 + x_2)^2(5x_1 + 4x_2)z^4 -\frac{27}{10}(x_1 + x_2)^2z^5-
  \frac{1}{2}(5x_1 + 4x_2)z^6 -\frac{9}{14}z^7},\\
\end{array}
\end{equation}
where $z$ is an algebraic function of $x_1,x_2,x_3$ defined by the equation
\begin{equation}
\label{equation:def-z--F4(1)}
2(2x_1 + x_2)(x_1^2 + 4x_1x_2 + x_2^2) - x_3 + 3(x_1 + x_2)^2z + z^3=0.
\end{equation}
In this case,
$w_1=w_2=1/3,\>w_3=w_4=1$.
As a  consequence, this gives
 the second example of semisimple Frobenius manifold having tri-Hamiltonian structure.

In the case of type $H_4$, we found the following example of an algebraic potential
corresponding to
a semisimple Frobenius manifold
having tri-Hamiltonian structure  with the aid of
the data by  T. Douvropoulos which will be mentioned later (see \S\ref{section:sec3}):
{\footnotesize
$$
\begin{array}{ll}
F=&
\displaystyle{\frac{1}{198}(-45669270x_1^{11} - 858257224x_1^{10}x_2 - 1955070535x_1^9x_2^2 
- 
        2004227280x_1^8x_2^3 - 1191552120x_1^7x_2^4 }\\
&- 458678880x_1^6x_2^5 - 
        120561210x_1^5x_2^6 - 21740400x_1^4x_2^7 - 2569050x_1^3x_2^8 - 
        178200x_1^2x_2^9\\
& - 4455x_1x_2^{10} + 198x_2x_3x_4 + 99x_1x_4^2)\\
& + 
        5x_1(x_1 + x_2)^2(3x_1 + x_2)^2 (9035x_1^5 + 32316x_1^4x_2 + 
         30270x_1^3x_2^2 + 11760x_1^2x_2^3 + 2115x_1x_2^4 + 180x_2^5)z\\
& -
     \frac{5}{6}(994315x_1^9 + 4563564x_1^8x_2 + 7356636x_1^7x_2^2 + 
      6239424x_1^6x_2^3 + 3321450x_1^5x_2^4 + 1252440x_1^4x_2^5\\
& +356940x_1^3x_2^6 + 
        73440x_1^2x_2^7 + 9315x_1x_2^8 + 540x_2^9)z^2\\
& + \frac{5}{3}(x_1 + x_2)(3x_1 + x_2)(15815x_1^6 + 37788x_1^5x_2 + 58125x_1^4x_2^2 + 
        42360x_1^3x_2^3 + 13185x_1^2x_2^4\\
& + 1260x_1x_2^5 - 45x_2^6)z^3\\
& - 
       100(1195x_1^7 + 2450x_1^6x_2 + 1283x_1^5x_2^2 - 650x_1^4x_2^3 - 
       1015x_1^3x_2^4 - 450x_1^2x_2^5 - 87x_1x_2^6 - 6x_2^7)z^4\\
& + 
       (-74420x_1^6 - 181952x_1^5x_2 - 174765x_1^4x_2^2 - 80920x_1^3x_2^3 - 
        17130x_1^2x_2^4 - 840x_1x_2^5 + 135x_2^6)z^5 \\
&- \frac{35}{3}
     (x_1 + x_2)(3x_1 + x_2)(305x_1^3 + 336x_1^2x_2 + 123x_1x_2^2 + 12x_2^3)z^6\\
& - 
       10(175x_1^4 + 344x_1^3x_2 + 236x_1^2x_2^2 + 64x_1x_2^3 + 5x_2^4)z^7\\
& - 
       5(95x_1^3 + 124x_1^2x_2 + 57x_1x_2^2 + 8x_2^3)z^8\\
& - 
       \frac{125}{9}(x_1 + x_2)(3x_1 + x_2)z^9 + \frac{1}{2}(-5x_1 - 4x_2)z^{10}
 -\frac{25}{33}z^{11}\\
\end{array}
$$
}
where $z$ is an algebraic function of $x_1,x_2,x_3$
defined by the equation
$$
\begin{array}{l}
843x_1^5 + 1160x_1^4x_2 + 530x_1^3x_2^2 + 120x_1^2x_2^3 + 15x_1x_2^4 - x_3 - 
15(x_1 + x_2)^2(3x_1 + x_2)^2z\\ 
+ 20x_1(13x_1^2 + 12x_1x_2 + 3x_2^2)z^2
 +  
10(x_1 + x_2)(3x_1 + x_2)z^3 + z^5=0.\\
\end{array}
$$
In this example, $w_1=w_2=1/5,\>w_3=w_4=1$.

It is interesting to  apply the  result by Romano explained above to  these three cases.
Then it is expected to show the correspondences between
 polynomial potentials in three dimensions constructed by B. Dubrovin
and four dimensional algebraic potentials having
tri-Hamiltonian structure:

(C-A3) the polynomial potential of type $A_3$ and the algebraic potential of type $D_4$,

(C-B3) the polynomial potential of type $B_3$ and the algebraic potential of tyep $F_4$,

(C-H3) the polynomial potential of type $H_3$ and the algebraic potential  of $H_4$.

The author studied these three cases in \cite{Se3}.

\section{Examples of algebraic potentials}
\label{section:sec3}

In this section, we collect examples of algebraic potentials related to exceptional real reflection groups.

\subsection{Examples  (I)}
\label{subsection:sec3-1}

In spite that it is not fully understood the meaning of that  the potential is algebraic,
the examples in the previous section
lead us to formulate a method of constructing algebraic potentials.
We explain it.
First, let $W$ be an irreducible real reflection group of rank $n$.
Let $d_1,d_2,\cdots,d_n$ be the degrees of the basic invariants of $W$
such that $0<d_1\le d_2\le \cdots\le d_n$.
We choose an integer $d_0$ $(d_1\le d_0\le d_{n-1}$) so that there is a duality among
$d_1,d_2,\cdots,d_{n-1},d_0$, namely,
there is an involution $\ast$ of the set $\{1,2,\cdots,n-1,0\}$
such that $d_j+d_{j^{\ast}}=d_1+d_{n-1}$ for all $j\in\{1,2,\cdots,n-1,0\}$.

We explain how to choose $d_0$ by taking the reflection group of type $F_4$ case as an example.
In this case $n=4$ and $d_1=2,d_2=6,d_3=8,d_4=12$.
Then it is easy to see that $d_0=4$.

We change the indices  of $d_1,\cdots, d_{n-1}, d_0$ to
 $d_1',\cdots,d_{n-1}',d_n'$ so that $d_1'\le d_2'\le \cdots\le d_n'$.
Setting $w_j=d_j'/d_{n}'\>(j\in\{1,2,\cdots,n\})$,
we take the variables
$x_1,x_2,\cdots,x_{n-1},x_n$ so that $w_j$ is the weight of the variable $x_j$ ($j=1,2,\cdots,n$).
Let $z$ 
be an algebraic function of $x_1,\cdots,x_{n-1}$ 
defined by the relation
\begin{equation}
\label{equation:def-z}
z^2+v(x_1,\cdots,x_{n-3},x_{n-2})-x_{n-1}=0,
\end{equation}
where $v(x_1,\cdots,x_{n-3},x_{n-2})$ is a weighted homogeneous polynomial of
$x_1,\cdots,x_{n-3},x_{n-2}$
and $z$ is of weight $w_{n-1}/2$.

To define a candidate of the potential, we recall the polynomial $F_0$ defined in (\ref{equation:def-F0}).
Let $Z$ be the variable of weight $d_{n-1}/2$
and let ${\tilde F}_1={\tilde F}_1(x_1,\cdots,x_{n-2},x_{n-1},Z)$
be a weighted homogeneous polynomial of
$(x_1,\cdots,x_{n-2},x_{n-1},Z)$ of weight $w_1+2$ which is the same as the weight of $F_0$.
Then $F=F_0+({\tilde F}_1|_{Z=z})$ is the candidate of the potential,
that is, the solution of the WDVV equation.

Under the assumption on ${\tilde F}_1$,
${\tilde F}_1$ may not contain $x_{n-1}$ because\\
$x_{n-1}=z^2+v(x_1,\cdots,x_{n-3},x_{n-2})$ holds.
Then we may assume that
${\tilde F}_1$ takes  the form
\begin{equation}
\label{equation:def-F1}
{\tilde F}_1=u_0+u_1Z+u_2Z^2+\cdots+u_pZ^p,
\end{equation}
where $u_0,u_1,\cdots,u_p$ are polynomials of $x_1,\cdots,x_{n-2}$.

To compute the partial derivatives of $F$ with respect to the variables $x_1,\cdots,x_n$,
we need some preparations.
We assume that
\begin{equation}
\label{equation:equation-u}
u_1=u_3=0.
\end{equation}
It is easy to show that if $u_1=0$, there exists a weighted homogeneous 
polynomial ${\tilde G}$ of $x_1,\cdots,x_{n-2},Z$ such that
$\partial_Z{\tilde F}_1=Z{\tilde G}$.
Moreover if $u_1=u_3=0$, then in addition to the existence of ${\tilde G}$,
there  exists a weighted homogeneous 
polynomial ${\tilde H}$ of $x_1,\cdots,x_{n-2},Z$ such that
$\partial_Z{\tilde G}=Z{\tilde H}$.
It follows from
(\ref{equation:def-z}) that
\begin{equation}
\partial_{x_j}z=-\frac{\partial_{x_j}v}{2z}\quad (j=1,\cdots,n-2),
\quad\partial_{x_{n-1}}z=\frac{1}{2z},
\quad
\partial_{x_n}z=0.
\end{equation}

\begin{remark}
There is no theoretical reason to assume 
(\ref{equation:equation-u}) for the polynomial ${\tilde F}_1$.
But under the assumption  (\ref{equation:equation-u}),
it is possible to show the existence of the polynomials ${\tilde H}$,
${\tilde G}$ of $x_1,\cdots,x_{n-2},Z$ and the polynomial ${\tilde G}_1$ of $x_1,\cdots,x_{n-2}$
such that 
$$
\int_0^ZZ{\tilde H}\>dZ={\tilde G},\quad
\int_0^ZZ({\tilde G}+{\tilde G}_1)dZ={\tilde F}_1-u_0.
$$
It is interesting to  understand the roles of  ${\tilde H},\>{\tilde G}$.
\end{remark}

\vspace{5mm}
In the below, we focus our attention to  the real reflection groups of types
$H_3,D_4,F_4,H_4,$ $E_6,E_7,E_8$.
We are going to determine 
$d_0$ in these cases.

\begin{center}
TABLE III
\end{center}

$$
\begin{array}{l||l||l|l|l}\hline
{\rm type}&&&d_0&\deg z\\\hline\hline
H_3&1,3,5&1,3&2&1\\
D_4&2,4,6,4&2,4,4&2&2\\
F_4&2,6,8,12&2,6,8&4&3\\
H_4&2,12,20,30&2,12,20&10&6\\
E_6&2,5,6,8,9,12&2,5,6,8,9&3&4\\
E_7&2,6,8,10,12,14,18&2,6,8,10,12,14&4&6\\
E_8&2,8,12,14,18,20,24,30&2,8,12,14,18,20,24&6&10\\\hline
\end{array}
$$

We are going to show  the algebraic potentials
satisfying the conditions above.
For this purpose, we put
$$
{\tilde{\sf{X}}}=\{H_3,\>D_4,\>F_4,\>H_4,\>E_6,\>E_7,\>E_8\}
$$
and for $X\in{\tilde X}$, put
$$
W_X=\{w_1,w_2,\ldots,w_n\}
$$
using the notation defined at the begining of this subsection.

\begin{theorem}
\label{theorem:theorem7}
For each $X\in{\tilde{\sf{X}}}$, there is an algebraic potential $F=F(x_1,\ldots,x_n)$ of type $X$ 
satisfying the assumption
(\ref{equation:equation-u}), that $w_j(\in W_X)$ is the weight of the variable $x_j$ $(j=1,\ldots,n)$ and
that $\det(T)$ is not zero
and no non-trivial multiple factor.
In particular the concrete form of $F$ is given as follows.

(i) \underline{The case $H_3$}

The result is the same as the case
$(H_3)''$ (cf. (\ref{equation:ex-(H3)''case})).

(ii) \underline{The case $D_4$}

This is the same as the case of $D_4$ treated in the previous section.

(iii) \underline{The case $F_4$}

This case is reduced to the polynomial potential of the reflection group
of type $B_4$.

(iv) \underline{The case $H_4$}
$$
\begin{array}{ll}
F=&\displaystyle{ x_2x_3x_4 + \frac{1}{2}x_1x_4^2} \\
&\displaystyle{+\frac{4 }{13167}{x_1} \left(51139 {x_1}^{20}-175560 {x_1}^{15}
   {x_2}-95760 {x_1}^{10} {x_2}^2+1790712 {x_1}^5
   {x_2}^3-13167 {x_2}^4\right)}\\
& \displaystyle{  -\frac{4}{45}  \left(17
   {x_1}^{15}+75 {x_1}^{10} {x_2}-900 {x_1}^5 {x_2}^2+30
   {x_2}^3\right)z^2}\\
   &\displaystyle{-\frac{1}{18} {x_1}^4  \left(17 {x_1}^5-90
   {x_2}\right)z^4+\frac{1}{6}{x_1}^3 z^6+\frac{1}{105}z^7}\\
   \end{array}
   $$
Here $z$ is an algebraic function of $x_1,x_2,x_3$ defined by the equation
$$
-4x_1(x_1^5 - 3x_2) - x_3 + z^2=0.
$$
In this case,
$w_1=1/10,w_2=1/2,w_3=3/5,w_4=1$.

(v) \underline{The case $E_6$}
{\footnotesize
$$
\begin{array}{ll}
F=&\displaystyle{ x_3x_4x_6 +
  x_2x_5x_6 +\frac{ x_1x_6^2}{2 }}\\[4mm]
&\displaystyle{+\frac{2106x_1^{10}}{5 }+\frac{ 8100x_1^7x_2^2}{7 }+
  348x_1^4x_2^4 +\frac{ 104x_1x_2^6}{3 }+
     54x_1^6x_2x_3 -
  84x_1^3x_2^3x_3 -\frac{ 44x_2^5x_3}{5 }+\frac{ 648x_1^5x_3^2}{5 }}\\[4mm]
&\displaystyle{+
     27x_1^2x_2^2x_3^2 + 6x_1x_2x_3^3 +\frac{ 3x_3^4}{8 }+
 54x_1^7x_4 -
  180x_1^4x_2^2x_4 -
     20x_1x_2^4x_4 +
 126x_1^3x_2x_3x_4 +4x_2^3x_3x_4}\\[4mm]
& \displaystyle{+
  18x_1^2x_3^2x_4 +
     27x_1^4x_4^2 + 10x_1x_2^2x_4^2 +x_2x_3x_4^2 +\frac{ 5x_1x_4^3}{3 }
-
  27x_1^6x_5 +
     90x_1^3x_2^2x_5 + 6x_2^4x_5}\\[4mm]
&\displaystyle{ -
  27x_1^2x_2x_3x_5 +\frac{ 9x_1x_3^2x_5 }{2}+ 9x_1^3x_4x_5 -
     4x_2^2x_4x_5 +\frac{ x_4^2x_5}{2 }+\frac{ 9x_1^2x_5^2}{4 }+\frac{ 
     2z^5}{5}}\\
\end{array}
$$
}
Here $z$ is an algebraic function of $x_1,x_2,x_3,x_4,x_5$
defined by the equation
\begin{equation}
\label{equation:def-ans-z-e6e}
v - x_5 + z^2=0,
\end{equation}
where
\begin{equation}
\label{equation:def-ans-z-e6v}
v=-9x_1^4 - 12x_1x_2^2 + 2x_2x_3 + 2x_1x_4.
\end{equation}
In this case,
$w_1=2/9,w_2=1/3,w_3=5/9,w_4=2/3,w_5=8/9,w_6=1$.

(vi) \underline{The case $E_7$}
{\footnotesize
$$
\begin{array}{ll}
&F\\
=&\frac{{x_7}^2 {x_1}}{2}+\frac{{x_4}^2
   {x_7}}{2}+{x_3} {x_5} {x_7}+{x_2} {x_6} {x_7}\\
&+\frac{3939238656 {x_1}^{15}}{1092455}-\frac{10368 {x_2}
   {x_1}^{13}}{7}+\frac{3456 {x_3} {x_1}^{12}}{7}+\frac{18166464 {x_2}^2
   {x_1}^{11}}{539}-\frac{1728 {x_4} {x_1}^{11}}{7}-\frac{19008}{7} {x_2}
   {x_3} {x_1}^{10}\\
   &-\frac{288 {x_5} {x_1}^{10}}{7}+\frac{216576 {x_2}^3
   {x_1}^9}{7}+\frac{643392 {x_3}^2 {x_1}^9}{49}-\frac{4608}{7} {x_2}
   {x_4} {x_1}^9-\frac{288 {x_6} {x_1}^9}{7}-\frac{1728}{7} {x_2}^2
   {x_3} {x_1}^8\\
   &+\frac{864}{7} {x_3} {x_4} {x_1}^8-\frac{432}{7}
   {x_2} {x_5} {x_1}^8+60264 {x_2}^4 {x_1}^7-\frac{115776}{7} {x_2}
   {x_3}^2 {x_1}^7+\frac{72360 {x_4}^2 {x_1}^7}{343}+3024 {x_2}^2
   {x_4} {x_1}^7\\
   &+\frac{144}{7} {x_3} {x_5} {x_1}^7-\frac{432}{7}
   {x_2} {x_6} {x_1}^7+\frac{42528 {x_3}^3 {x_1}^6}{7}+19152 {x_2}^3
   {x_3} {x_1}^6+\frac{36864}{7} {x_2} {x_3} {x_4} {x_1}^6+576
   {x_2}^2 {x_5} {x_1}^6\\
   &-\frac{72}{7} {x_4} {x_5}
   {x_1}^6+\frac{144}{7} {x_3} {x_6} {x_1}^6+\frac{252252 {x_2}^5
   {x_1}^5}{5}+18216 {x_2}^2 {x_3}^2 {x_1}^5+\frac{2844}{7} {x_2}
   {x_4}^2 {x_1}^5+\frac{942 {x_5}^2 {x_1}^5}{245}\\
   &+4824 {x_2}^3
   {x_4} {x_1}^5-\frac{10944}{7} {x_3}^2 {x_4} {x_1}^5+360 {x_2}
   {x_3} {x_5} {x_1}^5+576 {x_2}^2 {x_6} {x_1}^5-\frac{72}{7}
   {x_4} {x_6} {x_1}^5\\
   &-7008 {x_2} {x_3}^3 {x_1}^4+\frac{36}{7}
   {x_3} {x_4}^2 {x_1}^4+23940 {x_2}^4 {x_3} {x_1}^4+1728
   {x_2}^2 {x_3} {x_4} {x_1}^4+420 {x_2}^3 {x_5}
   {x_1}^4\\
   &+\frac{1632}{7} {x_3}^2 {x_5} {x_1}^4+\frac{456}{7} {x_2}
   {x_4} {x_5} {x_1}^4+360 {x_2} {x_3} {x_6}
   {x_1}^4-\frac{12}{7} {x_5} {x_6} {x_1}^4+38220 {x_2}^6
   {x_1}^3+3032 {x_3}^4 {x_1}^3\\
   &-\frac{102 {x_4}^3 {x_1}^3}{7}+4116
   {x_2}^3 {x_3}^2 {x_1}^3+444 {x_2}^2 {x_4}^2 {x_1}^3-\frac{27}{7}
   {x_2} {x_5}^2 {x_1}^3+\frac{54 {x_6}^2 {x_1}^3}{49}+5082 {x_2}^4
   {x_4} {x_1}^3\\
   &+264 {x_2} {x_3}^2 {x_4} {x_1}^3+156 {x_2}^2
   {x_3} {x_5} {x_1}^3-\frac{276}{7} {x_3} {x_4} {x_5}
   {x_1}^3+420 {x_2}^3 {x_6} {x_1}^3+48 {x_3}^2 {x_6} {x_1}^3\\
   &+24
   {x_2} {x_4} {x_6} {x_1}^3+1988 {x_2}^2 {x_3}^3 {x_1}^2-150
   {x_2} {x_3} {x_4}^2 {x_1}^2+\frac{25}{7} {x_3} {x_5}^2
   {x_1}^2+24402 {x_2}^5 {x_3} {x_1}^2\\
   &-544 {x_3}^3 {x_4}
   {x_1}^2+2436 {x_2}^3 {x_3} {x_4} {x_1}^2+357 {x_2}^4 {x_5}
   {x_1}^2-156 {x_2} {x_3}^2 {x_5} {x_1}^2-\frac{3}{7} {x_4}^2
   {x_5} {x_1}^2\\
   &+48 {x_2}^2 {x_4} {x_5} {x_1}^2+84 {x_2}^2
   {x_3} {x_6} {x_1}^2+\frac{156}{7} {x_3} {x_4} {x_6}
   {x_1}^2-\frac{6}{7} {x_2} {x_5} {x_6} {x_1}^2+5439 {x_2}^7
   {x_1}-812 {x_2} {x_3}^4 {x_1}\\
   &-6 {x_2} {x_4}^3 {x_1}+6468
   {x_2}^4 {x_3}^2 {x_1}+147 {x_2}^3 {x_4}^2 {x_1}+54 {x_3}^2
   {x_4}^2 {x_1}+\frac{5}{2} {x_2}^2 {x_5}^2 {x_1}-\frac{5}{14}
   {x_4} {x_5}^2 {x_1}\\
   &+\frac{9}{7} {x_2} {x_6}^2
   {x_1}+1764 {x_2}^5 {x_4} {x_1}+42
   {x_2}^2 {x_3}^2 {x_4} {x_1}+\frac{160}{3} {x_3}^3 {x_5}
   {x_1}+140 {x_2}^3 {x_3} {x_5} {x_1}\\
   &+20 {x_2} {x_3} {x_4}
   {x_5} {x_1}+294 {x_2}^4 {x_6} {x_1}-12 {x_2} {x_3}^2
   {x_6} {x_1}+\frac{9}{7} {x_4}^2 {x_6} {x_1}+36 {x_2}^2 {x_4}
   {x_6} {x_1}-\frac{10}{7} {x_3} {x_5} {x_6} {x_1}\\
   &+\frac{1484
   {x_3}^5}{15}+\frac{1666 {x_2}^3 {x_3}^3}{3}+{x_3}
   {x_4}^3+\frac{{x_5}^3}{84}+7 {x_2}^2 {x_3} {x_4}^2-\frac{1}{2}
   {x_2} {x_3} {x_5}^2+\frac{{x_3} {x_6}^2}{7}+1029 {x_2}^6
   {x_3}\\
   &+\frac{364}{3} {x_2} {x_3}^3 {x_4}+392 {x_2}^4 {x_3}
   {x_4}+\frac{49 {x_2}^5 {x_5}}{5}+21 {x_2}^2 {x_3}^2
   {x_5}-\frac{1}{2} {x_2} {x_4}^2 {x_5}+7 {x_2}^3 {x_4} {x_5}-6
   {x_3}^2 {x_4} {x_5}\\
   &-\frac{20 {x_3}^3 {x_6}}{3}+84 {x_2}^3
   {x_3} {x_6}+8 {x_2} {x_3} {x_4} {x_6}+2 {x_2}^2 {x_5}
   {x_6}+\frac{{x_4} {x_5} {x_6}}{7}+\frac{8}{105}z^5.\\
   \end{array}
$$
}
Here $z$ is an algebraic function of $x_1,x_2,x_3,x_4,x_5,x_6$ defined by
the equation
\begin{equation}
\label{equation:def-ans-z-e7e}
v-{x_6}+z^2=0,
\end{equation}
where
\begin{equation}
\label{equation:def-ans-z-e7v}
\begin{array}{lll}
v=-36 {x_1}^6-36 {x_1}^4 {x_2}+12 {x_1}^3 {x_3}-231 {x_1}^2
   {x_2}^2-6 {x_1}^2 {x_4}\\
\hspace{10mm}   -84 {x_1} {x_2} {x_3}
  -{x_1}
   {x_5}-49 {x_2}^3-7 {x_2} {x_4}-7 {x_3}^2.\\
   \end{array}
 \end{equation}

In this case,
$w_1=1/7,w_2=2/7,w_3=3/7,w_4=4/7,w_5=5/7,w_6=6/7,w_7=1$.

(vii) \underline{The case $E_8$}
{\footnotesize
$$
\begin{array}{ll}
&F\\
=& \frac{1}{2}{x_8}^2{x_1}  +{x_4}
   {x_5} {x_8}+{x_3} {x_6} {x_8}+{x_2} {x_7} {x_8}\\
&+\frac{472561 {x_1}^{25}}{883200}+\frac{5 {x_2} {x_1}^{22}}{6}+\frac{{x_3}
   {x_1}^{21}}{12}+\frac{14781 {x_2}^2 {x_1}^{19}}{304}+\frac{{x_4}
   {x_1}^{19}}{3}+\frac{65}{24} {x_2} {x_3} {x_1}^{18}-\frac{{x_5}
   {x_1}^{18}}{6}+\frac{6255 {x_3}^2 {x_1}^{17}}{1088}\\
   &-\frac{1225 {x_2}^3
   {x_1}^{16}}{6}-\frac{13}{2} {x_2} {x_4} {x_1}^{16}+\frac{{x_6}
   {x_1}^{16}}{3}+\frac{415}{8} {x_2}^2 {x_3} {x_1}^{15}+\frac{1}{12}
   {x_3} {x_4} {x_1}^{15}+\frac{19}{12} {x_2} {x_5}
   {x_1}^{15}-\frac{{x_7} {x_1}^{15}}{6}\\
   &+\frac{269}{8} {x_2} {x_3}^2
   {x_1}^{14}+\frac{1}{8} {x_3} {x_5} {x_1}^{14}+\frac{14879 {x_2}^4
   {x_1}^{13}}{12}+\frac{565 {x_3}^3 {x_1}^{13}}{96}+\frac{2135 {x_4}^2
   {x_1}^{13}}{104}-19 {x_2}^2 {x_4} {x_1}^{13}\\
   &-\frac{5}{2} {x_2}
   {x_6} {x_1}^{13}-\frac{766}{3} {x_2}^3 {x_3} {x_1}^{12}+\frac{329}{4}
   {x_2} {x_3} {x_4} {x_1}^{12}+22 {x_2}^2 {x_5}
   {x_1}^{12}+\frac{1}{2} {x_4} {x_5} {x_1}^{12}-\frac{1}{4} {x_3}
   {x_6} {x_1}^{12}\\
   &+\frac{5}{4} {x_2} {x_7} {x_1}^{12}+150 {x_2}^2
   {x_3}^2 {x_1}^{11}+\frac{199 {x_5}^2 {x_1}^{11}}{352}+\frac{373}{16}
   {x_3}^2 {x_4} {x_1}^{11}-\frac{43}{4} {x_2} {x_3} {x_5}
   {x_1}^{11}+\frac{1}{8} {x_3} {x_7} {x_1}^{11}\\
   &-\frac{21216 {x_2}^5
   {x_1}^{10}}{5}+\frac{1583}{48} {x_2} {x_3}^3 {x_1}^{10}-\frac{95}{2}
   {x_2} {x_4}^2 {x_1}^{10}+\frac{1727}{3} {x_2}^3 {x_4}
   {x_1}^{10}+\frac{109}{32} {x_3}^2 {x_5} {x_1}^{10}-49 {x_2}^2
   {x_6} {x_1}^{10}\\
   &-{x_4} {x_6} {x_1}^{10}+\frac{1559 {x_3}^4
   {x_1}^9}{192}+\frac{59}{2} {x_3} {x_4}^2 {x_1}^9+\frac{10795}{6}
   {x_2}^4 {x_3} {x_1}^9-\frac{239}{2} {x_2}^2 {x_3} {x_4}
   {x_1}^9-\frac{743}{6} {x_2}^3 {x_5} {x_1}^9\\
   &+\frac{105}{2} {x_2}
   {x_4} {x_5} {x_1}^9+21 {x_2} {x_3} {x_6} {x_1}^9+\frac{1}{2}
   {x_5} {x_6} {x_1}^9+\frac{49}{2} {x_2}^2 {x_7}
   {x_1}^9+\frac{1}{2} {x_4} {x_7} {x_1}^9-233 {x_2}^3 {x_3}^2
   {x_1}^8\\
   &-\frac{7}{8} {x_2} {x_5}^2 {x_1}^8+\frac{429}{4} {x_2}
   {x_3}^2 {x_4} {x_1}^8+\frac{207}{4} {x_2}^2 {x_3} {x_5}
   {x_1}^8+\frac{21}{2} {x_3} {x_4} {x_5} {x_1}^8-\frac{7}{4}
   {x_3}^2 {x_6} {x_1}^8-\frac{21}{2} {x_2} {x_3} {x_7}
   {x_1}^8\\
   &-\frac{1}{4} {x_5} {x_7} {x_1}^8+\frac{32150 {x_2}^6
   {x_1}^7}{3}+\frac{189}{2} {x_2}^2 {x_3}^3 {x_1}^7+\frac{163 {x_4}^3
   {x_1}^7}{6}+442 {x_2}^2 {x_4}^2 {x_1}^7+\frac{7}{8} {x_3} {x_5}^2
   {x_1}^7\\
   &+\frac{67 {x_6}^2 {x_1}^7}{28}-2491 {x_2}^4 {x_4}
   {x_1}^7+\frac{637}{24} {x_3}^3 {x_4} {x_1}^7+\frac{103}{8} {x_2}
   {x_3}^2 {x_5} {x_1}^7+\frac{584}{3} {x_2}^3 {x_6} {x_1}^7-26
   {x_2} {x_4} {x_6} {x_1}^7\\
   &+\frac{7}{8} {x_3}^2 {x_7}
   {x_1}^7+\frac{475}{24} {x_2} {x_3}^4 {x_1}^6+\frac{307}{2} {x_2}
   {x_3} {x_4}^2 {x_1}^6-5607 {x_2}^5 {x_3} {x_1}^6+\frac{3446}{3}
   {x_2}^3 {x_3} {x_4} {x_1}^6\\
   &+\frac{2659}{6} {x_2}^4 {x_5}
   {x_1}^6+\frac{131}{48} {x_3}^3 {x_5} {x_1}^6-\frac{1}{4} {x_4}^2
   {x_5} {x_1}^6-96 {x_2}^2 {x_4} {x_5} {x_1}^6-48 {x_2}^2
   {x_3} {x_6} {x_1}^6+\frac{39}{2} {x_3} {x_4} {x_6}
   {x_1}^6\\
   &+\frac{9}{2} {x_2} {x_5} {x_6} {x_1}^6-\frac{292}{3}
   {x_2}^3 {x_7} {x_1}^6+13 {x_2} {x_4} {x_7}
   {x_1}^6+\frac{1}{2} {x_6} {x_7} {x_1}^6+\frac{419 {x_3}^5
   {x_1}^5}{160}+\frac{2737}{2} {x_2}^4 {x_3}^2 {x_1}^5\\
   &+\frac{195}{4}
   {x_3}^2 {x_4}^2 {x_1}^5+\frac{51}{4} {x_2}^2 {x_5}^2
   {x_1}^5+\frac{19}{8} {x_4} {x_5}^2 {x_1}^5+\frac{17 {x_7}^2
   {x_1}^5}{80}+9 {x_2}^2 {x_3}^2 {x_4} {x_1}^5-147 {x_2}^3
   {x_3} {x_5} {x_1}^5\\
   &+33 {x_2} {x_3} {x_4} {x_5}
   {x_1}^5+\frac{9}{2} {x_2} {x_3}^2 {x_6} {x_1}^5+\frac{13}{4}
   {x_3} {x_5} {x_6} {x_1}^5+24 {x_2}^2 {x_3} {x_7}
   {x_1}^5-3 {x_3} {x_4} {x_7} {x_1}^5\\
   &-\frac{26144 {x_2}^7
   {x_1}^4}{3}-\frac{400}{3} {x_2}^3 {x_3}^3 {x_1}^4-\frac{64}{3} {x_2}
   {x_4}^3 {x_1}^4-\frac{{x_5}^3 {x_1}^4}{48}-640 {x_2}^3 {x_4}^2
   {x_1}^4+\frac{13}{8} {x_2} {x_3} {x_5}^2 {x_1}^4\\
   &+8 {x_2}
   {x_6}^2 {x_1}^4+4956 {x_2}^5 {x_4} {x_1}^4+65 {x_2} {x_3}^3
   {x_4} {x_1}^4+30 {x_2}^2 {x_3}^2 {x_5} {x_1}^4+\frac{51}{4}
   {x_3}^2 {x_4} {x_5} {x_1}^4\\
   &-470 {x_2}^4 {x_6}
   {x_1}^4+\frac{29}{12} {x_3}^3 {x_6} {x_1}^4+14 {x_4}^2 {x_6}
   {x_1}^4+122 {x_2}^2 {x_4} {x_6} {x_1}^4+\frac{5}{8} {x_3}
   {x_5} {x_7} {x_1}^4+\frac{185}{8} {x_2}^2 {x_3}^4
   {x_1}^3\\
   &+\frac{86}{3} {x_3} {x_4}^3 {x_1}^3+190 {x_2}^2 {x_3}
   {x_4}^2 {x_1}^3+{x_3}^2 {x_5}^2 {x_1}^3+\frac{7}{2} {x_3}
   {x_6}^2 {x_1}^3+\frac{10672}{3} {x_2}^6 {x_3} {x_1}^3+\frac{383}{48}
   {x_3}^4 {x_4} {x_1}^3\\
   &-1900 {x_2}^4 {x_3} {x_4} {x_1}^3-452
   {x_2}^5 {x_5} {x_1}^3+\frac{41}{12} {x_2} {x_3}^3 {x_5}
   {x_1}^3+34 {x_2} {x_4}^2 {x_5} {x_1}^3+198 {x_2}^3 {x_4}
   {x_5} {x_1}^3\\
   &+162 {x_2}^3 {x_3} {x_6} {x_1}^3+18 {x_2}
   {x_3} {x_4} {x_6} {x_1}^3-15 {x_2}^2 {x_5} {x_6}
   {x_1}^3+4 {x_4} {x_5} {x_6} {x_1}^3+220 {x_2}^4 {x_7}
   {x_1}^3\\
   &-\frac{1}{12} {x_3}^3 {x_7} {x_1}^3+\frac{13}{2} {x_4}^2
   {x_7} {x_1}^3-43 {x_2}^2 {x_4} {x_7} {x_1}^3+{x_2} {x_6}
   {x_7} {x_1}^3+\frac{25}{16} {x_2} {x_3}^5 {x_1}^2-492 {x_2}^5
   {x_3}^2 {x_1}^2\\
   &+36 {x_2} {x_3}^2 {x_4}^2 {x_1}^2-9 {x_2}^3
   {x_5}^2 {x_1}^2-{x_2} {x_7}^2 {x_1}^2+302 {x_2}^3 {x_3}^2
   {x_4} {x_1}^2+\frac{25}{32} {x_3}^4 {x_5} {x_1}^2+\frac{27}{2}
   {x_3} {x_4}^2 {x_5} {x_1}^2\\
   &+124 {x_2}^4 {x_3} {x_5}
   {x_1}^2-3 {x_2}^2 {x_3} {x_4} {x_5} {x_1}^2-15 {x_2}^2
   {x_3}^2 {x_6} {x_1}^2+\frac{1}{2} {x_5}^2 {x_6}
   {x_1}^2+\frac{21}{2} {x_3}^2 {x_4} {x_6} {x_1}^2\\
   &+9 {x_2}
   {x_3} {x_5} {x_6} {x_1}^2-69 {x_2}^3 {x_3} {x_7}
   {x_1}^2
   +9 {x_2} {x_3} {x_4} {x_7} {x_1}^2+\frac{15}{2}
   {x_2}^2 {x_5} {x_7} {x_1}^2-\frac{1}{2} {x_4} {x_5} {x_7}
   {x_1}^2\\
   &+{x_3} {x_6} {x_7} {x_1}^2+1240 {x_2}^8 {x_1}+\frac{11
   {x_3}^6 {x_1}}{96}+\frac{35 {x_4}^4 {x_1}}{3}+\frac{76}{3} {x_2}^4
   {x_3}^3 {x_1}+20 {x_2}^2 {x_4}^3 {x_1}+\frac{1}{6} {x_2}
   {x_5}^3 {x_1}\\
   &+360 {x_2}^4 {x_4}^2 {x_1}+8 {x_3}^3 {x_4}^2
   {x_1}+\frac{3}{2} {x_2}^2 {x_3} {x_5}^2 {x_1}+{x_3} {x_4}
   {x_5}^2 {x_1}+10 {x_2}^2 {x_6}^2 {x_1}+5 {x_4} {x_6}^2
   {x_1}\\
   &+\frac{1}{8} {x_3} {x_7}^2 {x_1}
   -\frac{3200}{3} {x_2}^6 {x_4} {x_1}-5 {x_2}^2 {x_3}^3
   {x_4} {x_1}-7 {x_2}^3 {x_3}^2 {x_5} {x_1}+\frac{21}{2} {x_2}
   {x_3}^2 {x_4} {x_5} {x_1}
   +96 {x_2}^5 {x_6} {x_1}\\
   &+\frac{5}{2}
   {x_2} {x_3}^3 {x_6} {x_1}+20 {x_2} {x_4}^2 {x_6} {x_1}-80
   {x_2}^3 {x_4} {x_6} {x_1}+\frac{5}{4} {x_3}^2 {x_5} {x_6}
   {x_1}+\frac{15}{2} {x_2}^2 {x_3}^2 {x_7} {x_1}
   +\frac{1}{8}
   {x_5}^2 {x_7} {x_1}\\
   &+\frac{3}{4} {x_3}^2 {x_4} {x_7}
   {x_1}-\frac{3}{2} {x_2} {x_3} {x_5} {x_7}
   {x_1}+\frac{{x_2}^3 {x_3}^4}{3}
   +\frac{44}{3} {x_2}
   {x_3} {x_4}^3+\frac{{x_3} {x_5}^3}{24}-44 {x_2}^3 {x_3}
   {x_4}^2+\frac{1}{4} {x_2} {x_3}^2 {x_5}^2\\
   &+{x_2} {x_3}
   {x_6}^2+\frac{{x_5} {x_6}^2}{2}-\frac{704 {x_2}^7
   {x_3}}{3}+\frac{7}{4} {x_2} {x_3}^4 {x_4}
   +176 {x_2}^5 {x_3}
   {x_4}+\frac{160 {x_2}^6 {x_5}}{3}+\frac{1}{2} {x_2}^2 {x_3}^3
   {x_5}+\frac{2 {x_4}^3 {x_5}}{3}\\
   &+2 {x_2}^2 {x_4}^2 {x_5}-32
   {x_2}^4 {x_4} {x_5}+\frac{3}{4} {x_3}^3 {x_4}
   {x_5}+\frac{{x_3}^4 {x_6}}{8}
   +6 {x_3} {x_4}^2 {x_6}-8 {x_2}^4
   {x_3} {x_6}+12 {x_2}^2 {x_3} {x_4} {x_6}\\
   &+4 {x_2}^3 {x_5}
   {x_6}+2 {x_2} {x_4} {x_5} {x_6}-\frac{168 {x_2}^5
   {x_7}}{5}-\frac{1}{4} {x_2} {x_3}^3 {x_7}
   -4 {x_2} {x_4}^2
   {x_7}+24 {x_2}^3 {x_4} {x_7}+\frac{1}{8} {x_3}^2 {x_5}
   {x_7}\\
   &-4 {x_2}^2 {x_6} {x_7}+{x_4} {x_6} {x_7}+\frac{1}{15}z^5.\\

   \end{array}
$$
}
Here $z$ is an algebraic function of $x_1,x_2,x_3,x_4,x_5,x_6,x_7$ defined by the equation
\begin{equation}
\label{equation:def-ans-z-e8e}
v-{x_7}+z^2=0,
\end{equation}
where
\begin{equation}
\label{equation:def-ans-z-e8v}
\begin{array}{lll}
v=-{x_1}^{10}+5 {x_1}^7 {x_2}+\frac{{x_1}^6 {x_3}}{2}-78 {x_1}^4
   {x_2}^2+2 {x_1}^4 {x_4}+20 {x_1}^3 {x_2} {x_3}-{x_1}^3
   {x_5}-\frac{3 {x_1}^2 {x_3}^2}{2}&&\\
 \hspace{20mm}  +88 {x_1} {x_2}^3-24 {x_1}
   {x_2} {x_4}+2 {x_1} {x_6}-12 {x_2}^2 {x_3}+2 {x_2}
   {x_5}+2 {x_3} {x_4}.\\
   \end{array}
 \end{equation}
In this case,
$w_1=1/12,w_2=1/4,w_3=1/3,w_4=1/2,w_5=7/12,w_6=3/4,w_7=5/6,w_8=1$.
\end{theorem}

\begin{remark}
 Y. Dinar discussed the existence of algebraic potentials corresponding to
subregular nilpotent elements of simple Lie algebras in \cite{D3}.
It is expected that the algebraic potentials shown in this section  coincide  those
corresponding to simple Lie algebras in \cite{D3} except $H_3$ and $H_4$ cases.
In particular,
 Y. Diner computed the algebraic potential corresponding to the subregular nilpotent element  of the Lie algebra of type $E_8$.
The author confirmed the coincidence of the potential of type $E_8$ obtained by him
and that given in  Theorem \ref{theorem:theorem7} (vii).
For the details, see \cite{DiSe}.
\end{remark}

{\bf Proof:}
We will only prove (v). (It is possible to prove the remaining cases
by an argument parallel to this case
in spite that for example, it is complicated to treat the case (vii).)

Let $x_1,\ldots,x_6$ be the variables with the weights
$w_1=2/9,w_2=1/3,w_3=5/9,w_4=2/3,w_5=8/9,w_6=1$,
respectively
and let $z$ be an algebraic function of $x_1,\ldots,x_5$
defined by the equation
\begin{equation}
\label{equation:def-z-e6}
z^2 + v - x_5=0,
\end{equation}
where $v=v(x_1,x_2,x_3,x_4)$ is a weighted homogeneous polynomial of weight 8/9.
Then
\begin{equation}
\label{equation:def-v-e6-ori}
   v= v_{0} {x_1} {x_4}+v_{1} {x_2} {x_3}+v_{2}
   {x_1} {x_2}^2+v_{3} {x_1}^4
\end{equation}
for some constants $v_j\>(j=0,1,2,3)$

From the assumption (\ref{equation:def-F1}) it follows that
the candidate of the potential is
\begin{equation}
F_A=\frac{{x_1} {x_6}^2}{2}+{x_6}
   ({x_2} {x_5}+{x_3} {x_4})
 +   {u_{0}}
       +{u_{1}} {z}+{u_{2}} {z}^2+{u_{3}} {z}^3+{u_{4}} {z}^4
   +u_{5} {z}^5,
\end{equation}
where $u_j=u_j(x_1,x_2,x_3,x_4)$ is a weighted homogeneous polynomial  of weight $4(5-j)/9$ ($j=0,1,\ldots,5$).
In particular $u_0$ is a constant.
We assume that $u_1=u_3=0$  (cf. (\ref{equation:equation-u})).
Moreover we may assume that
$u_0,u_2,u_4$ take forms
$$
\begin{array}{lll}
   u_{0}&=& s_{1} {x_1} {x_4}^3+s_{2} {x_2} {x_3} {x_4}^2+s_{3} {x_1} {x_2}^2
   {x_4}^2+s_{4} {x_1}^4 {x_4}^2+s_{5} {x_1}^2 {x_3}^2 {x_4}+s_{6} {x_2}^3 {x_3} {x_4}+s_{7} {x_1}^3
   {x_2} {x_3} {x_4}\\
   &&+s_{8} {x_1} {x_2}^4 {x_4}+s_{9} {x_1}^4 {x_2}^2 {x_4}+s_{10} {x_1}^7 {x_4}+s_{11}
   {x_3}^4+s_{12} {x_1} {x_2} {x_3}^3+s_{13} {x_1}^2 {x_2}^2 {x_3}^2+s_{14} {x_1}^5 {x_3}^2\\
   &&+s_{15} {x_2}^5
   {x_3}+s_{16} {x_1}^3 {x_2}^3 {x_3}+s_{17} {x_1}^6 {x_2} {x_3}+s_{18} {x_1} {x_2}^6+s_{19} {x_1}^4
   {x_2}^4+s_{20} {x_1}^7 {x_2}^2+s_{21} {x_1}^{10},\\
 u_{2}&=& {r_1} {x_4}^2+{r_2} {x_2}^2 {x_4}+{r_3} {x_1}^3
   {x_4}+{r_4} {x_1} {x_3}^2+{r_5} {x_1}^2 {x_2} {x_3}+{r_6} {x_2}^4+{r_7} {x_1}^3 {x_2}^2+{r_8}
   {x_1}^6,\\
 u_{4}&=& {u_{40}} {x_1}^2,\\
 \end{array}
   $$
 where
 $s_i\>(i=1,\ldots,21),\>r_j\>(j=1,\ldots, 8)$ and $u_{40}$ are constants to be determined.
 
We assume  $u_5=0$. In virtue of 
 (\ref{equation:def-z-e6}), we find that $F_A$ turns out to be a polynomial of $x_1,\ldots,x_6$.
 This doesn't match our purpose.
Accordingly, we may assume 

(A.1) $u_5\not=0$.

We  compute  $P,\>C,\>{\tilde B}_j,\>T$ introduced in \S\ref{section:sec1}.
We start the construction of the algebraic potential with
studying a basic property of $\det(T)$.
In virtue of
(\ref{equation:def-z-e6}), $\det(T)$ is regarded as a polynomial of $x_1,\ldots,x_4,x_6,z$.
It is straightforward to show that
$$
\det(T)|_{x_1=x_2=x_3=x_4=z=0}=x_6^2(x_6^2-cx_5^3)^2,
$$
where
\begin{equation}
\label{equation:def-c-e6}
c=\frac{8}{9}(3s_1+2r_1s_2-3r_1v_0-2r_1^2v_1).
\end{equation}

We assume that

(A.2) $c\not=0.$

Under the assumptions (A.1), (A.2), we show that the constants $s_i,r_j,u_{40},v_0,v_1,v_2,v_3$ are uniquely determined
if $F$ is a solution of the WDVV equation.

Let
$B_{pq}=z[{\tilde B}_p,\>{\tilde B}_q]$  and let $B_{pq}(i,j)$ be 
the $(i,j)$ entry of $B_{pq}$.
Since each matrix entry of $C$ and   $z\partial_{x_j}z$  are polynomials of $x_1,\ldots,x_4,z$,
so is $B_{pq}(i,j)$.

Our purpose is to find the conditions among the constants
$s_i,r_j,u_{40},v_0,v_1,v_2,v_3$ so that $B_{pq}=O$ for all $p,q=1,\ldots,6$.
Since ${\tilde B}_6$ is the identity matrix, it is sufficient to show that $B_{pq}=O$ for $p,q<6$.

We focus our attention to $B_{45}$.

It is easy to show that $B_{45}(6,j)=0$ for $j=1,\ldots,6$.

We compute $B_{45}(5,j)$.
It is straightforward to show that $B_{45}(5,j)=0$ for $j=1,2$.
From the identity equation $B_{45}(5,3)=0$ for $x_1,x_2,x_3,x_4,z$, we obtain the two equations
$$
s_2 - 2 r_1 v_1=0,
\quad
2 s_3 - s_2v_0- 3 s_1v_1 + 2 r_1 v_0 v_1 - 4 r_1 v_2=0,
$$
which imply
\begin{equation}
s_2=2r_1v_1,
\quad
s_3=\frac{1}{2}( s_2v_0+ 3 s_1v_1 -2 r_1 v_0 v_1 + 4 r_1 v_2).
\end{equation}
In this manner, $s_2,s_3$ are eliminated in the matrices $B_{pq}$.
From the identity equation $B_{45}(5,4)=0$ for $x_1,x_2,x_3,x_4,z$, we obtain four equations which imply
\begin{equation}
v_0=2r_1v_1,\quad
s_6=\frac{1}{3}(r_2v_1+4r_1{v_1}^2+v_2),\quad
u_{40}=r_1r_4,\quad
s_7=r_3v_1+4r_1s_5v_1+4v_3.
\end{equation}
As a consequence,
$v_0,\>s_6,\>u_{40},\>s_7$ are eliminated in the matrices $B_{pq}$.
From the identity equation $B_{45}(5,5)=0$ for $x_1,x_2,x_3,x_4,z_0$, we obtain six equations  which imply

\begin{equation}
\left\{
\begin{array}{lll}
r_2&=&-4r_1v_1,\\
v_2&=&-6r_1{v_1}^2,\\
s_8&=&\frac{1}{6}{v_1}^2(3s_1+104{r_1}^2v_1),\\
r_3&=&\frac{2}{3}(-r_1r_5+8{r_1}^2r_4v_1),\\
s_9&=&-4r_1v_1(-r_1r_5v_1-r_1s_5v_1+8{r_1}^2r_4{v_1}^2-2v_3)\\
\end{array}
\right.
\end{equation}
and the equation
\begin{equation}
\label{equation:equation-v-s1}
v_1(3s_1-16{r_1}^2v_1)=0.
\end{equation}
At this stage, it is necessary to treat the cases $v_1=0$ and $3s_1-16{r_1}^2v_1=0$, separately.
To avoid to do, we check the remaining matrix entries of $B_{45}$.
In spite that there are no terms of $B_{45}(5,6),\>B_{45}(4,j)\>(j=1,2,3,4)$
whose coefficients are 
related to (\ref{equation:equation-v-s1}),
we find the equation
 $$3s_1-16{r_1}^2v_1=0
 $$
by considering $B_{45}(4,5)$.
Then
 (\ref{equation:equation-v-s1}) is reduced to
\begin{equation}
s_1=\frac{16}{3}{r_1}^2v_1.
\end{equation}
Moreover, by substituting $s_1,s_2,v_0$ in the definition of c
(cf. (\ref{equation:def-c-e6})), we find that
$$
c=\frac{32}{3}{r_1}^2v_1.
$$
Then the assumption (A.2) implies that

(A.2)' $r_1v_1\not=0.$

Under the assumption (A.1), (A.2)', the remaining constants $r_i,\>s_j,\>v_k$ are uniquely expressed in terms of $u_5,r_1,v_1$, which we are going to explain.
We have already shown that
$$s_2,s_3,v_0,s_6,u_{40},s_7,
r_2,v_2,s_8,r_3,s_9,s_1
$$
are expressed by the  constants
 $$s_4,s_5,s_i(i=10,\ldots,21),r_1,r_j(j=4,\ldots,8),v_1,v_3,u_5.
 $$
The identity equation $B_{45}(5,6)=0$
implies that $r_4,r_5,r_6,v_3,s_4$ are expressed by the  constants
 $$s_5,s_i(i=10,\ldots,21),r_1,r_7,r_8,v_1,u_5.
 $$
Then by the following order,  we find that $s_5,r_7,r_8,s_{10}$
are expressed by the  constants
 $s_i(i=11,\ldots,21),r_1,v_1,u_5$.

$
B_{45}(4,4)\Longrightarrow s_5,\quad
B_{45}(4,6)\Longrightarrow r_7,\quad
B_{45}(3,6)\Longrightarrow r_8,\quad
B_{45}(2,6)\Longrightarrow s_{10}.
$

As a consequence, it follows that $B_{45}=O$.

We next treat $B_{35}$.
By the following order,  we find that $s_{12},s_{13},s_{15},s_{16},s_{11},s_{17},s_{14}$
are expressed by the  constants
$s_{18},s_{19},s_{20},s_{21},r_1,v_1,u_5$.

$
B_{35}(5,4)\Longrightarrow s_{12},s_{13},\quad
B_{35}(5,5)\Longrightarrow s_{15},s_{16},\quad
B_{35}(5,6)\Longrightarrow s_{11},s_{17}\quad
B_{35}(4,6)\Longrightarrow s_{14}.
$

As a consequence, it follows that $B_{35}=O$.
 
 We next treat $B_{25}$.
By the following order,  we find that $s_{18},s_{19},s_{20}$
are expressed by the constants
$u_5,r_1,v_1,s_{21}$.

$
B_{25}(5,5)\Longrightarrow s_{18},s_{19},\quad
B_{25}(5,6)\Longrightarrow s_{20}.
$

 As a consequence, it follows that $B_{25}=O$.

 To determine $s_{21}$, it is sufficient to compute $B_{15}(2,6)$ and 
 we conclude that the all the constants
  $r_i,\>s_j,\>v_k,u_{40}$ are expressed by $u_5,r_1,v_1$.
Moreover it is straightforward to show that $ B_{15}=0$.
Then it is easy to confirm that $B_{pq}=O\>(p,q<5)$.
 
 The potential $F_A$ has parameters $u_5,r_1,v_1$.
By the coordinate the transformation
 $$
 x_j\to m_jx_j\>(j=1,\ldots,6),
 $$
 where
 $$
 \begin{array}{l}
 m_1=2,\>
 m_2=5r_1^{1/2}u_5v_1,\>
 m_3=10r_1^{3/2}u_5v_1,\\
 m_4=\frac{25}{2}r_1u_5^2v_1^3,\>
 m_5=5^{1/2}r_1^{1/2}u_5^{1/2}v_1,\>
 m_6=\frac{125}{2}r_1^{5/2}u_5^3v_1^5,\\
 \end{array}
$$
  $F_A$ is transformed to
  {\footnotesize
 $$
 \begin{array}{ll}
 &F_B\\
 =&
 \frac{1}{2}{x_1} {x_6}^2+{x_2} {x_5} {x_6}+{x_3} {x_4} {x_6}\\
&   +
 \frac{16929}{20} {x_1}^{10}+\frac{8100}{7} {x_1}^7 {x_2}^2-162 {x_1}^7
   {x_4}+162 {x_1}^6 {x_2} {x_3}+\frac{891}{10} {x_1}^5 {x_3}^2-462 {x_1}^4
 {x_2}^4-180 {x_1}^4 {x_2}^2 {x_4}\\
&   +\frac{99}{2} {x_1}^4 {x_4}^2+312
   {x_1}^3 {x_2}^3 {x_3}+108 {x_1}^3 {x_2} {x_3} {x_4}-72 {x_1}^2
   {x_2}^2 {x_3}^2+27 {x_1}^2 {x_3}^2 {x_4}-\frac{112}{3} {x_1} {x_2}^6+40
   {x_1} {x_2}^4 {x_4}\\
&   -4 {x_1} {x_2}^2 {x_4}^2+15 {x_1} {x_2}
   {x_3}^3+\frac{8}{3} {x_1} {x_4}^3+\frac{16}{5} {x_2}^5
   {x_3}-4 {x_2}^3 {x_3} {x_4}+2 {x_2} {x_3} {x_4}^2+\frac{3}{8} {x_3}^4\\
& +    (-\frac{135}{2} {x_1}^6+36
   {x_1}^3 {x_2}^2+18 {x_1}^3 {x_4}-18 {x_1}^2 {x_2} {x_3}+\frac{9}{2} {x_1}
   {x_3}^2+6 {x_2}^4-4 {x_2}^2 {x_4}+\frac{1}{2}{x_4}^2)z^2\\
&   +\frac{9}{4}
   {x_1}^2z^4+\frac{2}{5}z^5\\
   \end{array}
   $$
   }
 and
 (\ref{equation:def-v-e6-ori}) turns out to be
 (\ref{equation:def-ans-z-e6v}).
 
 It is straightforward to show that $F_B$ coincides with $F$ defined in (v) under the condition
  (\ref{equation:def-ans-z-e6e}).

From now on, we show that $\det(T)$ is not only  identically zero but also
has no multiple factor as a polynomial of $x_6$.

It is straightforward to show that
if $T_0$ is the restriction of the matrix $T$ for the subspace defined by $x_1=x_2=x_3=0$,
then
$$T_0=\left(
\begin{array}{cccccc}
 {x_6} & -\frac{10 {x_4} {z}}{3} & \frac{4}{3} \left(5 {x_4}^2-2 {z}^3\right) & 0 & 0 & \frac{8}{3} {z} \left(3 {z}^3+4
   {x_4}^2\right) \\
 \frac{8 {z}^2}{9} & {x_6} & 0 & \frac{4}{3} \left({x_4}^2-2 {z}^3\right) & -\frac{112 {x_4} {z}^2}{9} & 0 \\
 \frac{2 {x_4}}{3} & 0 & {x_6} & 0 & \frac{4}{3} \left({x_4}^2-2 {z}^3\right) & 0 \\
 0 & \frac{2 {x_4}}{3} & \frac{8 {z}^2}{9} & {x_6} & 0 & \frac{4}{3} \left(5 {x_4}^2-2 {z}^3\right) \\
 0 & \frac{2 {z}}{3} & \frac{2 {x_4}}{3} & 0 & {x_6} & -\frac{10 {x_4} {z}}{3} \\
 0 & 0 & 0 & \frac{2 {x_4}}{3} & \frac{8 {z}^2}{9} & {x_6} \\
\end{array}
\right).
$$
Accordingly, we find that
{\footnotesize
$$
\begin{array}{ll}
&\det(T_0)\\
=&\frac{1}{729} \left(16 {x_4}^3-3 {x_6}^2\right) \left(20736 {x_4}^4 {z}^3+1296 {x_4}^3 {x_6}^2-18432 {x_4}^2 {z}^6-5184
   {x_4} {x_6}^2 {z}^3-243 {x_6}^4+4096 {z}^9\right)\\
   \end{array}
$$
}
It is provable that $\det(T_0)$ has neither multiple factor as a polynomial of $x_6$,
which implies that $\det(T)$ has nor  multiple factor
and (v)  follows.

\subsection{The data by  T. Douvropoulos}
\label{subsection:sec3-2}

In the previous subsection, we showed algebraic potentials of a special kind
in which the defining equation of $z_0$ is quadratic.
As in the case of $F_4$ and $H_4$ treated in \S\ref{subsection:sec1-4},
there are algebraic potentials such that the degree of the defining equation of $z$ is 
three or grater than three.
It is difficult to construct the algebraic potentials in such cases.
To treat these cases, the data calculated by
T. Douvropoulos is very helpful.
He determines primitive conjugacy classes, eigenvalues of representatives,
the extension degrees of the corresponding algebraic potentials
for the cases
$H_3,D_4,F_4,H_4,E_6,E_7,E_8$.
The data are given below:

\vspace{5mm}
$H_3$(0) (Coxeter): [2,6,10], degree: 1 (polynomial prepotential)

$H_3$(1): [2,4,6], degree: 2

$H_3$(2): [6,8,10], degree: 4

\vspace{5mm}
$D_4$(0) (Coxeter): [2,4,4,6], degree: 1 (polynomial prepotential)

$D_4$(1): [2,2,4,4], degree: 2

\vspace{5mm}
$F_4$(0) (Coxeter): [2,6,8,12], degree: 1 (polynomial prepotential)

$F_4$(1): [2,2,6,6], degree: 3

\vspace{5mm}
$H_4$(0) (Coxeter): [2,12,20,30], degree: 1 (polynomial prepotential)

$H_4$(1): [2,10,12,20], degree: 2

$H_4$(2): [2,5,12,15], degree: 3

$H_4$(3): [2,6,8,12], degree: 4

$H_4$(4): [2,2,10,10], degree: 5

$H_4$(5): This is a quasi-Coxeter element of order 10, but not regular.

$H_4$(6): [4,9,10,15], degree: 9

$H_4$(7): [6,10,16,20], degree: 8

$H_4$(8): [2,2,6,6], degree: 15

$H_4$(9): [14,20,24,30], degree: 14

$H_4$(10): [6,6,10,10], degree: 45

\vspace{5mm}

$E_6$(0) (Coxeter): [2,5,6,8,9,12], degree: 1 (polynomial prepotential)

$E_6$(1): [2,3,5,6,8,9], degree: 2

$E_6$(2): [2,2,3,5,6,6], degree: 5

\vspace{5mm}

$E_7$(0) (Coxeter): [2,6,8,10,12,14,18], degree: 1 (polynomial prepotential)

$E_7$(1): [2,4,6,8,10,12,14], degree: 2

$E_7$(2): This is a quasi-Coxeter element of order 12 that is not regular.

$E_7$(3): This is a quasi-Coxeter element of order 10 that is not regular.

$E_7$(4): [2,2,2,4,6,6,6], degree: 18

\vspace{5mm}

$E_8$(0) (Coxeter): [2,8,12,14,18,20,24,30], degree: 1 (polynomial prepotential)

$E_8$(1): [2,6,8,12,14,18,20,24], degree: 2

$E_8$(2): [2,4,8,10,12,14,18,20], degree: 3

$E_8$(3): [2,2,6,6,8,8,12,12], degree: 8 (this might be related to F4)

$E_8$(4): This is a quasi-Coxeter element of order 18 that is not regular.

$E_8$(5): [2,3,5,8,9,12,14,15], degree: 7

$E_8$(6): [2,2,4,4,8,8,10,10], degree: 20

$E_8$(7): This is a quasi-Coxeter element of order 12 that is not regular.

$E_8$(8): [2,2,2,2,6,6,6,6]: degree: 135 (=5*3*3*3)

\vspace{5mm}

\subsection{Examples  (II)}
\label{subsection:sec3-3}
Based on the data given in \S\ref{subsection:sec3-2}, we apply
the idea explained in \S\ref{subsection:sec3-1} to the primitive conjugacy classes
of reflection groups.
Then we succeeded to construct algebraic potentials in some cases.
We report here the algebraic potentials.

\subsubsection{Algebraic potentials related to the reflection group $H_4$}
\label{subsubsection:sec3-3-1}

\underline{$H_4(0)$ case}

This case corresponds to the polynomial potential of the reflection group of type $H_4$.
$$
\begin{array}{lll}
F_{H_4(0)}&=&\frac{{x_1} {x_4}^2}{2}+{x_2} {x_3} {x_4}+\frac{32 {x_1}^{31}}{22395255890625}+\frac{2 {x_1}^{19}
   {x_2}^2}{15582375}+\frac{{x_1}^{13} {x_2}^3}{72900}+\frac{8 {x_1}^{11}
   {x_3}^2}{22275}\\
   &&+\frac{1}{810} {x_1}^9 {x_2}^2 {x_3}+\frac{{x_1}^7
   {x_2}^4}{1080}
   +\frac{1}{15} {x_1}^5 {x_2} {x_3}^2+\frac{1}{18} {x_1}^3
   {x_2}^3 {x_3}+\frac{{x_1} {x_2}^5}{240}+\frac{2 {x_1}
   {x_3}^3}{3}.\\
   \end{array}
   $$
(Cf. \cite{Du}, Lecture 4, Example 4.4.)

\underline{$H_4$(1) case}

See Theorem \ref{theorem:theorem7}, 
 (iv).

\underline{$H_4$(2) case}
{\footnotesize
$$
\begin{array}{lll}
F_{H_4(2)}&=&\displaystyle{\frac{1}{2} {x_1} {x_4}^2
   +{x_2} {x_3}
   {x_4}}\\
     &&\displaystyle{+\frac{1}{73920}(-66084040 {x_1}^{16}+143564400
   {x_1}^{11} {x_2}^2-40727610 {x_1}^6 {x_2}^4-392931 {x_1}
   {x_2}^6)}\\[3mm]
&&\displaystyle{
   -\frac{3}{4} {x_1}^4 \left(2288 {x_1}^{10}-1620 {x_1}^5
   {x_2}^2-27 {x_2}^4\right)z
-760 {x_1}^{12} z^2}\\[3mm]
   &&\displaystyle{+\frac{1}{48}  \left(1744 {x_1}^{10}-4860
   {x_1}^5 {x_2}^2-81 {x_2}^4\right)z^3
+140 {x_1}^8 z^4+24 {x_1}^6 z^5-\frac{53}{6} {x_1}^4
   z^6-\frac{10}{7} {x_1}^2 z^7
   +\frac{1}{4}z^8.}\\
   \end{array}
$$
}
Here $z$ is an algebraic function of
$x_1,x_2,x_3$ defined by the equation
$$
-11 {x_1}^6-12 {x_1}^4 z+\frac{27 {x_1} {x_2}^2}{4}-{x_3}+z^3=0.
$$

\underline{$H_4$(3) case}
{\footnotesize
 $$
 \begin{array}{ll}  
F_{H_4(3)}=&  \frac{1}{2} {x_1}{x_4}^2+x_2x_3x_4\\
& +  \frac{125}{4004} {x_1} \left(441124488 {x_1}^{12}+123891768 {x_1}^9 {x_2}+12602304 {x_1}^6
   {x_2}^2-496496 {x_1}^3 {x_2}^3-9009 {x_2}^4\right)\\
   &
   -\frac{100}{3} {x_1}^3 z \left(14
   {x_1}^3-{x_2}\right) \left(5453 {x_1}^6+10266 {x_1}^3 {x_2}-660
   {x_2}^2\right)\\
   &+360 {x_1}^2 z^2 \left(119 {x_1}^3-5 {x_2}\right) \left(14
   {x_1}^3-{x_2}\right)^2-400 {x_1} z^3 \left(14
   {x_1}^3-{x_2}\right)^3\\
   &-\frac{10}{3} z^4 \left(11011 {x_1}^9-13794 {x_1}^6
   {x_2}+912 {x_1}^3 {x_2}^2-6 {x_2}^3\right)\\
   &-72 {x_1}^2 z^5 \left(119
   {x_1}^3-5 {x_2}\right) \left(14 {x_1}^3-{x_2}\right)+120 {x_1} z^6 \left(14
   {x_1}^3-{x_2}\right)^2\\
   &+\frac{272}{35} z^7 \left(14
   {x_1}^3-{x_2}\right)^2+\frac{18}{5} {x_1}^2 z^8 \left(119 {x_1}^3-5
   {x_2}\right)
   -12 {x_1} z^9\left(14 {x_1}^3-{x_2}\right)\\
   &-\frac{16}{15} z^{10}
   \left(14 {x_1}^3-{x_2}\right)+\frac{2 {x_1} z^{12}}{5}+\frac{256 z^{13}}{4875}.\\
   \end{array}
   $$
   }
Here $z$ is an algebraic function of
$x_1,x_2,x_3$ defined by the equation

$$
75 {x_1}^4-140 {x_1}^3 z+25 {x_1} {x_2}+10 {x_2} z-{x_3}+z^4=0.
$$

\underline{$H_4$(4) case}

See \S\ref{section:sec5}.

\underline{$H_4$(5) case} Not regular.

\underline{$H_4$(6) case}
{\footnotesize
$$
\begin{array}{lll}
F_{H_4(6)}&=&   \frac{1}{2}{x_1} {x_4}^2
   +{x_2} {x_3} {x_4}
   -\frac{1875 {x_1}^8
   {z}}{16}+\frac{1825 {x_1}^7 {z}^3}{8}-\frac{625
   {x_1}^6 {x_3}}{16}-\frac{475 {x_1}^6
   {z}^5}{8}+\frac{395}{8} {x_1}^5 {x_3}
   {z}^2-\frac{895 {x_1}^5 {z}^7}{8}\\
   &&-\frac{415}{16}
   {x_1}^4 {x_3} {z}^4
   +\frac{655 {x_1}^4
   {z}^9}{12}-2 {x_1}^3 {x_3}^2 {z}-\frac{59}{4}
   {x_1}^3 {x_3} {z}^6+\frac{259 {x_1}^3
   {z}^{11}}{8}+2 {x_1}^2 {x_3}^2
   {z}^3+\frac{361}{16} {x_1}^2 {x_3}
   {z}^8\\
   &&-\frac{1617 {x_1}^2 {z}^{13}}{104}-\frac{2
   {x_1} {x_3}^3}{15}+\frac{2}{5} {x_1} {x_3}^2
   {z}^5+\frac{159}{40} {x_1} {x_3}
   {z}^{10}-\frac{153 {x_1}
   {z}^{15}}{40}+\frac{4 {x_3}^3
   {z}^2}{5}
   +\frac{66 {x_3}^2 {z}^7}{35}-\frac{333
   {x_3} {z}^{12}}{80}\\
   &&+\frac{2187
   {z}^{17}}{1360}.\\
   \end{array}
   $$
}
Here $z$ is an algebraic function of
$x_1,x_2,x_3$ defined by the equation
   $$
{x_2}^2-z(x_1-{z}^2)^2(5 {x_1}^2- {z}^4)-
x_3({x_1}-z^2)^2=0.
$$ 

\underline{$H_4$(7) case}
{\footnotesize
$$
\begin{array}{lll}
F_{H_4(7)}&=&\frac{1}{2}{x_1}{x_4}^2+{x_2} {x_3} {x_4}

   -\frac{4096}{135} {x_1} {x_2} \left(32 {x_1}^5+315 {x_2}^3\right)
    +\frac{32768}{1125} {x_1}^2 \left(2 {x_1}^5+75 {x_2}^3\right)z^2
    \\&&
   -\frac{32768}{225} {x_2}  \left(2 {x_1}^5+75 {x_2}^3\right)z^3
    -\frac{16384}{5625} {x_1} \left(14 {x_1}^5+375
   {x_2}^3\right)z^5
  +\frac{34816}{225} {x_1}^4 {x_2}z^6
    -\frac{116736}{175} {x_1}^2 {x_2}^2z^7
    \\&&
+\frac{256}{75}  \left(3 {x_1}^5+220 {x_2}^3\right)z^8
 -\frac{118784}{945} {x_1}^3 {x_2} z^9
   +\frac{44544}{175} {x_1} {x_2}^2 z^{10}
    -\frac{5632}{1575} {x_1}^4 z^{11}
     +\frac{832}{225} {x_1}^2 {x_2} z^{12}
\\&&
 +\frac{17664 {x_2}^2 z^{13}}{455}
  -\frac{352}{315} {x_1}^3z^{14}
     +\frac{1568}{225} {x_1} {x_2}z^{15}
   +\frac{496 {x_1}^2 z^{17}}{2975}
 +\frac{496 {x_2} z^{18}}{945}
    +\frac{71 {x_1}z^{20}}{1575}
  +\frac{16z^{23}}{7245}.\\
   \end{array}
   $$
   }
Here $z$ is an algebraic function of
$x_1,x_2,x_3$ defined by the equation
 $$  
   -\frac{64 {x_1}^2 z^2}{5}+128 {x_1}
   {x_2}+\frac{32 {x_1} z^5}{5}+64 {x_2}
   z^3-{x_3}+z^8=0.
   $$

\underline{$H_4$(9) case}
{\footnotesize
$$
\begin{array}{lll}
F_{H_4(9)}&=& \frac{{x_1}
   {x_4}^2}{2}+{x_2}{x_3} {x_4}
 -\frac{189 {x_1}^6}{800 z^5}-\frac{189 {x_1}^5
   z^2}{32}-\frac{63 {x_1}^4 {x_2}}{40 z}-\frac{3017 {x_1}^4
   z^9}{52}-\frac{273}{16} {x_1}^3 {x_2} z^6-\frac{53865
   {x_1}^3 z^{16}}{416}\\
  & &-\frac{63}{50} {x_1}^2 {x_2}^2
   z^3
   -\frac{1988}{65} {x_1}^2 {x_2} z^{13}-\frac{50246
   {x_1}^2 z^{23}}{345}+\frac{{x_1}
   {x_2}^3}{50}-\frac{217}{100} {x_1} {x_2}^2
   z^{10}-\frac{20321}{780} {x_1} {x_2} z^{20}\\
    &&-\frac{16324 {x_1} z^{30}}{195}-\frac{2
   {x_2}^3 z^7}{25}-\frac{308 {x_2}^2 z^{17}}{221}
    -\frac{4984 {x_2} z^{27}}{585}-\frac{9072
   z^{37}}{481}.\\
   \end{array}
   $$
   }
Here $z$ is an algebraic function of
$x_1,x_2,x_3$ defined by the relation
   $$
\frac{3}{20} {x_1}^2+\frac{1}{5}{x_2} z^4+\frac{8}{5} {x_1}
   z^7-{x_3}z^2+z^{14}=0.
   $$

\underline{$H_4$(8), $H_4$(10) cases} To be determined.

\subsubsection{$E_8(2)$ case}

By a little tedious computation, it is possible to construct the algebraic potential in this case.
We need a preparation to show the result.
We define $g_1,g_2,g_3$ by
$$
\begin{array}{lll}
{g_1}&=& 294 {x_1}^7 {x_2}-\frac{1323 {x_1}^5 {x_2}^2}{2}+6 {x_1}^5 {x_3}+\frac{405
   {x_1}^3 {x_2}^3}{2}-\frac{45}{7} {x_1}^3 {x_2} {x_3}+\frac{5}{2} {x_1}^2 {x_2} {x_4}-49
   {x_1}^2 {x_6}\\
   &&-\frac{405 {x_1} {x_2}^4}{16}+\frac{135}{98} {x_1} {x_2}^2 {x_3}-27 {x_1}
   {x_2} {x_5}-\frac{9 {x_1} {x_3}^2}{686}-\frac{3 {x_2}^2 {x_4}}{7}+\frac{21 {x_2}
   {x_6}}{2}+\frac{{x_3} {x_4}}{98}-{x_7},\\
   {g_2}&=& -196 {x_1}^6-\frac{135 {x_1}^2
   {x_2}^2}{2}-{x_1} {x_4}+\frac{135 {x_2}^3}{28}-\frac{27 {x_2} {x_3}}{98}+9 {x_5},\\
   {g_3}&=&
   -\frac{614656 {x_1}^{10} {x_2}}{9}+16464 {x_1}^8 {x_2}^2-\frac{8624 {x_1}^8 {x_3}}{9}-18914
   {x_1}^6 {x_2}^3+\frac{728}{3} {x_1}^6 {x_2} {x_3}-\frac{4900}{27} {x_1}^5 {x_2}
   {x_4}\\
   &&+\frac{96040 {x_1}^5 {x_6}}{27}+\frac{2835 {x_1}^4 {x_2}^4}{2}-196 {x_1}^4 {x_2}^2
   {x_3}+\frac{8 {x_1}^4 {x_3}^2}{21}+\frac{70}{9} {x_1}^3 {x_2}^2 {x_4}-\frac{6860}{9} {x_1}^3
   {x_2} {x_6}\\
   &&-\frac{80}{27} {x_1}^3 {x_3} {x_4}-270 {x_1}^2 {x_2}^5+\frac{81}{7} {x_1}^2
   {x_2}^3 {x_3}-\frac{12}{49} {x_1}^2 {x_2} {x_3}^2+3 {x_1} {x_2}^3 {x_4}-196 {x_1}
   {x_2}^2 {x_6}\\
   &&+\frac{2}{21} {x_1} {x_2} {x_3} {x_4}+\frac{28 {x_1} {x_3}
   {x_6}}{3}+\frac{2295 {x_2}^6}{112}-\frac{297 {x_2}^4 {x_3}}{196}+\frac{9 {x_2}^2
   {x_3}^2}{686}+\frac{{x_2} {x_4}^2}{9}+\frac{2 {x_3}^3}{2401}-\frac{98 {x_4} {x_6}}{27}.\\
   \end{array}
   $$
   Then the algebraic potential reads
   {\footnotesize
$$
\begin{array}{ll}
&F_{E_8(2)}\\
= &{x_4} {x_5} {x_8}+{x_3}
   {x_6} {x_8}+{x_2} {x_7} {x_8}+\frac{{x_8}^2 {x_1}}{2}\\
   &
+\frac{44220964928 {x_1}^{21}}{855}+\frac{2921459968 {x_2}^2 {x_1}^{17}}{51}-\frac{30118144 {x_4} {x_1}^{16}}{135}-21916328 {x_2}^3
   {x_1}^{15}+\frac{998816}{15} {x_2} {x_3} {x_1}^{15}\\
   &+\frac{30118144 {x_5} {x_1}^{15}}{15}
   +16682148 {x_2}^4 {x_1}^{13}+\frac{344960
   {x_3}^2 {x_1}^{13}}{39}-343000 {x_2}^2 {x_3} {x_1}^{13}+\frac{3476648}{9} {x_2}^2 {x_4} {x_1}^{12}\\
   &-\frac{9411920}{9} {x_2}
   {x_6} {x_1}^{12}
   -\frac{18111429 {x_2}^5 {x_1}^{11}}{5}-6272 {x_2} {x_3}^2 {x_1}^{11}+\frac{1133272 {x_4}^2
   {x_1}^{11}}{891}+264012 {x_2}^3 {x_3} {x_1}^{11}\\
   &-2996448 {x_2}^2 {x_5} {x_1}^{11}
   -\frac{70658}{3} {x_2}^3 {x_4}
   {x_1}^{10}+\frac{20776}{3} {x_2} {x_3} {x_4} {x_1}^{10}+\frac{76832}{3} {x_4} {x_5} {x_1}^{10}+\frac{7731220}{3} {x_2}^2
   {x_6} {x_1}^{10}\\
   &-\frac{192080}{9} {x_3} {x_6} {x_1}^{10}
   +\frac{4994325 {x_2}^6 {x_1}^9}{8}+\frac{136 {x_3}^3 {x_1}^9}{7}+2748
   {x_2}^2 {x_3}^2 {x_1}^9+364952 {x_5}^2 {x_1}^9-\frac{57015}{2} {x_2}^4 {x_3} {x_1}^9\\
   &-889056 {x_2}^3 {x_5}
   {x_1}^9
   -15288 {x_2} {x_3} {x_5} {x_1}^9+\frac{73941}{2} {x_2}^4 {x_4} {x_1}^8-\frac{152}{9} {x_3}^2 {x_4}
   {x_1}^8-\frac{2800}{3} {x_2}^2 {x_3} {x_4} {x_1}^8\\
   &-1166886 {x_2}^3 {x_6} {x_1}^8
   +24696 {x_2} {x_3} {x_6}
   {x_1}^8+\frac{300591 {x_2}^7 {x_1}^7}{8}-\frac{612}{49} {x_2} {x_3}^3 {x_1}^7+\frac{495}{7} {x_2}^3 {x_3}^2
   {x_1}^7\\
   &+\frac{29204}{27} {x_2}^2 {x_4}^2 {x_1}^7+\frac{7395080 {x_6}^2 {x_1}^7}{27}-8235 {x_2}^5 {x_3} {x_1}^7+52920
   {x_2}^4 {x_5} {x_1}^7-288 {x_3}^2 {x_5} {x_1}^7\\
   &-26712 {x_2}^2 {x_3} {x_5} {x_1}^7-\frac{211288}{27} {x_2} {x_4}
   {x_6} {x_1}^7-\frac{196 {x_4}^3 {x_1}^6}{729}-\frac{76251}{20} {x_2}^5 {x_4} {x_1}^6+\frac{44}{7} {x_2} {x_3}^2 {x_4}
   {x_1}^6\\
   &+393 {x_2}^3 {x_3} {x_4} {x_1}^6-11074 {x_2}^2 {x_4} {x_5} {x_1}^6+\frac{229467}{2} {x_2}^4 {x_6}
   {x_1}^6+\frac{392}{3} {x_3}^2 {x_6} {x_1}^6
   -2450 {x_2}^2 {x_3} {x_6} {x_1}^6\\
   &+240100 {x_2} {x_5} {x_6}
   {x_1}^6-\frac{98901 {x_2}^8 {x_1}^5}{32}+\frac{192 {x_3}^4 {x_1}^5}{2401}+\frac{1602}{343} {x_2}^2 {x_3}^3
   {x_1}^5
   -\frac{33615}{196} {x_2}^4 {x_3}^2 {x_1}^5-\frac{1253}{18} {x_2}^3 {x_4}^2 {x_1}^5\\
   &+\frac{164}{27} {x_2} {x_3}
   {x_4}^2 {x_1}^5+34398 {x_2}^2 {x_5}^2 {x_1}^5-\frac{840350}{9} {x_2} {x_6}^2 {x_1}^5
   +\frac{16929}{5} {x_2}^6 {x_3}
   {x_1}^5+\frac{23247}{10} {x_2}^5 {x_5} {x_1}^5\\
   &-\frac{792}{7} {x_2} {x_3}^2 {x_5} {x_1}^5-\frac{1568}{9} {x_4}^2 {x_5}
   {x_1}^5+6318 {x_2}^3 {x_3} {x_5} {x_1}^5
   -686 {x_2}^2 {x_4} {x_6} {x_1}^5-\frac{6272}{27} {x_3} {x_4} {x_6}
   {x_1}^5\\
   &+588 {x_4} {x_5}^2 {x_1}^4+\frac{6435}{8} {x_2}^6 {x_4} {x_1}^4-\frac{64 {x_3}^3 {x_4}
   {x_1}^4}{1029}-\frac{32}{49} {x_2}^2 {x_3}^2 {x_4} {x_1}^4-\frac{1095}{28} {x_2}^4 {x_3} {x_4} {x_1}^4+\frac{4137}{2}
   {x_2}^3 {x_4} {x_5} {x_1}^4\\
   &-158 {x_2} {x_3} {x_4} {x_5} {x_1}^4+\frac{81585}{4} {x_2}^5 {x_6} {x_1}^4-10
   {x_2} {x_3}^2 {x_6} {x_1}^4-\frac{4515}{2} {x_2}^3 {x_3} {x_6} {x_1}^4-5145 {x_2}^2 {x_5} {x_6} {x_1}^4\\
   &-2940
   {x_3} {x_5} {x_6} {x_1}^4+\frac{111537 {x_2}^9 {x_1}^3}{224}-\frac{18 {x_2} {x_3}^4 {x_1}^3}{2401}-\frac{675}{686}
   {x_2}^3 {x_3}^3 {x_1}^3-10584 {x_5}^3 {x_1}^3+\frac{8262}{343} {x_2}^5 {x_3}^2 {x_1}^3\\
   &+\frac{121}{6} {x_2}^4 {x_4}^2
   {x_1}^3+\frac{74 {x_3}^2 {x_4}^2 {x_1}^3}{1323}-\frac{23}{63} {x_2}^2 {x_3} {x_4}^2 {x_1}^3-14742 {x_2}^3 {x_5}^2
   {x_1}^3-324 {x_2} {x_3} {x_5}^2 {x_1}^3\\
   &+\frac{33614}{3} {x_2}^2 {x_6}^2 {x_1}^3+1372 {x_3} {x_6}^2
   {x_1}^3-\frac{136323}{784} {x_2}^7 {x_3} {x_1}^3-10935 {x_2}^6 {x_5} {x_1}^3-\frac{288}{343} {x_3}^3 {x_5}
   {x_1}^3\\
   &+\frac{1458}{49} {x_2}^2 {x_3}^2 {x_5} {x_1}^3-\frac{2025}{7} {x_2}^4 {x_3} {x_5} {x_1}^3+\frac{1960}{3} {x_2}^3
   {x_4} {x_6} {x_1}^3-\frac{140}{9} {x_2} {x_3} {x_4} {x_6} {x_1}^3+\frac{19}{54} {x_2}^2 {x_4}^3
   {x_1}^2\\
   &-\frac{4802}{9} {x_4} {x_6}^2 {x_1}^2+\frac{20817}{224} {x_2}^7 {x_4} {x_1}^2+\frac{19 {x_2} {x_3}^3 {x_4}
   {x_1}^2}{2401}+\frac{9}{196} {x_2}^3 {x_3}^2 {x_4} {x_1}^2+\frac{675}{392} {x_2}^5 {x_3} {x_4} {x_1}^2\\
   &-603 {x_2}^4
   {x_4} {x_5} {x_1}^2+\frac{24}{49} {x_3}^2 {x_4} {x_5} {x_1}^2+\frac{81}{7} {x_2}^2 {x_3} {x_4} {x_5}
   {x_1}^2-\frac{79947}{16} {x_2}^6 {x_6} {x_1}^2+\frac{26}{49} {x_3}^3 {x_6} {x_1}^2\\
   &-\frac{171}{14} {x_2}^2 {x_3}^2 {x_6}
   {x_1}^2-\frac{245}{27} {x_2} {x_4}^2 {x_6} {x_1}^2+\frac{2025}{4} {x_2}^4 {x_3} {x_6} {x_1}^2-882 {x_2}^3 {x_5}
   {x_6} {x_1}^2+882 {x_2} {x_3} {x_5} {x_6} {x_1}^2\\
   &+\frac{85293 {x_2}^{10} {x_1}}{125440}+\frac{279 {x_3}^5
   {x_1}}{8235430}-\frac{675 {x_2}^2 {x_3}^4 {x_1}}{941192}+\frac{7 {x_4}^4 {x_1}}{2916}+\frac{5265 {x_2}^4 {x_3}^3
   {x_1}}{134456}-\frac{2187 {x_2}^6 {x_3}^2 {x_1}}{10976}+\frac{117}{28} {x_2}^5 {x_4}^2 {x_1}\\
   &-\frac{2 {x_2} {x_3}^2
   {x_4}^2 {x_1}}{1029}-\frac{1}{49} {x_2}^3 {x_3} {x_4}^2 {x_1}+\frac{4131}{2} {x_2}^4 {x_5}^2 {x_1}+\frac{162}{49} {x_3}^2
   {x_5}^2 {x_1}-\frac{486}{7} {x_2}^2 {x_3} {x_5}^2 {x_1}+1029 {x_2}^3 {x_6}^2 {x_1}\\
   &-196 {x_2} {x_3} {x_6}^2
   {x_1}-4802 {x_5} {x_6}^2 {x_1}+\frac{2187 {x_2}^8 {x_3} {x_1}}{6272}-\frac{13851}{112} {x_2}^7
   {x_5} {x_1}+\frac{216 {x_2} {x_3}^3 {x_5} {x_1}}{2401}
   +\frac{81}{98} {x_2}^3 {x_3}^2 {x_5} {x_1}\\
   &-\frac{19}{2} {x_2}^2
   {x_4}^2 {x_5} {x_1}-\frac{6075}{196} {x_2}^5 {x_3} {x_5} {x_1}-189 {x_2}^4 {x_4} {x_6} {x_1}-\frac{8}{21} {x_3}^2
   {x_4} {x_6} {x_1}+6 {x_2}^2 {x_3} {x_4} {x_6} {x_1}\\
   &+98 {x_2} {x_4} {x_5} {x_6} {x_1}+\frac{73 {x_2}^3
   {x_4}^3}{1512}+\frac{16807 {x_6}^3}{27}+54 {x_2}^2 {x_4} {x_5}^2+\frac{343}{6} {x_2} {x_4} {x_6}^2+\frac{243 {x_2}^8
   {x_4}}{6272}-\frac{5 {x_3}^4 {x_4}}{235298}\\
   &+\frac{15 {x_2}^2 {x_3}^3 {x_4}}{67228}+\frac{81 {x_2}^6 {x_3}
   {x_4}}{10976}-\frac{2 {x_4}^3 {x_5}}{81}-\frac{459}{112} {x_2}^5 {x_4} {x_5}-\frac{9}{686} {x_2} {x_3}^2 {x_4}
   {x_5}-\frac{45}{49} {x_2}^3 {x_3} {x_4} {x_5}+\frac{243 {x_2}^7 {x_6}}{224}\\
   &-\frac{15}{686} {x_2} {x_3}^3
   {x_6}-\frac{91}{36} {x_2}^2 {x_4}^2 {x_6}+\frac{2}{27} {x_3} {x_4}^2 {x_6}-1323 {x_2} {x_5}^2 {x_6}-\frac{81}{112}
   {x_2}^5 {x_3} {x_6}+\frac{945}{8} {x_2}^4 {x_5} {x_6}-\frac{9}{7} {x_3}^2 {x_5} {x_6}\\
  &+ \frac{14}{27} {g_1}^2 {x_1} \left(28 {x_1}^2-9 {x_2}\right)-\frac{2}{63}{g_1} {g_2} \left(7546 {x_1}^4 {x_2}-1323 {x_1}^2 {x_2}^2+84 {x_1}^2
   {x_3}+252 {x_2}^3-9 {x_2} {x_3}\right)\\
   &+\frac{224
   }{135}{g_1}{g_2}^2+{g_1}{g_3}-\frac{7}{5} {g_1} {g_2} z^2+\frac{119
   {g_2}^3 z}{135}+z^7.\\
 \end{array}
   $$
   }
Here $z$ is an algebraic function of $x_1,x_2,\cdots,x_7$ defined by the equation:
$$z^3 + g_2z + g_1=0.
$$


\section{Three families associated to isolated singularities
of types $E_{6},E_{7},E_{8}$}
\label{section:sec5}

In this section, we introduce three families of hypersurfaces in ${\bf C}^3$
which are deformations of simple singularities of type  $E_n$ ($n=6,7,8$)  denoted by
${\cal F}_{E_n(1)}$ and study the relationship between these families and the algebraic
potentials of Theorem \ref{theorem:theorem7}  (v), (vi), (vii) which
are denoted by $E_n(1)$ ($n=6,7,8$), respectively.
The family plays the role of the versal family of $E_n$-singularities
when the polynomial potential of type $E_n$ is  replaced by  the algebraic potential $E_n(1)$.

\subsection{Three families associated to  isolated singularities
of types $E_{12},E_{13},E_{14}$}
\label{subsection:sec5-1}

In the notes \cite{SL}, H. Shiga
introduced a family of polarized K3 surfaces defined by the equation
$$
S(t):z^2=y^3+(t_4+t_{10}x^3)y+(x^7+t_6x^6+t_{12}x^5+t_{18}x^4).
$$
The notes \cite{SL} are based on the paper by Nagano and Shiga \cite{NS}.
In page 4 of \cite{NS}, the  polynomial $d_{90}(t)$ is introduced.
$d_{90}(t)$ is a polynomial of $t=(t_4,t_6,t_{10}t_{12},t_{18})$,
 is the discriminant of $S(t)$ and  is regarded as the discriminant of Burkhardt group
which is nothing but the complex reflection group ST33 in \cite{ST}.

The definition of the family $S(t)$ leads us to
the existence of three families of hypersurfaces which are deformations
of the hypersurface isolated singularities
of types $E_{12},E_{13},E_{14}$ in the sense of Arnol'd.
They read
$$
\begin{array}{ll}
{\cal G}_{E_{12}R}:&z^2=y^3+t_{9}x^4+t_{6}x^5+t_3x^6+x^7+y(t_{5}x^3+t_2x^4)+s_4xy^2,\\
{\cal G}_{E_{13}R}:&z^2=y^3+t_{7}x^4+t_{5}x^5+t_3x^6+t_1x^7+y(t_{4}x^3+t_2x^4+x^5)+s_3xy^2,\\
{\cal G}_{E_{14}R}:&z^2=y^3+t_{12}x^4+t_{9}x^5+t_6x^6+t_3x^7+x^8+y(t_{7}x^3+t_4x^4+t_1x^5)+s_5xy^2.\\
\end{array}
$$
The family ${\cal G}_{E_{12}}$ is slightly different from  $S(t)$, since it contains the new parameter $s_4$.

\subsection{Three families associated to  isolated singularities
of types $E_{6},E_{7},E_{8}$}
\label{subsection:sec4-2}
In page 30 of \cite{SL}, the following family of elliptic surfaces is introduced:
$$
\begin{array}{lll}
x_2^2&=&GN(Y,Z),\\
GN(Y,Z)&=&Y^3+q_{18}Z^2+q_{12}Z^4+q_6Z^6+Z^8+Y(r_{10}Z^2+r_4Z^4).\\
\end{array}
$$
The definition of the above family 
leads us to
construct  three families of hypersurfaces which are deformations
of the hypersurface simple singularities
of types $E_{6},E_{7},E_{8}$.
$$
\begin{array}{ll}
{\cal F}_{E_6(1)}:&z^2=f_{E_6(1)},\\
{\cal F}_{E_7(1)}:&z^2=f_{E_7(1)},\\
{\cal F}_{E_8(1)}:&z^2=f_{E_8(1)},\\
\end{array}
$$
where
$$
\begin{array}{lll}
f_{E_6(1),(t_2,t_3,t_5,t_6,t_9,s_4)}(x,y)&=&y^3+t_{9}x+t_{6}x^2+t_3x^3+x^4+y(t_{5}x+t_2x^2)+s_4y^2,\\
f_{E_7(1),(t_1,\ldots,t_5,t_7,s_3)}(x,y)&=&y^3+t_{7}x+t_{5}x^2+t_3x^3+t_1x^4+x^3y+y(t_{4}x+t_2x^2)+s_3y^2,\\
f_{E_8(1),(t_1,\ldots,t_{12},s_5)}(x,y)&=&y^3+t_{12}x+t_{9}x^2+t_6x^3+t_3x^4+x^5+y(t_{7}x+t_4x^2+t_1x^3)+s_5y^2.\\
\end{array}
$$

We find the following relationship between
${\cal F}_{E_6(1)}$ and ${\cal F}_{E_{12R}}$ 

\begin{tabular}{l}
``The resultant of ${\cal F}_{E_6(1)}$ w.r.t. $y$'' $=xR_6(s,t)$\\
``The resultant of ${\cal F}_{E_{12R}}$ w.r.t. $y$'' $=x^7R_6(s,t)$\\
\end{tabular}

Here $R_6=27x^7+54t_3x^6+\cdots$.

Similarly we find the following relationship between
${\cal F}_{E_7(1)}$ and ${\cal F}_{E_{13R}}$ 

\begin{tabular}{l}
``The resultant of ${\cal F}_{E_7(1)}$ w.r.t. $y$'' $=xR_7(s,t)$\\
``The resultant of ${\cal F}_{E_{13R}}$ w.r.t. $y$'' $=x^7R_7(s,t)$\\
\end{tabular}

Here $R_7=4x^8+12t_3x^7+\cdots$.

Also we find the following relationship between
${\cal F}_{E_8(1)}$ and ${\cal F}_{E_{14R}}$ 

\begin{tabular}{l}
``The resultant of ${\cal F}_{E_8(1)}$ w.r.t. $y$'' $=xR_8(s,t)$\\
``The resultant of ${\cal F}_{E_{14R}}$ w.r.t. $y$'' $=x^7R_8(s,t)$\\
\end{tabular}

Here $R_8=27x^9+54t_3x^8+\cdots$.

\subsection{An observation}
\label{subsection:sec5-3}
We focus our attention on the cases $E_6,\>E_7,\>E_8$ in TABLE III.
The degrees of basic invariants of the reflection groups $W(E_n)\>(n=6,7,8)$.
 were collected in TABLE III  (cf. \S3).

Eliminate the greatest number  in the second column of TABLE III
and include a number so that the duality among the new set of 
numbers is satisfied.
Then the numbers in the second column of TABLE IV is obtained.
The number in the last column of TABLE IV is the half of the second greatest number of the second column.

\vspace{5mm}
TABLE IV\\
$
\begin{array}{l|l|l}\hline
E_6&2,\underline{3},5,6,8,9&4\\
E_7&2,\underline{4},6,8,10,12,14&6\\
E_8&2,\underline{6},8,12,14,18,20,24&10\\\hline
\end{array}
$

\vspace{5mm}
Eliminating the second greatest number in the column of TABLE IV,
we obtain the numbers in the second column in TABLE V.
The numbers in TABLE V are the same as the indices of the parameters $t_{j_1},\ldots,t_{j_n}$
of the polynomial $f_{E_n(1),(t,s)}$.

\vspace{5mm}
TABLE V\\
$
\begin{array}{l|l|l}\hline
ST33&2,3,5,6,9&4\\
ST34&2,4,6,8,10,14&6\\
{\rm No}&2,8,12,14,18,24&10\\\hline
\end{array}
$

\vspace{5mm}
It is underlined here that the degrees of the basic invariants of the groups $ST33$  and $ST34$  are restored by
this procedure.
Here $ST33$ and $ST34$ mean  the finite complex reflection groups denoted No.33 and No.34 in the list of
Shephard-Todd \cite{ST}.


\subsection{From TABLE IV to TABLE V}
\label{subsection:sec5-5}

In this subsection,  $n$ is always assumed one of the numbers $6,7,8$
and take  the variables $x_1,x_2,\cdots,x_n$  with the weights
$w_1,w_2,\cdots,w_n$ which are the same as the variables
taken  for the algebraic potential $E_{n}(1)\>(n=6,7,8)$.

In each of the cases $E_{6}(1),\>E_{7}(1),\>E_{8}(1)$, there is a weighted homogeneous polynomial $v(x_1,\cdots,x_{n-2})$
such that the algebraic function $z$ is defined by the equation
$$
z^2-x_{n-1}+v(x_1,\cdots,x_{n-2})=0.
$$
From the algebraic potential $F$,
we define $P$ and $C$ as in \S2.1.
It follows from the definition that $\partial_{x_n}C=I_n$.
Moreover it is possible to show that
$$
[\partial_{x_j}C,\>\partial_{x_k}C]=O\quad (\forall j,k).
$$
Let $E$ be the Euler vector field, that is,
$$
E=\sum_{j=1}^nw_jx_j\partial_{x_j}.
$$
By using $E$ and $C$, we define
the $n\times n$ matrix $T=EC$.
It is clear from the definition that $\Delta_{E_n(1)}=\det(T)$ is an algebraic function and is regarded a monic
polynomial of $x_n$.

Noting the relation $x_{n-1}=z^2-v(x_1,\cdots,x_{n-2})$,
we regard $z$ a variable and $x_{n-1}$ a polynomial of $x_1,\ldots,x_{n-2},x_n, z$ and define
the function
$$
g_{E_n}(x_1,\cdots,x_{n-2},x_n,z)=\Delta_{E_n(1)}|_{x_{n-1}=z^2-v}.
$$

We are going to penetrate another nature of $g_{E_n}$.
It is clear from the definition that $g_{E_n}$ is a polynomial of $x_1,\cdots,x_{n-2},x_n,z$
and is weighted homogeneous of weights
$w_1,w_{n-2},w_n,w_{n-1}/2$.
It is expected and actually will be shown later in this section that $g_{E_n}=0$ defines a free divisor.

We compute the restriction of $g_{E_n}$ to the hyperplane $z=0$.
For this purpose, we introduce the $n\times n$ matrix $U=I_n+U_0$, where $U_0$ is defined  as follows.
Let $u_{ij}$ denote the $(i,j)$-entry of $U_0$.
Then
$$
u_{ij}=\left\{\begin{array}{ll}
\partial_{x_i}v&(i<n-1,\>j=n-1)\\
0&({\rm others)}\\
\end{array}\right.
$$
Using $U$, we put ${\tilde T}=UTU^{-1}|_{z=0}$.
If ${\tilde T}_{ij}$ is the $(i,j)$-entry of ${\tilde T}$,
it can be shown that
${\tilde T}_{n-1\>j}=0\>(j\not=n-1)$ and ${\tilde T}_{n-1\>n-1}=\varphi$ for a polynomial $\varphi$
of $x_1,\ldots,x_{n-2},x_n$.

Since $\det({\tilde T})=\det(T)|_{z=0}=g_{E_n}|_{z=0}$,
it follows that $g_{E_n}|_{z=0}=\varphi\cdot \delta_{E_n}$, where
$\delta_{E_n}=\delta_{E_n}(x_1,\cdots,x_{n-2},x_n)$
is a polynomial of $x_1,\cdots,x_{n-2},x_n$ and
takes the form
$$
\delta_{E_n}={x_n}^{n-1}+q_1{x_n}^{n-2}+\ldots+q_{n-1},
$$
where
$q_1,\cdots,q_{n-1}$ are polynomials of $x_1,\cdots,x_{n-2}$.

The following theorem gives the answer to the question why the degrees of the basic invariants
of the groups $ST33,\>ST34$ appear in TABLE IV.

\begin{theorem}
\label{theorem:dsc-st33-st34}
(1)
$\delta_{E_6}$ is regarded as  the discriminant of the complex reflection group ST33
up to a constant factor.

(2)
$\delta_{E_7}$ is regarded as the discriminant of the complex reflection group ST34
up to a constant factor.
\end{theorem}

{\bf Proof.} 
We prove (2) here.
The proof of (1) is shown by an argument to that of  (2).

We reproduce here the potential vector field
$P_{ST34}=(h_1,h_2,\ldots,h_6)$  of the group $ST34$
given in \cite{KMS1} for the sake of
completeness.  It reads
{\footnotesize
$$
\begin{array}{ll}
h_1&\hspace{-2.5mm}=(20y_1^4y_2^2 + 120y_1^2y_2^3 - 60y_2^4 - 12y_1^5y_3 + 
    60y_1^3y_2y_3
 - 180y_1y_2^2y_3 + 135y_1^2y_3^2 + 135y_2y_3^2 + 
    120y_1^4y_4 \\
&\hspace{2mm}- 
        180y_1^2y_2y_4 + 540y_2^2y_4 + 270y_1y_3y_4 + 405y_4^2 + 
    180y_1^3y_5 + 540y_1y_2y_5 + 405y_3y_5 + 405y_1y_6)/405, \\
h_2&\hspace{-2.5mm}= (64y_1^9 - 288y_1^7y_2 - 1728y_1^5y_2^2 + 
    1728y_1^3y_2^3 - 2592y_1y_2^4 + 432y_1^6y_3
 + 3888y_1^4y_2y_3 + 
    3888y_1^2y_2^2y_3 \\
&\hspace{2mm}+ 
        1944y_2^3y_3 + 972y_1^3y_3^2 + 729y_3^3 - 1296y_1^5y_4 + 
    7776y_1^3y_2y_4 + 8748y_1^2y_3y_4 - 2916y_2y_3y_4 - 
    5832y_1y_4^2\\
&\hspace{2mm} + 3888y_1^4y_5 - 
        2916y_1^2y_2y_5 - 2916y_2^2y_5 + 4374y_1y_3y_5 + 4374y_4y_5 + 
    4374y_2y_6)/4374, \\
 h_3&\hspace{-2.5mm}= (64y_1^{10} + 3456y_1^8y_2 + 6624y_1^6y_2^2 + 
    2160y_1^4y_2^3 + 3888y_2^5 + 7776y_1^7y_3 + 19440y_1^5y_2y_3 + 
    32400y_1^3y_2^2y_3\\
&\hspace{2mm}
 + 
        7776y_1y_2^3y_3 + 19440y_1^4y_3^2 + 26244y_1^2y_2y_3^2 + 
    2916y_2^2y_3^2 + 7290y_1y_3^3 + 9504y_1^6y_4 + 38880y_1^4y_2y_4\\
&\hspace{2mm} - 
        11664y_1^2y_2^2y_4 - 3888y_2^3y_4 + 25272y_1^3y_3y_4 + 
    34992y_1y_2y_3y_4 + 8748y_3^2y_4 + 5832y_1^2y_4^2 + 
    8748y_2y_4^2 \\
&\hspace{2mm}+ 11664y_1^5y_5 + 
        25272y_1^3y_2y_5 + 17496y_1y_2^2y_5 + 26244y_1^2y_3y_5 + 
    8748y_2y_3y_5 + 26244y_1y_4y_5\\
&\hspace{2mm} + 6561y_5^2 + 6561y_3y_6)/6561, \\
h_4&\hspace{-2.5mm}=(1152y_1^{11} + 832y_1^9y_2 + 23616y_1^7y_2^2 - 
    13824y_1^5y_2^3 + 34560y_1^3y_2^4 + 7776y_1^8y_3 + 
    42768y_1^6y_2y_3 + 
        16200y_1^4y_2^2y_3 \\
&\hspace{2mm}+ 7776y_1^2y_2^3y_3 - 5832y_2^4y_3 + 
    17496y_1^5y_3^2 + 27216y_1^3y_2y_3^2 + 23328y_1y_2^2y_3^2 + 
    13122y_1^2y_3^3 + 
        2916y_2y_3^3 \\
&\hspace{2mm}+ 12960y_1^7y_4 - 25920y_1^5y_2y_4 + 
    103680y_1^3y_2^2y_4 - 31104y_1y_2^3y_4 + 19440y_1^4y_3y_4 + 
    69984y_1^2y_2y_3y_4\\
&\hspace{2mm} - 
        5832y_2^2y_3y_4 + 17496y_1y_3^2y_4 + 38880y_1^3y_4^2 - 
    46656y_1y_2y_4^2 - 4374y_3y_4^2 + 2592y_1^6y_5 + 
    21384y_1^4y_2y_5\\
&\hspace{2mm} - 
        5832y_1^2y_2^2y_5 + 11664y_2^3y_5 + 26244y_1^3y_3y_5 + 
    17496y_1y_2y_3y_5 + 6561y_3^2y_5 - 8748y_1^2y_4y_5 + 
    17496y_2y_4y_5 \\
&\hspace{2mm}+ 
        13122y_1y_5^2 + 13122y_4y_6)/13122,\\
h_5&\hspace{-2.5mm}= (10496y_1^{12} + 
    70656y_1^{10}y_2 + 86976y_1^8y_2^2 + 233856y_1^6y_2^3 - 
    25920y_1^4y_2^4 - 20736y_2^6 + 
        71808y_1^9y_3 \\
&\hspace{2mm}
+ 264384y_1^7y_2y_3 + 393984y_1^5y_2^2y_3 + 
    129600y_1^3y_2^3y_3 + 93312y_1y_2^4y_3 + 165888y_1^6y_3^2 + 
    408240y_1^4y_2y_3^2\\
&\hspace{2mm} + 
        221616y_1^2y_2^2y_3^2 + 134136y_1^3y_3^3 + 87480y_1y_2y_3^3 + 
    10935y_3^4 + 22464y_1^8y_4 + 209088y_1^6y_2y_4 + 
    38880y_1^4y_2^2y_4 \\
&\hspace{2mm}+ 
        233280y_1^2y_2^3y_4 + 23328y_2^4y_4 + 241056y_1^5y_3y_4 + 
    272160y_1^3y_2y_3y_4 + 139968y_1y_2^2y_3y_4 + 
    209952y_1^2y_3^2y_4\\
&\hspace{2mm} + 
        87480y_2y_3^2y_4 + 19440y_1^4y_4^2 + 349920y_1^2y_2y_4^2 - 
    69984y_2^2y_4^2 + 139968y_1y_3y_4^2 - 17496y_4^3 + 10368y_1^7y_5\\
&\hspace{2mm} + 
        38880y_1^5y_2y_5 + 93312y_1^3y_2^2y_5 - 23328y_1y_2^3y_5 + 
    134136y_1^4y_3y_5 + 227448y_1^2y_2y_3y_5 + 52488y_2^2y_3y_5 \\
&\hspace{2mm}+ 
        104976y_1y_3^2y_5 + 198288y_1^3y_4y_5 + 104976y_1y_2y_4y_5 + 
    78732y_3y_4y_5 + 78732y_1^2y_5^2 + 26244y_2y_5^2\\
&\hspace{2mm} + 39366y_5y_6)/39366, \\
 h_6&\hspace{-2.5mm}= (109056y_1^{14} + 433664y_1^{12}y_2 + 1983744y_1^{10}y_2^2 - 
    400512y_1^8y_2^3 + 2784768y_1^6y_2^4 - 282240y_1^4y_2^5\\
&\hspace{2mm} + 
        967680y_1^2y_2^6 + 207360y_2^7 + 403200y_1^{11}y_3 + 
    2395008y_1^9y_2y_3 + 3709440y_1^7y_2^2y_3 + 
    6096384y_1^5y_2^3y_3\\
&\hspace{2mm} - 846720y_1^3y_2^4y_3 - 
        725760y_1y_2^5y_3 + 1611792y_1^8y_3^2 + 
    4445280y_1^6y_2y_3^2 + 3492720y_1^4y_2^2y_3^2 + 
    1360800y_1^2y_2^3y_3^2\\
&\hspace{2mm} + 462672y_2^4y_3^2 + 
        1496880y_1^5y_3^3 + 2857680y_1^3y_2y_3^3 + 
    734832y_1y_2^2y_3^3 + 489888y_1^2y_3^4 + 91854y_2y_3^4 \\
&\hspace{2mm}+ 
    177408y_1^{10}y_4 + 80640y_1^8y_2y_4 + 
        3499776y_1^6y_2^2y_4 - 3870720y_1^4y_2^3y_4 + 
    2540160y_1^2y_2^4y_4 - 653184y_2^5y_4 \\
&\hspace{2mm}+ 1439424y_1^7y_3y_4 + 
    7039872y_1^5y_2y_3y_4 + 
        3991680y_1^3y_2^2y_3y_4 + 1741824y_1y_2^3y_3y_4 + 
    2721600y_1^4y_3^2y_4 \\
&\hspace{2mm}+ 1388016y_1^2y_2y_3^2y_4 + 
    244944y_2^2y_3^2y_4 + 
        857304y_1y_3^3y_4 + 1874880y_1^6y_4^2 - 
    3175200y_1^4y_2y_4^2\\
&\hspace{2mm} + 5225472y_1^2y_2^2y_4^2 + 
    925344y_1^3y_3y_4^2 + 2939328y_1y_2y_3y_4^2 + 
        306180y_3^2y_4^2 + 1469664y_1^2y_4^3 - 816480y_2y_4^3 \\
&\hspace{2mm}- 
    8064y_1^9y_5 + 254016y_1^7y_2y_5 + 217728y_1^5y_2^2y_5 + 
    362880y_1^3y_2^3y_5 + 
        544320y_1y_2^4y_5 + 707616y_1^6y_3y_5\\
&\hspace{2mm} + 
    1632960y_1^4y_2y_3y_5 + 2122848y_1^2y_2^2y_3y_5 - 
    326592y_2^3y_3y_5 + 1714608y_1^3y_3^2y_5 + 
        1224720y_1y_2y_3^2y_5\\
&\hspace{2mm} + 183708y_3^3y_5 + 816480y_1^5y_4y_5 + 
    4572288y_1^3y_2y_4y_5 - 1959552y_1y_2^2y_4y_5 + 
    2694384y_1^2y_3y_4y_5\\
&\hspace{2mm} + 
        979776y_2y_3y_4y_5 - 244944y_1y_4^2y_5 + 1102248y_1^4y_5^2 + 
    979776y_1^2y_2y_5^2 + 489888y_2^2y_5^2 + 734832y_1y_3y_5^2\\
&\hspace{2mm} + 
    367416y_4y_5^2 + 
        137781y_6^2)/275562.\\
\end{array}
$$
}
The $6\times 6$ matrix
$$
C_{ST34}={}^t(\partial_{y_1}P_{ST34},\partial_{y_2}P_{ST34},\cdots,\partial_{y_6}P_{ST34})
$$
is defined similarly to the definition of $C$.
Let $E_{ST34}$ be the Euler vector field defined by
$E_{ST34}=\frac{1}{7}(\sum_{j=1}^5jy_j\partial_{y_j}+7y_6\partial_{y_6}).$
Then $T_{ST34}=E_{ST34}C_{ST34}$ plays the role  similar to $T$.
It is shown in \cite{S8} that $\det(T_{ST34})$ is the discriminant of the complex reflection group $ST34$.

On the other hand, it can be shown by direct computation that
by the change of variables
$$
\begin{array}{lll}
 x_1 &=&\frac{2^3\cdot 3^4}{5^2\cdot 7}y_1,\\[1mm]
   x_2 &=& -\frac{2^6\cdot 3^8}{5^4\cdot 7^3}(y_1^2 + 12y_2),\\[1mm]
 x_3 &=&\frac{2^9\cdot 3^{11}}{5^6 \cdot 7^4}(11y_1^3 + 36y_1y_2 + 54y_3), \\[1mm]
   x_4 &=& -\frac{2^{11}\cdot 3^{15}}{5^8\cdot 7^5}(77y_1^4 - 216y_1^2y_2 + 432y_2^2 + 
         216y_1y_3 + 1296y_4),\\[1mm]
   x_5 &=& -\frac{2^{15}\cdot 3^{19}}{5^{10}\cdot 7^6}(175y_1^5 - 168y_1^3y_2 + 1296y_1y_2^2 - 
         1188y_1^2y_3 - 
                1296y_2y_3 - 1296y_1y_4 - 3888y_5), \\[1mm]
x_7 &=&\frac{2^{21}\cdot 3^{26}} {5^{13}\cdot 7^9}
(1645y_1^7 + 
       9828y_1^5y_2 + 
             3024y_1^3y_2^2 + 108864y_1y_2^3 - 38178y_1^4y_3 + 
       27216y_1^2y_2y_3\\[1mm]
       && - 
             69984y_2^2y_3 
             + 64152y_1y_3^2 + 117936y_1^3y_4 - 
       139968y_1y_2y_4 + 
             69984y_3y_4 + 128304y_1^2y_5\\[1mm]
             && + 139968y_2y_5 + 
       209952y_6), \\
\end{array}
$$
$\delta_{E_7}$ turns out to be $\det(T_{ST34})$.
Thus the theorem is proved.[]

\vspace{5mm}

Concerning  the case $E_8(1)$, we have the follwing result in spite that
there is no complex reflection group corresponding to $E_8(1)$.
To state the result, we introduce the following coordinate transformation:
{\footnotesize
\begin{equation}
\label{equation:E8(1)PV}
\left\{
\begin{array}{lll}
{y_1}&=& 2 {x_1},\\
{y_2}&=& -
\frac{2}{3} \left({x_1}^3-12 {x_2}\right),\\
{y_3}&=& -\frac{3}{2} \left({x_1}^4+8 {x_1} {x_2}-4
   {x_3}\right),\\
   {y_4}&=& \frac{2}{3} \left({x_1}^6+30 {x_1}^3 {x_2}-9 {x_1}^2 {x_3}-36 {x_2}^2+36 {x_4}\right),\\
   {y_5}&=&
   \frac{3}{14} \left(15 {x_1}^7-84 {x_1}^4 {x_2}+336 {x_1} {x_2}^2-56 {x_1} {x_4}-56 {x_2} {x_3}+28
   {x_5}\right),\\
   {y_6}&=& \frac{1}{36} \left(-763 {x_1}^9+576 {x_1}^6 {x_2}+1080 {x_1}^5 {x_3}-24192 {x_1}^3 {x_2}^2-432
   {x_1}^3 {x_4}+6480 {x_1}^2 {x_2} {x_3}\right.\\
   &&\left.-648 {x_1}^2 {x_5}-648 {x_1} {x_3}^2+10368 {x_2}^3-5184 {x_2} {x_4}+2592
   {x_6}\right),\\
   {y_7}&=&\frac{1}{11880}( 582775 {x_1}^{12}+4051608 {x_1}^9 {x_2}-2061180 {x_1}^8 {x_3}-218592 {x_1}^6 {x_2}^2\\
   &&+3069792
   {x_1}^6 {x_4}-2965248 {x_1}^5 {x_2} {x_3}+270864 {x_1}^5 {x_5}+1104840 {x_1}^4 {x_3}^2\\
   &&+65577600 {x_1}^3
   {x_2}^3
   -3991680 {x_1}^3 {x_2} {x_4}-1425600 {x_1}^3 {x_6}-21384000 {x_1}^2 {x_2}^2 {x_3}\\
   &&+1710720 {x_1}^2 {x_2}
   {x_5}-855360 {x_1}^2 {x_3} {x_4}+2138400 {x_1} {x_2} {x_3}^2-427680 {x_1} {x_3} {x_5}\\
   &&-18817920 {x_2}^4+10264320
   {x_2}^2 {x_4}-1710720 {x_2} {x_6}-71280 {x_3}^3-855360 {x_4}^2+855360 {x_8})\\
   \end{array}
   \right.
   \end{equation}
   }
  Using this transformation, we define
  $$
  J_{y,x}=\frac{\partial(y_1,\ldots,y_7)}{\partial(x_1,\ldots,x_6,x_8)}=
  \left(
  \begin{array}{cccc}
  \frac{\partial y_1}{\partial x_1}&\cdots& \frac{\partial y_6}{\partial x_1}&
 \frac{\partial y_7}{\partial x_1}\\
 \vdots&\cdots&\vdots&\vdots\\
  \frac{\partial y_1}{\partial x_6}&\cdots& \frac{\partial y_6}{\partial x_6}&
 \frac{\partial y_7}{\partial x_6}\\
   \frac{\partial y_1}{\partial x_8}&\cdots& \frac{\partial y_6}{\partial x_8}&
 \frac{\partial y_7}{\partial x_8}\\
 \end{array}
 \right)
 $$
Moreover we define the $7\times 7$ matrix $T_A$ from ${\tilde T}$ which is an $8\times 8$ matrix in this case
by
$$
{T_A}'=\left(
\begin{array}{cccc}
{\tilde T}_{11}&\cdots&{\tilde T}_{16}&{\tilde T}_{18}\\
\vdots&\cdots&\vdots&\vdots\\
{\tilde T}_{61}&\cdots&{\tilde T}_{66}&{\tilde T}_{68}\\
{\tilde T}_{81}&\cdots&{\tilde T}_{86}&{\tilde T}_{68}\\
\end{array}
\right).
$$
Since each matrix entry of $T_A$ is a weighted homogeneous polynomial of $x_1,\ldots,x_6,x_8$,
it is easy to show that there is a $7\times 7$ matrix $C_A$
such that
$$
\left(\sum_{j=1}^6w_jx_j\partial_{x_j}+x_7\partial_{x_7}\right)C_A=T_A.
$$

\begin{theorem}
There exist  seven weighted homogeneous
polynomials $h_j\>(j=1,2,\cdots,7)$
of $y_1,\cdots,y_6,y_8$  with weights $w_1,\cdots,w_6,w_8$ satisfying the following conditions.

(M1)
$P_B=(h_1,\cdots,h_7)$ is a potential vector field.

(M2)  ${}^t(\partial_{y_1}P_B,\cdots,\partial_{y_6}P_B,\partial_{y_8}P_B)=J_{y,x}^{-1}C_AJ_{y,x}$.

Here each matrix entry of $C_A$ is regarded as a polynomial of $y_1,\cdots,y_7$
under the transformation (\ref{equation:E8(1)PV}).

\end{theorem}

\vspace{5mm}
{\bf Proof.}
The theorem is shown by direct computation.
The most difficult part of the theorem is to determine the coordinate transformation
(\ref{equation:E8(1)PV}).


\subsection{The algebraic potentials of type $E_{n}(1) $ and the family ${\cal F}_{E_{n}(1)}$ ($n=6,7,8$) }
\label{subsection:sec5-6}

The study in \S\ref{subsection:sec5-5} 
leads us to the question 
on the relationship between
the polynomial $\Delta_{E_n(1)}$ and the discriminant of ${\cal F}_{E_{n}(1)}$ ($n=6,7,8$)

 \vspace{5mm}
 The answer to the question above is summarized in the following theorem.
 
\begin{theorem}
\label{theorem:theorem-4}

(i) \underline{$E_6$ case}

By the coordinate transformation 
$$
\begin{array}{lll}
{x_1}&=& -\frac{9}{4}
   {c_6}^5
   {t_2},\\
   {x_2}&=&
   \frac{3}{4} {c_6}^3
   {t_3},\\
   {x_3}&=&
   -\frac{27}{8} {c_6}^8
   (7 {t_2} {t_3}-18
   {t_5}),\\
   {x_4}&=&
   \frac{1}{16} {c_6}^6
   \left(108 {s_4}
   {t_2}-11 {t_2}^3+99
   {t_3}^2-324
   {t_6}\right),\\
   {x_6}&=&\frac{1}{5832}\left(540 {s_4}
   {t_2} {t_3}-1944
   {s_4} {t_5}-115
   {t_2}^3 {t_3}+270
   {t_2}^2 {t_5}+345
   {t_3}^3-1620 {t_3}
   {t_6}+5832
   {t_9}\right),\\
    {z}&=&
   {c_6}
   {s_4}\\
\end{array}
$$   
where $c_6$ is a constant satisfying
$
{c_6}^9=\frac{4}{2187},
$
$t_6\Delta_{E_{6}(1)}$ coincides with the discriminant of ${\cal F}_{E_{6}(1)}$
up to a non-zero constant factor.

(ii) \underline{$E_7$ case}
   
By  the coordinate transformation      
{\footnotesize
  $$ 
  \begin{array}{lll}
{x_1}&=& -3 {c_7}^2 {t_1},\\
{x_2}&=& -\frac{3}{7} {c_7}^4 \left(39
   {t_1}^2+28 {t_2}\right),\\
   {x_3}&=& \frac{3}{7} {c_7}^6 \left(-392
   {s_3}-1581 {t_1}^3-1428 {t_1} {t_2}+1176 {t_3}\right),\\
   {x_4}&=&
   \frac{1}{64} {c_7} \left(-32 {s_3} {t_1}-105 {t_1}^4-120 {t_1}^2
   {t_2}+96 {t_1} {t_3}-16 {t_2}^2+64 {t_4}\right),\\
   {x_5}&=& -\frac{3
   {c_7}^3}{10976} \left(215600 {s_3} {t_1}^2-21952 {s_3} {t_2}+585915
   {t_1}^5+816200 {t_1}^3 {t_2}-646800 {t_1}^2 {t_3}+215600 {t_1} {t_2}^2\right.\\
 &&\left.  -241472 {t_1} {t_4}-241472 {t_2} {t_3}+307328
   {t_5}\right),\\
   {x_7}&=& \frac{1}{843308032}
   (32269440 {s_3}^2 {t_1}-747369560
   {s_3} {t_1}^4-481736640 {s_3} {t_1}^2 {t_2}+408746240 {s_3}
   {t_1} {t_3}\\
   &&+32269440 {s_3} {t_2}^2-180708864 {s_3} {t_4}-1498061973
   {t_1}^7-2850681036 {t_1}^5 {t_2}+2242108680 {t_1}^4 {t_3}\\
   &&-1494739120
   {t_1}^3 {t_2}^2+686109760 {t_1}^3 {t_4}+2058329280 {t_1}^2 {t_2}
   {t_3}-613119360 {t_1}^2 {t_5}-160578880 {t_1} {t_2}^3\\
   &&+408746240
   {t_1} {t_2} {t_4}-613119360 {t_1} {t_3}^2+204373120 {t_2}^2
   {t_3}-301181440 {t_2} {t_5}-301181440 {t_3} {t_4}\\
   &&+843308032{t_7}),\\
   {z}&=& -2352 {c_7}^6 {s_3}\\
   \end{array}
$$
}
where $c_7$ is a constant defined by
${c_7}^7=\frac{1}{14112}$,
$t_7\Delta_{E_{7}(1)}$ coincides with the discriminant of ${\cal F}_{E_{7}(1)}$
up to a non-zero constant factor.

(iii) \underline{$E_8$ case}

By the coordinate transformation
\begin{equation}
\label{equation:trans-e8-case}
\left\{
\begin{array}{lll}
{x_1}&=& \frac{9}{16}  {c_8}^7t_1,\\
{x_2}&=& -\frac{1}{1024}{c_8}^9 \left(23 {t_1}^3+432 {t_3}\right),\\
{x_3}&=& \frac{1}{27648}{c_8}^4 \left(481
   {t_1}^4+11232 {t_1} {t_3}-20736 {t_4}\right),\\
   {x_4}&=& -\frac{1}{7962624}{c_8}^6 \left(1492992 {s_5} {t_1}+9061 {t_1}^6+301104 {t_1}^3
   {t_3}-528768 {t_1}^2 {t_4}\right.\\
   &&\left.+1586304 {t_3}^2-4478976 {t_6}\right),\\
   {x_5}&=&\begin{array}{l} \frac{1}{214990848}{c_8} \left(20901888 {s_5} {t_1}^2+134695
   {t_1}^7+5155920 {t_1}^4 {t_3}
   -8999424 {t_1}^3 {t_4}\right.\\
\left.  \hspace{18mm} +40497408 {t_1} {t_3}^2
-62705664 {t_1} {t_6}-62705664 {t_3} {t_4}+214990848
   {t_7}\right),\\\end{array}\\
   {x_6}&=& -{c_8}^3 (6897623040 {s_5} {t_1}^4+85136375808 {s_5} {t_1} {t_3}-185752092672 {s_5}
   {t_4}+40346383 {t_1}^9\\
   &&+1954725696 {t_1}^6 {t_3}-3391331328 {t_1}^5 {t_4}+25434984960 {t_1}^3 {t_3}^2-20692869120 {t_1}^3 {t_6}\\
   &&-62078607360
   {t_1}^2 {t_3} {t_4}+42568187904 {t_1}^2 {t_7}+42568187904 {t_1} {t_4}^2+62078607360 {t_3}^3\\
   &&-255409127424 {t_3} {t_6}+557256278016
   {t_9})/{743008370688},\\
   {x_8}&=& (-267483013447680 {s_5}^2 {t_1}^2-4240658644992 {s_5} {t_1}^7-118958736015360 {s_5} {t_1}^4
   {t_3}\\
   &&+215472427499520 {s_5} {t_1}^3 {t_4}
   -484812961873920 {s_5} {t_1} {t_3}^2+802449040343040 {s_5} {t_1} {t_6}\\
   &&+802449040343040 {s_5}
   {t_3} {t_4}-3851755393646592 {s_5} {t_7}-20239023875 {t_1}^{12}\\
   &&-1290950892000 {t_1}^9 {t_3}
   +2230763141376 {t_1}^8 {t_4}-26769157696512
   {t_1}^6 {t_3}^2\\
   &&+12721975934976 {t_1}^6 {t_6}
  +76331855609856 {t_1}^5 {t_3} {t_4}
   -23791747203072 {t_1}^5 {t_7}\\
   &&
   -59479368007680 {t_1}^4 {t_4}^2
   -190829639024640 {t_1}^3 {t_3}^3+356876208046080 {t_1}^3 {t_3} {t_6}\\
 &&  -161604320624640 {t_1}^3 {t_9}
   +535314312069120 {t_1}^2 {t_3}^2
   {t_4}-484812961873920 {t_1}^2 {t_3} {t_7}\\
  && -484812961873920 {t_1}^2 {t_4} {t_6}
   -484812961873920 {t_1} {t_3} {t_4}^2+802449040343040 {t_1}
   {t_4} {t_7}\\
 &&  +11555266180939776 {t_{12}}
   -267657156034560 {t_3}^4+1454438885621760 {t_3}^2 {t_6}\\
 &&  -2407347121029120 {t_3} {t_9}
   +133741506723840
   {t_4}^3\\
   &&-1203673560514560 {t_6}^2)/{11555266180939776},\\
   {z}&=& -\frac{81}{8}  {c_8}^{11}s_5,\\
   \end{array}
\right.
\end{equation}
where $c_8$ is a constant satisfying
   $
{c_8}^{12} = -16/243,
$
$t_8\Delta_{E_{8}(1)}$ coincides with the discriminant of ${\cal F}_{E_{8}(1)}$
up to a non-zero constant factor.
\end{theorem}

{\bf Proof:}
It is hard to prove the theorem without  help of the software Mathematica.
We shall explain the outline of the proof of  (iii) which is the most complicated case to accomplish.

\underline{Step 1} The determination of the discriminant of ${\cal F}_{E_{8}(1)}$.

We put $f=f_{E_{8}(1)}(t,s)$ for simplicity.
Let $S$ be the hypersurface of ${\bf C}^3$ defined by $z^2=f(x,y)$.
Then the singular locus of $S$ is reduced to that of the curve $S_0=\{(x,y)\in {\bf C}^2|\>f(x,y)=0\}$
of ${\bf C}^2$.
Let $C_f$ be the singular locus of $S_0$.
Then it follows from the definition that
$$
C_f=\{(x,y)\in{\bf C}^2|\>f=f_x=f_y=0\}.
$$
Our purpose is to construct the defining equation of $C_f$.
For this purpose,
we introduce polynomials $g_j\>(j=1,2,3)$ of $x,y$ with coefficients $a_i,b_j,c_k$ by
$$
\begin{array}{lll}
g_1&=&a_1 + a_2 x + a_3 x^2 + 
      a_4x^3 + (a_5 + a_6x + a_7x^2 + a_8x^3) y + 
      (a_9 + a_{10}x + a_{11}x^2 + a_{12}x^3) y^2,\\
g_2&=&b_1 + b_2 x + b_3 x^2 + b_4 x^3 + 
      b_5x^4 + (b_6 + b_7 x + b_8 x^2 + b_9 x^3 + b_{10}x^4) y\\
      && + (b_{11} + 
         b_{12} x + b_{13}x^2 + b_{14}x^3 + b_{15}x^4) y^2,\\
g_3&=&c_1 + c_2 x + c_3 x^2 + c_4 x^3 + 
      c_5x^4 + (c_6 + c_7 x + c_8 x^2 + c_9 x^3 + c_{10}x^4) y \\
      &&+ (c_{11} + 
         c_{12} x + c_{13}x^2 + c_{14}x^3 + c_{15}x^4) y^2\\
         \end{array}
         $$
and define $H=g_1f+g_2f_x+g_3f_3$.
By direct computation, we find that
the constants 
$b_j,c_k$ are so taken as the linear combinations of $a_j$ whose coefficients are polynomials of $t,s$
that
$$
H=p_1x^3y+{\tilde p}_2+p_3x^2y+p_4y+p_5x^3+p_6xy+p_7x^2+p_8x.
$$
Noting that ${\tilde p}_2$ has a factor $t_{12}$, we put $p_2={\tilde p}_2/t_{12}$.
Then
there are polynomials $p_{ij}\>(i,j=1,\ldots,8)$ of $t,s$ such that
$$
p_i=p_{i1}a_7+p_{i2}a_4+p_{i3}a_6+p_{i4}a_8+p_{i5}a_3+p_{i6}a_5+p_{i7}a_2+p_{i8}a_1.
$$
Let $P_A$ be the $8\times 8$ matrix whose $(i,j)$-entry is $p_{ij}$
and let ${\tilde P}_A$ be the matrix which is obtained by
changing the second row vector  of $P_A$ by $t_{12}(p_{21}\>\>\cdots\>p_{28})$
 and the remaining row vector  of $P_A$ left as it is.
Moreover let ${\bf v}$ be the vector ${}^t(x^3y\>\>1\>\>x^2y\>\>y\>\>x^3\>\>xy\>\>x^2\>\>x)$.
Then
we have
$$
(a_1\>a_2\>\cdots\>a_8){\tilde P}_A{\bf v}=p_1x^3y+t_{12}p_2+p_3x^2y+p_4y+p_5x^3+p_6xy+p_7x^2+p_8x.
$$
Since this equation holds for any $a_1,\ldots,a_8$, we conclude that if there is  a vector ${\bf v}\not={\bf o}$
such that $(a_1\>a_2\>\cdots\>a_8){\tilde P}_A{\bf v}={\bf o}$,
then $\det({\tilde P}_A)=0$.
As a consequence, if $C_f\not=\empty$ if and only if  $t_{12}\det(P_A)=0$.

\underline{Step 2} Modification of the matrix $P_A$ (I)

Let $P_j$ be the $j$-th row vector of $P_A$ and
define
$$
P_4'=\sum_{j=1}^8k_jP_j,
$$
where
$$
\left\{
\begin{array}{lll}
{k_1}&=& \frac{1}{25} \left(3 {t_1}^3+20 {t_3}\right),\\
{k_2}&=& \frac{1}{75} \left(50 {s_5}+12 {t_1}^2 {t_3}-25 {t_1}
   {t_4}\right),\\
   {k_3}&=& \frac{1}{25} \left(2 {t_1}^2 {t_4}+15 {t_6}\right),\\
   {k_4}&=& 1,\\
   {k_5}&=& \frac{1}{75} \left(40 {s_5}
   {t_3}+9 {t_1}^2 {t_6}-20 {t_1} {t_7}-10 {t_4}^2\right),\\
   {k_6}&=& \frac{1}{25} \left({t_1}^2 {t_7}+10
   {t_9}\right),\\
   {k_7}&=& \frac{1}{25} \left(10 {s_5} {t_6}+2 {t_1}^2 {t_9}-5 {t_4} {t_7}\right),\\
   {k_8}&=& \frac{1}{15} \left(4
   {s_5} {t_9}-{t_7}^2\right).\\
   \end{array}
   \right.
   $$
It is straightforward to show that $P_4'$ is divisible by $t_{12}$.
So we put $P_4''=\frac{1}{t_{12}}P_4'$
and define 
the matrix $P_B$ which is obtained from $P_A$ by changing the 4th row vector from$P_4$ with $P_4''$.
Then $\det (P_A)=t_{12}\cdot \det(P_B)$.

\underline{Step 3} Modification of the matrix $P_A$ (II)

We put
$Q_1=\partial_{t_{12}}(P_B)$
and
$$
Q_2=
\left(
\begin{array}{cccccccc}
 1 & 0 & 0 & 0 & 0 & 0 & 0 & 0 \\
 0 & -5 & -3 {t_1} & \frac{3 {t_7}}{2 {s_5}} & -\frac{8 {s_5} {t_3}-{t_1} {t_7}}{2 {s_5}} & -2 {t_4} & -\frac{6 {s_5}
   {t_6}-{t_4} {t_7}}{2 {s_5}} & -\frac{4 {s_5} {t_9}-{t_7}^2}{2 {s_5}} \\
 0 & 0 & 1 & 0 & 0 & 0 & 0 & 0 \\
 0 & 0 & 0 & -\frac{3}{2 {s_5}} & -\frac{{t_1}}{2 {s_5}} & 0 & -\frac{{t_4}}{2 {s_5}} & -\frac{{t_7}}{2 {s_5}} \\
 0 & 0 & 0 & 0 & 1 & 0 & 0 & 0 \\
 0 & 0 & 0 & 0 & 0 & 1 & 0 & 0 \\
 0 & 0 & 0 & 0 & 0 & 0 & 1 & 0 \\
 0 & 0 & 0 & 0 & 0 & 0 & 0 & 1 \\
\end{array}
\right).
$$
Using $Q_1$ and $Q_2$, we  define
$P_C={Q_1Q_2}^{-1}P_BQ_2$.

\underline{Step 4}
Let $T$ be the $8\times 8$ matrix obtained from the potential of $E_{8(1)}$.
Then $T$ is regarded as a matrix whose entries are polynomials of $t,s$ by the substitution
(\ref{equation:trans-e8-case}).
We modify $T$ so that the resulting matrix coincides with $P_C$.

For this purpose, we first introduce the Jacobi matrix between the variables $x$ and $t,s$.
Namely we define
$$
J_{x,t}=\frac{\partial(x_1,\ldots,x_7,x_8)}{\partial(t_1,t_3,t_4,s_5,t_6,t_7,t_9,t_{12})}=
\left(\begin{array}{cccc}
\frac{\partial x_1}{\partial t_1}&\frac{\partial x_2}{\partial t_1}&\cdots&\frac{\partial x_8}{\partial t_1}\\
\frac{\partial x_1}{\partial t_3}&\frac{\partial x_2}{\partial t_3}&\cdots&\frac{\partial x_8}{\partial t_{12}}\\
\vdots&\vdots&\cdots&\vdots\\
\frac{\partial x_1}{\partial t_{12}}&\frac{\partial x_2}{\partial t_{12}}&\cdots&\frac{\partial x_8}{\partial t_{12}}\\
\end{array}
\right)
$$
and
$$
J_{t,x}=\frac{\partial(t_1,t_3,t_4,s_5,t_6,t_7,t_9,t_{12})}{\partial(x_1,\ldots,x_7,x_8)}=(J_{x,t})^{-1}.
$$
In (\ref{equation:trans-e8-case}), $x_7$ doesn't appear.
But by the relation 
 (\ref{equation:def-ans-z-e8e}), we find that
 {\footnotesize
$$
\begin{array}{lll}
x_7&=& -\frac{{c_8}^{10}}{330225942528} \left(2229025112064 {s_5}^2+5486745600 {s_5} {t_1}^5+96745881600 {s_5} {t_1}^2 {t_3}-185752092672
   {s_5} {t_1} {t_4}\right.\\
&&\left.   +29988035 {t_1}^{10}+1605970800 {t_1}^7 {t_3}-2781475200 {t_1}^6 {t_4}+25033276800 {t_1}^4
   {t_3}^2-16460236800 {t_1}^4 {t_6}\right.\\
   &&\left.-65840947200 {t_1}^3 {t_3} {t_4}+32248627200 {t_1}^3 {t_7}+48372940800 {t_1}^2
   {t_4}^2+98761420800 {t_1} {t_3}^3\right.\\
   &&\left.-290237644800 {t_1} {t_3} {t_6}+278628139008 {t_1} {t_9}-145118822400 {t_3}^2
   {t_4}+278628139008 {t_3} {t_7}\right.\\
   &&\left.+278628139008 {t_4} {t_6}\right).\\
   \end{array}
   $$
   }
We put $T_Z=J_{x,t}TJ_{t,x}$ and
also regard $T_Z$ as a matrix whose entries are polynomials of $t,s$.
By direct computation, we find that $T_Z$ coincides $P_C$.

\underline{Step 5}

The conclusion of Step 4 implies that
$
t_{12}\Delta_{E_8(1)}$
 coincides with $\det(P_A)$ up to a non-zro constant factor.
 
 \begin{remark}
 \label{remark:Arnold}
 There are many problems called ``Arnolds Problem''.
 We cite here one of them.
In the book \cite{An}, p.20, it is written that

\vspace{5mm}
``1974-5 Find and application of the (Shephard-Todd) complex reflection groups to singularity theory''

\vspace{5mm}
This is one of Arnolds problems.
The families ${\cal F}_{E_{6(1)}}|_{s_4=0}$,  ${\cal F}_{E_{7(1)}}|_{s_3=0}$ are combined with
groups ST33 and ST34, respectivelyby Theorem \ref{theorem:dsc-st33-st34}
and Theorem \ref{theorem:theorem-4}.
As a consequence, the study of these families is expected to give a new insight 
on the role of the groups $ST33$ and $ST34$
in the deformation theory of simple hypersurface singularities of types $E_6$ and $E_7$.

There are many studies on the Arnolds problem
(cf. comments by V. V. Goryunov in \cite{An}, Goryunov \cite{Gr}).
It is noted here that  the Arnold's problem in the case of the Valentiner group (ST27)
was treated  in \cite{Sval}.
On the contrary to the cases treated by Goryunov and the cases treated here, the family corresponding to
the group $ST27$ is that of hypersurfaces
which are deformations
of the hypersurface $z^2=x^5-x^2y^3$
This has non-isolated  singularities.
\end{remark}
 
 \begin{remark}
We mention here the case of the family of hypersurfaces  corresponding to
the algebraic potential $F_{E_8(2)}$ (cf. \S3.3.2).
We introduce the family of hypersurfaces on the $xyz$-space
defined by
$$
{\cal H}_{E_8(2)}:\>z^2=f_{E_8(2)},
$$
where
$$
\begin{array}{lll}
f_{E_8(2)}(x,y)&=&x^5+y^3+
y \left({t_1} x^3+{t_{10}}+{t_4} x^2+{t_7} x\right)+y^2 ({t_2}
   x+{t_5})+\frac{{t_6}^2 x}{4}+{t_6} x^3\\[3mm]
&&+s_3\left(\frac{ {t_6}^2}{4}+ {t_6} x^2+x^4\right)\\
\end{array}
$$
is a polynomial with parameters $(t_1,t_2,t_4,t_5,t_6,t_7,t_{10},s_3)$.
Then the discriminant of the family ${\cal H}_{E_8(2)}$
is divided by 
$$
\frac{{t_1}^2 {t_6}^3}{8}-\frac{1}{2} {t_1} {t_6}^2 {t_7}+{t_{10}}^2-{t_{10}} {t_4} {t_6}+\frac{{t_4}^2
   {t_6}^2}{4}+\frac{{t_6} {t_7}^2}{2}
   $$
and its quotient coincides with $\det(T_{F_{E_8(2)}})$ up to a constant factor
by a certain coordinate
transformation between
$(t_1,t_2,t_4,t_5,t_6,t_7,t_{10},s_3)$
and $(x_1,x_2,x_4,x_5,x_6,x_7,x_{10},z)$.
The proof of this result is similar to that of the theorem and is left to the reader.

\end{remark}

\section{Flat coordinates and algebraic functions}
\label{section:sec4}

In the case of the reflection group of type $H_3$, there are one polynomial potential
$F_{H_3}$ (cf. (\ref{equation:potential-H3})) and two algebraic potentials $F_{(H_3)'}$, $F_{(H_3)''}$
(cf. (\ref{equation:ex-(H3)'case}),  (\ref{equation:ex-(H3)''case})).
It is known that the flat coordinates of a polynomial potential are regarded as basic invariants of the
corresponding reflection group.
Then it is a basic question whether  the flat coordinates and the algebraic functions
are written in terms of
 the basic invariants of $W(H_3)$ or not.
In this section, we discuss   this question.
Another question is to construct a family of hypersurfaces
which corresponds to algebraic potentials.
We discuss this problem in the cases $F_{H_3}$ and $F_{(H_3)''}$.
As to the former, this problem is already answered by T. Yano \cite{Ya}.
The questions similar to those explained above for the the reflection group
of type $H_4$ are also discussed.

\subsection{$H_3$ case}
In this subsection,
 the flat coordinates and algebraic functions to write $F_{(H_3)'}$
are denoted by $y_1,y_2,y_3$ and $w$ instead of $x_1,x_2,x_3$ and $z$,
and
 the flat coordinates and algebraic functions to write $F_{(H_3)''}$
are denoted by $u_1,u_2,u_3$ and $v$ instead of $x_1,x_2,x_3$ and $z$ to avoid the confusion.

There is a deep relationship between the flat coordinates of the potentials  
$F_{H_3},\>F_{(H_3)'},\>F_{(H_3)''}$ defined by
(\ref{equation:potential-H3}), (\ref{equation:ex-(H3)'case}), (\ref{equation:ex-(H3)''case}).

As is known that $\Delta_{H_3}=\det(T_{F_{H_3}})$ coincides with the discriminant of the reflection group of type $H_3$
up to a constant factor, regarding $x_1,x_2,x_3$
as basic invariants of degrees 2,6,10, respectively.

We discuss the relationship among $\Delta_{H_3},\>\Delta_{(H_3)'}(=\det(T_{F_{(H_3)'}})),\>\Delta_{(H_3)''}(=\det(T_{F_{(H_3)''}}))$.

Noting that $\Delta_{(H_3)'}$ is a polynomial of $y_1,y_3,w$ by eliminating $y_2$
in virtue of  the relation $y_2=y_1w+\frac{1}{3}w^4$,
we find by direction computation  that
$\Delta_{(H_3)'}$ coincides with $\Delta_{H_3}$
by substituting
\begin{equation}
\left\{
\begin{array}{lll}
y_1&=&\frac{1}{3}(x_1^3+3x_2),\\
y_3&=&\frac{1}{2}(x_1^2x_2+2x_3),\\
w&=&-x_1.\\
\end{array}
\right.
\end{equation}
The definition of $y_2$ shows that
$
y_2=-x_1x_2
$.
As a consequence, the flat coordinate $(y_1,y_2,y_3)$ of the potential $F_{(H_3)'}$ 
is described by the basic invariants of type $H_3$.

To treat the relationship between $\Delta_{H_3}$ and $\Delta_{(H_3)''}$,
we first note that
\begin{equation}
\begin{array}{lll}
\Delta_{H_3}&=&\frac{1}{1000}(-{x_1}^{15}-10 {x_1}^{12} {x_2}-10 {x_1}^{10} {x_3}+80 {x_1}^9 {x_2}^2+
20 {x_1}^6 {x_2}^3+100 {x_1}^5
   {x_3}^2\\
&&   -1200 {x_1}^4 {x_2}^2 {x_3}+920 {x_1}^3 {x_2}^4+1000 {x_1}^2 {x_2} {x_3}^2-1800 {x_1} {x_2}^3 {x_3}+216
   {x_2}^5+1000 {x_3}^3)\\
   \end{array}
\end{equation}
\begin{equation}
\begin{array}{ll}
&\Delta_{(H_3)''}\\
=&\frac{1}{19683}
\left(7221032
   {u_1}^9-11787048 {u_1}^7 v^2-1248156 {u_1}^6 {u_3}+4285008 {u_1}^6 v^3+5544504 {u_1}^5 v^4\right.
\\
&\left.+775656
   {u_1}^4 {u_3} v^2-101088 {u_1}^4 v^5-106434 {u_1}^3 {u_3}^2-104976 {u_1}^3 {u_3} v^3-1353240
   {u_1}^3 v^6\right.\\
   &\left.+119556 {u_1}^2 {u_3} v^4-314928 {u_1}^2 v^7+144342 {u_1} {u_3}^2 v^2+314928 {u_1}
   {u_3} v^5+101088 {u_1} v^8\right.\\
   &\left.+19683 {u_3}^3-78732 {u_3}^2 v^3-128304 {u_3} v^6-46656 v^9\right)\\
   \end{array}
\end{equation}
Then by direct computation, the coordinate transformation
\begin{equation}
\left\{
\begin{array}{lll}
{x_1}&=& {u_1}-3 v,\\
{x_2}&=& \frac{1}{27} \left(-701 {u_1}^3+297 {u_1}^2 v+1017 {u_1} v^2+108
   {u_3}-621 v^3\right),\\
   {x_3}&=& \frac{1}{270} \left(38409 {u_1}^5+64445 {u_1}^4 v-80710 {u_1}^3 v^2-5670
   {u_1}^2 {u_3}-25830 {u_1}^2 v^3-9180 {u_1} {u_3} v\right.\\
   &&\left.+93645 {u_1} v^4+9450 {u_3} v^2-36423
   v^5\right)\\
 \end{array}
 \right.
 \end{equation}
 implies that
 \begin{equation}
 \Delta_{H_3}=\frac{1}{864^2}
\left(6499 {u_1}^3+4725 {u_1}^2 v+1233 {u_1} v^2-864 {u_3}-513 v^3\right)^2 \Delta_{(H_3)''}.
\end{equation}

It is known that there is a group homomorphism of the reflection group $W(H_3)$ of type $H_3$
to the reflection group $W(D_6)$ of type $D_6$.
This suggests the existence of the family of hypersurfaces in ${\bf C}^3$ related to $W(H_3)$
consisting of the deformations of $D_6$-singularity
and actually the family is realized by T. Yano \cite{Ya}.
We employ here a family slightly different from the original one due to T. Yano.
It reads
$$
{\cal G}_{H_3}:\>z^2=f_{H_3}(x,y),
$$
where
$$
f_{H_3,(t_1,t_3,t_5)}=x^5+5xy^2+t_1x^4+\frac{4}{15}{t_1}^2x^3-\frac{5}{4}t_3x^2-t_1t_3x-\frac{5}{2}t_3y+t_5.
$$
It is straightforward to show that
by the transformation
\begin{equation}
\left\{
\begin{array}{lll}
   {x_1}&=&- \frac{1}{75}\gamma^2
   {t_1},\\
   {x_2}&=& \frac{1}{2700}\gamma (675{t_3}-4{t_1}^3),\\
 {x_3}&=& -\frac{119 }{168750}{t_1}^5+\frac{7 }{60}{t_1}^2t_3+{t_5},\\  \end{array}
   \right.
  \end{equation}
 ( where  $\gamma=-15^{3/5}$,)
  the discriminant of the family ${\cal G}_{H_3}$
  coincides with $\Delta_{H_3}$.
  
  If $f_{H_3}$ is the family corresponding to $\Delta_{H_3}$, what is the family corresponding to $\Delta_{(H_3)''}$?
  As an answer to this question, we construct the family
 $$
 {\cal G}_{(H_3)''}:\>z^2=f_{(H_3)''},
 $$
 where
$$
\begin{array}{lll}
f_{(H_3)'',(r_1,r_3,s_1)}&=&\frac{3}{2} {s_1}^3 \left(40 {r_1}^2-25 {r_1} {s_1}+4
   {s_1}^2\right)-{r_3}\{(-3 {s_1} (4 {r_1}-{s_1})-2
   y\}\\
   &&+\frac{1}{2} x (6 {r_1}+x) \left(5 {r_1} x^2+2 {r_3}+2
   x^3\right)-{r_1} x^4+5 x y^2\\
   \end{array}
   $$
We don't show the relationship between the discriminant of the family $f_{(H_3)''}$ and $\Delta_{(H_3)''}$
and is left to the reader.

\subsection{$H_4$ case}
 \label{subsection:sec5-2}

It is known that there is a group homomorphism of $W(H_4)$ to $W(E_8)$.
This suggests  the existence of a family of hypersurfaces which  are deformations of the 
hypersurface of $E_8$-singularity corresponding to the algebraic potential $F_{H_4(j)}$.
In this subsection, we construct such families only
for the cases $F_{H_4(0)}$ and $F_{H_4(1)}$.
The result is given in the following.

\begin{theorem}
The families
$$
\begin{array}{lll}
{\cal G}_{H_4(0)}&:&z^2=f_{H_4(0)},\\
{\cal G}_{H_4(1)}&:&z^2=f_{H_4(1)},\\
\end{array}
$$
where
{\small
$$
\begin{array}{lll}
f_{H_4(0)}&=&x^5+{t_{15}}+5 {t_6}^2 x+5 {t_6} x^3+
y \left(10 {t_1} {t_6} x+5 {t_1} x^3+{t_{10}}\right)+5{t_1}^2xy^2+y^3,\\
f_{H_4(1)}&=&x \left(5 {t_3}^2+5 {t_3} x+x^2\right)^2+y \left\{5 {t_1}
   ({t_3}+x) \left(5 {t_3}^2+5 {t_3}
   x+x^2\right)+{t_{10}}\right\}+y^2 \left(5 {t_1}^2
   x+{t_5}\right)+y^3,\\
   \end{array}
   $$
   }
are those whose discriminants coincide with $\Delta_{H_4(0)}=\det(T_{F_{H_4(0)}})$
and
 $\Delta_{H_4(1)}=\det(T_{F_{H_4(1)}})$, respectively.
\end{theorem}

We are going to explain not only the outline of the proof of this theorem
but also its consequences.

We first treat the case $H_4(0)$.
Put $f=f_{H_4(0)}$ and let $S$ be the hypersurface of ${\bf C}^3$ defined by $z^2=f(x,y)$.
Our argument below is parallel to that in the proof of Theorem \ref{theorem:theorem-4}.
Let $C_f$ be the singular locus of $S_0=\{(x,y)\in{\bf C}^2|\>f(x,y)=0\}$, namely,
$$
C_f=\{(x,y)\in{\bf C}^2|\>f=f_x=f_y=0\}.
$$
There is a polynomial $\delta(t)$ of 
the parameter $t=(t_1,t_6,t_{10},t_{15})$ such that $\{t\in{\bf C}^4|\>C_f\not=\emptyset\}$ 
coincides with $\{t\in{\bf C}^4|\>\delta(t)=0\}$.
To obtain $\delta(t)$,
we introduce polynomials $g_j\>(j=1,2,3)$ of $x,y$ with coefficients $a_i,\>b_j,\>c_k$
and define $H=g_1f+g_2f_x+g_3f_y$ as in the proof of Theorem \ref{theorem:theorem-4}.
Regarding $H$ as a polynomial of $x,y$, we eliminate the constants $b_j,\>c_k$
so that $H$ takes the form
$$p_1+p_2x+p_3y+p_4x^2+p_5xy+p_6x^3+p_7x^2y+p_8x^3y.
$$
From this condition, we find that each $p_j$ is a linear combination of $a_1,\ldots,a_8$ whose coefficients are
polynomials of $t_1,t_6,t_{10},t_{15}$.
Namely, there are polynomials $p_{ij}\>(i,j=1,\ldots,8)$ of 
 $t_1,t_6,t_{10},t_{15}$ such that
 $$
 p_j=\sum_{j=1}^8p_{ij}a_j.
 $$
 Let $P_{H_4(0)}$ be the $8\times 8$ matrix whose $(i,j)$ entry is $p_{ij}$.
 Then it follows that for the parameter $(t_1,t_6,t_{10},t_{15})$,
  $C_f\not=\emptyset$ if and only if $\det(P_{H_4(0)})=0$.
 On the other hand, by direct computation, we find that
$
 \det(P_{H_4(0)})$
 coincides with
 $(D_{H_4(0)})^2$ up to a non-zero constant factor, where
 {\small
 $$
 \begin{array}{ll}
 &D_{H_4(0)}\\
=& -3125 {t_1}^{12} {t_{10}}^3 {t_6}^3-27 {t_1}^{10} {t_{10}}^5-9375 {t_1}^{10} {t_{10}}^2 {t_6}^5-105 {t_1}^8 {t_{10}}^4
   {t_6}^2-9375 {t_1}^8 {t_{10}} {t_6}^7\\
   &+13695 {t_1}^6 {t_{10}}^3 {t_6}^4-3125 {t_1}^6 {t_6}^9+120 {t_1}^4 {t_{10}}^5
   {t_6}+36885 {t_1}^4 {t_{10}}^2 {t_6}^6+2280 {t_1}^2 {t_{10}}^4 {t_6}^3\\
   &+35640 {t_1}^2 {t_{10}} {t_6}^8+16 {t_{10}}^6-864
   {t_{10}}^3 {t_6}^5+11664 {t_6}^{10}\\
& +  3 {t_1} t_{15}\left(3125 {t_1}^{12} {t_{10}}^2 {t_6}^2+6250 {t_1}^{10} {t_{10}} {t_6}^4+25
   {t_1}^8 {t_{10}}^3 {t_6}+3125 {t_1}^8 {t_6}^6-14150 {t_1}^6 {t_{10}}^2 {t_6}^3-4 {t_1}^4 {t_{10}}^4\right.\\
   &\left.-24075 {t_1}^4
   {t_{10}} {t_6}^5-1980 {t_1}^2 {t_{10}}^3 {t_6}^2-11700 {t_1}^2 {t_6}^7+4320 {t_{10}}^2 {t_6}^4\right)\\
   &-3{ t_{15}}^2\left(3125 {t_1}^{14}
   {t_{10}} {t_6}+3125 {t_1}^{12} {t_6}^3+125 {t_1}^{10} {t_{10}}^2-14625 {t_1}^8 {t_{10}} {t_6}^2-11250 {t_1}^6 {t_6}^4-2430
   {t_1}^4 {t_{10}}^2 {t_6}\right.\\
   &\left.+7830 {t_1}^2 {t_{10}} {t_6}^3-72 {t_{10}}^3-1944 {t_6}^5\right)\\
 & + 25 {t_1}^3{t_{15}}^3 \left(125 {t_1}^{12}-675
   {t_1}^6 {t_6}-81 {t_1}^2 {t_{10}}+729 {t_6}^2\right)+729{t_{15}}^4.
   \\
   \end{array}
   $$
   }
We recall the definition of the polynomial potential $T_{H_4(0)}$ given in
\S\ref{subsubsection:sec3-3-1}.
Then $\Delta_{H_4(0)}=\det(T_{F_{H_4(0)}})$ is a polynomial of $(x_1,x_2,x_3,x_4)$.
It is straightforward to show that by the coordinate transformation
$$
\left\{
\begin{array}{lll} 
 {x_1}&=& {t_1},\\
 {x_2}&=& \frac{4}{225} \left(10 {t_1}^6-27
   {t_6}\right),\\
   {x_3}&=& \frac{2}{50625} \left(500 {t_1}^{10}-2025 {t_1}^4 {t_6}-486
   {t_{10}}\right),\\
   {x_4}&=& \frac{8}{11390625} \left(8125 {t_1}^{15}-47250 {t_1}^9
   {t_6}-7290 {t_1}^5 {t_{10}}+54675 {t_1}^3 {t_6}^2+13122
   {t_{15}}\right),\\
   \end{array}
   \right.
   $$
 $\Delta_{H_4(0)}$ coincides with $D_{H_4(0)}$ up to a constant factor.
 
We next treat the case $H_4(1)$. Put $f=f_{H_4(1)}$ and
let $C_f$ the singular locus of $S_0=\{(x,y)\in{\bf C}^2|\>f(x,y)=0\}$
similar the the case $H_4(0)$.
Then
$$
C_f=\{(x,y)\in{\bf C}^2|\>f=f_x=f_y=0\}.
$$
There is a polynomial $\delta(t)$ of 
the parameter $t=(t_1,t_3,t_5,t_{10})$ such that $\{t\in{\bf C}^4|\>C_f\not=\emptyset\}$ 
coincides with $\{t\in{\bf C}^4|\>\delta(t)=0\}$.
By an argument similar to that in the case $H_4(0)$,
we find that the defining equation of 
 the set $\{t\in{\bf C}^4|\>C_f\not=\emptyset\}$ is
 $t_{10}D_{H_4(1)}=0$,
 where
$$
\begin{array}{ll} 
&D_{H_4(1)}\\
=& {t_3}^5 \left(25 {t_1}^2 {t_3}-4 {t_5}\right)^2
\left(125 {t_1}^9
   {t_3}^2-25 {t_1}^7 {t_3} {t_5}-1225 {t_1}^6 {t_3}^3+{t_1}^5
   {t_5}^2+370 {t_1}^4 {t_3}^2 {t_5}\right.\\
   &\left.+335 {t_1}^3 {t_3}^4-35 {t_1}^2
   {t_3} {t_5}^2-45 {t_1} {t_3}^3 {t_5}-27
   {t_3}^5+{t_5}^3\right)\\
&   -{t_3}^3 \left(10625 {t_1}^{12} {t_3}^3-4750
   {t_1}^{10} {t_3}^2 {t_5}-91750 {t_1}^9 {t_3}^4+605 {t_1}^8 {t_3}
   {t_5}^2+54350 {t_1}^7 {t_3}^3 {t_5}\right.\\
   &\left.-83175 {t_1}^6 {t_3}^5-20 {t_1}^6
   {t_5}^3-10548 {t_1}^5 {t_3}^2 {t_5}^2+25600 {t_1}^4 {t_3}^4 {t_5}+22520
   {t_1}^3 {t_3}^6+805 {t_1}^3 {t_3} {t_5}^3\right.\\
   &\left.-2930 {t_1}^2 {t_3}^3
   {t_5}^2-2880 {t_1} {t_3}^5 {t_5}-20 {t_1} {t_5}^4-1728 {t_3}^7+136
   {t_3}^2 {t_5}^3\right)t_{10}\\
& +(  -1225 {t_1}^{11} {t_3}^3+255 {t_1}^9 {t_3}^2
   {t_5}+13895 {t_1}^8 {t_3}^4-30 {t_1}^7 {t_3} {t_5}^2-3900 {t_1}^6
   {t_3}^3 {t_5}
   -19976 {t_1}^5 {t_3}^5\\
   &+{t_1}^5 {t_5}^3+550 {t_1}^4
   {t_3}^2 {t_5}^2+3700 {t_1}^3 {t_3}^4 {t_5}+840 {t_1}^2 {t_3}^6-40
   {t_1}^2 {t_3} {t_5}^3-200 {t_1} {t_3}^3 {t_5}^2\\
   &+288 {t_3}^5
   {t_5}+{t_5}^4){t_{10}}^2\\
   &+(-27 {t_1}^{10}+360 {t_1}^7 {t_3}-36 {t_1}^5 {t_5}-920
   {t_1}^4 {t_3}^2+160 {t_1}^2 {t_3} {t_5}+480 {t_1} {t_3}^3-8
   {t_5}^2){t_{10}}^3+
   16{t_{10}}^4.\\
   \end{array}
   $$
 It is straightforward to show that
 by the coordinate transformation  
   $$
   \begin{array}{lll}
 {x_1}&=& {t_1},\\
 {z}&=& 8 {t_3},\\
 {x_5}&=& \frac{3}{5} \left(17
   {t_1}^5-160 {t_1}^2 {t_3}+16 {t_5}\right),\\
   {x_{10}}&=& \frac{6}{25}
   \left(987 {t_1}^{10}-15120 {t_1}^7 {t_3}+1512 {t_1}^5 {t_5}+50400
   {t_1}^4 {t_3}^2-8960 {t_1}^2 {t_3} {t_5}-12800 {t_1} {t_3}^3\right.\\
   &&\left.-2560
   {t_{10}}+448 {t_5}^2\right),
\\
\end{array}
 $$
 $\Delta_{H_4(1)}$ coincides with $D_{H_4(1)}$ up to a constant factor.
   
  We now explain the relationship between $D_{H_4(1)}$ and $D_{H_4(0)}$.
    
   Put
   $$
   \begin{array}{lll}
   x&=&x_1-2t_3,\\
   y&=&y_1+\frac{1}{3}(10{t_1}^2t_3-t_5).\\
   \end{array}
   $$
   Then $f_{H_4(1),(t_1,t_3,t_3,t_{10})}(x,y)$
   is transformed to
   $ f_{H_4(0),(u_1,u_6,u_{10},u_{15})}(x_1,y_1)$ by the substitution
 \begin{equation}
 \label{equation:H4-u-t}
 \left\{
 \begin{array}{lll}  
  u_1&=&    {t_1},\\
  u_6&=&\frac{1}{3} \left(10 {t_1}^3 {t_3}-{t_1} {t_5}-3 {t_3}^2\right),\\
  u_{10}&=&\frac{1}{3} \left(-100
   {t_1}^4 {t_3}^2+20 {t_1}^2 {t_3} {t_5}+15 {t_1} {t_3}^3+3
   {t_{10}}-{t_5}^2\right),\\
   u_{15}&=&\frac{1}{27} \left(-2000 {t_1}^6 {t_3}^3+600 {t_1}^4 {t_3}^2 {t_5}+450
   {t_1}^3 {t_3}^4+90 {t_1}^2 {t_{10}} {t_3}-60 {t_1}^2 {t_3} {t_5}^2\right.\\
   &&\left.-45 {t_1} {t_3}^3
   {t_5}-9 {t_{10}} {t_5}-54 {t_3}^5+2 {t_5}^3\right).\\
   \end{array}
   \right.
   \end{equation}
Moreover, by direct computation, we find that
$D_{H_4(0)}(u_1,u_6,u_{10},u_{15})$ coincides with\\
 $t_{10}^2D_{H_4(1)}(t_1,t_3,t_5,t_{10})$ up to a constant factor by the transformation
   (\ref{equation:H4-u-t}).

\subsection{More concerning $H_4$}

There are at least seven algebraic potentials $H_4(j)\>(j=0,1,2,3,6,7,9)$ are knwon
related to the reflection group of type $H_4$ (cf. \S3.3.1).
We observe that for each of pairs $(F_{H_4(0)},F_{H_4(9)})$,  $(F_{H_4(1)},F_{H_4(7)})$,  $(F_{H_4(2)},F_{H_4(6)})$,
the flat coordinates of the latter are expressed as polynomials of the former.
The results are summarized as follows.

\vspace{5mm}
\underline{{\bf (Case 1):} $H_4(0)$ and $H_4(9)$}

Let $F_{H_4(0)}$ be the polynomial potential of the reflection group of type $H_4$ (cf. \S3.3.1)
and let $F_{H_{4(9)}}$ be the algebraic potential defined in \S3.3.1.
We change of the variables of $F_{H_{4(9)}}$ by $y_1,y_2,y_3,y_4,$ and the algebraic function of $w$
instead of $x_1,x_2,x_3,x_4$ and $z$, respectively.
Then, by the transformation
$$
\begin{array}{lll}
y_1&=&\frac{1}{6}x_1(4{x_1}^6+31x_2),\\
y_2&=&\frac{3675}{8}({x_1}^4x_2+18x_3),\\
y_4&=&\frac{46305}{32}(5{x_1}^3{x_2}^2+8{x_1}^5x_3+90x_4),\\
w&=&-x_1.\\
y_3&=&\frac{147}{16}(2{x_1}^6x_2+45{x_2}^2+180{x_1}^2x_3),\\
\end{array}
$$
$\det(T_{H_4(9)})$ coincides with
$\det(T_{H_4(0)})$ up to a constant factor.

\underline{{\bf (Case 2):} $H_4(1)$ and $H_4(7)$}

Let $F_{H_4(1)}$ be the algebraic potential of type $H_{4(1)}$
and let $F_{H_{4(7)}}$ be the algebraic potential type $H_{4(7)}$ defined in \S3.3.1.
We change of the variables of $F_{H_{4(7)}}$ by $y_1,y_2,y_3,y_4,$ and the algebraic function of $w$
instead of $x_1,x_2,x_3,x_4$ and $z$, respectively.
Then, by the transformation
$$
\begin{array}{lll}
{y_1}&=& \frac{1}{81} \left(-4{x_1}^3-3 {z}\right),\\
{y_2}&=& -\frac{3}{800} \left(2 {x_1}^5+2x_2-{x_1}^2{z}\right),\\
{y_4}&=&
 -  \frac{9}{5000} \left(24 {x_1}^{10}-36{x_1}^5 {x_2}-12  {x_2} ^2+x_4+2{x_1}^7{z}-2 {x_1}^2 {x_2} {z}-
   5 {x_1}^4{z}- {x_1}{z}^3\right),\\
   w&=&x_1\\
   \end{array}
   $$
$\det(T_{H_4(7)})$ coincides with
$\det(T_{H_4(1)})$ up to a constant factor.

\underline{{\bf (Case 3):} $H_4(2)$ and $H_4(6)$}

Let $F_{H_4(2)}$ be the algebraic potential of type $H_{4(2)}$
and let $F_{H_{4(6)}}$ be the algebraic potential type $H_{4(6)}$ defined in \S3.3.1.
We change of the variables of $F_{H_{4(6)}}$ by $y_1,y_2,y_3,y_4,$ and the algebraic function of $w$
instead of $x_1,x_2,x_3,x_4$ and $z$, respectively.
Then, by the transformation
$$
\begin{array}{lll}
{y_1}&=& -\frac{5}{12} \left({x_1}^2-2 {z}\right),\\
{y_2}&=& \frac{125}{48} {x_2} \left(8 {x_1}^2-{z}\right),\\
{y_4}&=&
   \frac{15625}{6912} \left(168 {x_1}^5 {x_2}+228 {x_1}^3 {x_2} {z}-48 {x_1} {x_2} {z}^2-27 {x_2}^3+4 {x_4}\right),\\
   w&=&
   -\frac{5 {x_1}}{2}\\
   \end{array}
   $$
$\det(T_{H_4(6)})$ coincides with
$\det(T_{H_4(2)})$ up to a constant factor.

\vspace{5mm}

{\large{\bf Appendix:
Comparison between the flat coordinate $(x_1,\ldots,x_6,x_8,z)$ and\\
 $(t_1,t_3,t_4,s_5,t_6,t_7,t_9,t_{12})$}}

In this appendix, we explain how to obtain the transformation
(\ref{equation:trans-e8-case}) for the sake of completeness.

We define the coordinate transformation of 
 $(t_1,t_3,t_4,s_5,t_6,t_7,t_9,t_{12})$ to
 $(x_1,\ldots,x_6,x_8,z)$ by
\begin{equation}
\label{equation:trans-t-x}
\left\{
\begin{array}{lll}
x_1&=&d_1t_1,\\
x_2&=&a_1t_3+d_2{t_1}^3,\\
x_3&=&a_2t_4+d_3t_1t_3+d_4{t_1}^4,\\
x_4&=&a_3t_6+b_1{t_3}^2+d_5{t_1}^2t_4+d_6{t_1}^3t_3+d_7{t_1}^6+h_1s_5t_1,\\
x_5&=&a_4t_7+b_2{t_3}t_4+d_8t_1t_6+d_9{t_1}^3t_4+d_{10}t_1{t_3}^2+d_{11}{t_1}^4t_3+d_{12}{t_1}^7+h_2s_5{t_1}^2,\\
x_6&=&a_5t_9+b_3t_3t_6+b_4{t_3}^3+d_{13}{t_1}^2t_7+d_{14}{t_1}^3t_6+
d_{15}{t_1}{t_4}^2+d_{16}{t_1}^2t_3t_4\\
&&+d_{17}{t_1}^5t_4
+d_{18}{t_1}^3{t_3}^3+d_{19}{t_1}^6t_3+d_{20}{t_1}^9
+s_5(h_3t_4+h_4t_1t_3+h_5{t_1}^4),\\
x_8&=&t_{12}+b_5t_3t_9+b_6{t_6}^2+b_7{t_3}^2t_6+b_8{t_4}^3+b_9{t_3}^4
+d_{21}{t_1}^3t_9+d_{22}t_1t_4t_7+d_{23}{t_1}^2t_3t_7\\
&&+d_{24}{t_1}^5t_7+
d_{25}{t_1}^2t_4t_6+d_{26}{t_1}^3t_4t_6+d_{27}{t_1}^6t_6
+d_{28}t_1t_3{t_4}^2+d_{29}{t_1}^4{t_4}^2\\
&&+d_{30}{t_1}^2{t_3}^2t_4+d_{31}{t_1}^5t_3t_4+d_{32}{t_1}^8t_4
+d_{33}{t_1}^3{t_3}^3
+d_{34}{t_1}^6{t_3}^2+d_{35}{t_1}^9t_3+d_{36}{t_1}^{12}\\
&&+
s_5(h_6s_5{t_1}^2+h_8t_1t_6+h_9t_3t_4+h_{10}{t_1}^3t_4+h_{11}t_1{t_3}^2+h_{12}{t_1}^4t_3+h_{13}{t_1}^7),\\
z&=&h_0s_5,\\
\end{array}
\right.
\end{equation}
where $a_i,b_j,d_k,h_l$ are constants to be determined.
By the condition that (\ref{equation:trans-t-x}) actually defines a transformation,
we assume that each of $a_1,\dots,a_5,\>h_0$ are not zero.
As the normalization, we assume that the coefficient of $t_{12}$ is 1.

Let $T$ be the $8\times 8$ matrix defined by the potential for $E_{8(1)}$.
Then $\det(T)$ is regarded as a polynomial
of  $(x_1,\ldots,x_6,x_8,z)$.
By the transformation
(\ref{equation:trans-t-x}), $\det(T)$ turns out to be a polynomial of
$(t_1,t_3,t_4,s_5,t_6,t_7,t_9,t_{12})$ whose coefficients are polynomials of $a_i,b_j,d_k,h_l$.

We are going to determine the constants $a_i,b_j,d_k$ so that the identity equation
\begin{equation}
\label{equation:id-eq-T-P}
c_0\det(T)=\det(P_C)
\end{equation}
holds, where $c_0$ is a non-zero constant.

(Step A-1) Consider the case $t_1=s_5=0$.

For a moment we assume $t_1=s_5=0$ without any comment.
It follows from (\ref{equation:trans-t-x}) that
the transformation of $(t_3,t_4,t_6,t_7,t_9,t_{12})$ to  $(x_2,\ldots,x_6,x_8)$ takes the form
\begin{equation}
\label{equation:trans-t-x-x1s5}
\left\{
\begin{array}{lll}
x_2&=&a_1t_3,\\
x_3&=&a_2t_4,\\
x_4&=&a_3t_6+b_1{t_3}^2,\\
x_5&=&a_4t_7+b_2{t_3}t_4,\\
x_6&=&a_5t_9+b_3t_3t_6+b_4{t_3}^3,\\
x_8&=&t_{12}+b_5t_3t_9+b_6{t_6}^2+b_7{t_3}^2t_6+b_8{t_4}^3+b_9{t_3}^4.\\
\end{array}
\right.
\end{equation}
We first observe that $\det(T)$ has a factor of the form
$$
L(x)=x_8+r_1x_2x_6+r_2{x_4}^2+r_3{x_2}^2x_4+r_4{x_3}^3+r_4{x_2}^4
$$
by substituting special values in $\det(T)$ and actually
$$
L_0(x)=x_8-10x_2x_6-5{x_4}^2+100{x_2}^2x_4-\frac{5}{12}{x_3}^3-250{x_2}^4
$$
is the candidate of a factor of $\det(T)$.
On the other hand, we also observe that
$\det(P_E)$ is  divided by $t_{12}$.
Then
$L_0(x)$ coincides with $t_{12}$ under the condition
(\ref{equation:trans-t-x-x1s5}).
Then we have
\begin{equation}
\left\{
\begin{array}{lll}
b_5&=&10a_1a_5,\\
b_6&=&5{a_3}^2,\\
b_7&=&-10(10{a_1}^2a_3-a_3b_1-a_1b_3),\\
b_8&=&\frac{5}{12}{a_2}^3,\\
b_9&=&5(50{a_1}^4-20{a_1}^2b_1+{b_1}^2+2a_1b_4).\\
\end{array}
\right.
\end{equation}
We normalize $\det(P_E)$ so that
$$
c_0=\frac{65536}{1953125}.
$$
By the substitution $t_3=t_4=t_6=t_9=0$ in (\ref{equation:id-eq-T-P}),
we obtain
\begin{equation}
\label{equaton:def-a4}
{a_4}^{12}=-\frac{16}{243}.
\end{equation}
By the substitution $t_3=t_4=t_6=0$ in (\ref{equation:id-eq-T-P}), we obtain
$
a_5=-\frac{3}{4}{a_4}^3.
$
By the substitution $t_3=t_4=0$ in (\ref{equation:id-eq-T-P}), we obtain
$
a_3=\frac{9}{16}{a_4}^6.
$
By the substitution $t_3=0$ in (\ref{equation:id-eq-T-P}), we obtain
$
a_2=-\frac{3}{4}{a_4}^4.
$
By the substitution $t_4=t_6=t_7=0$  in (\ref{equation:id-eq-T-P}) under ${t_3}^2\equiv 0$, we obtain
$a_1=-\frac{27}{64}{a_4}^9$.
By the substitution $t_4=t_6=t_7=0$  in (\ref{equation:id-eq-T-P}) under ${t_3}^3\equiv 0$, we obtain
$b_1=-\frac{51}{256}{a_4}^6$.
By the substitution $t_4=t_6=t_7=0$  in (\ref{equation:id-eq-T-P}) under ${t_3}^4\equiv 0$, we obtain
$b_4=-\frac{385}{4608}{a_4}^3$.
By the substitution $t_4=t_7=t_9=0$  in (\ref{equation:id-eq-T-P}) under ${t_3}^3\equiv 0$, we obtain
$b_3=\frac{11}{32}{a_4}^3$.
By the substitution $t_6=t_7=t_9=0$  in (\ref{equation:id-eq-T-P}) under ${t_3}^5\equiv 0$, we obtain
$b_2=-\frac{7}{24}{a_4}$.

Accordingly, all the constants $a_1,\ldots,a_5,b_1,\ldots,b_9$ are determined.

(Step A-2) Consider the case $s_5=0$.

For a moment we assume $s_5=0$ without any comment.
Then the constants $h_l$'s do not appear  in (\ref{equation:trans-t-x})
 and our job is reduced to determine $d_k$'s.

In the below, we explain the procedure of the determination of $d_k$'s.
Our aim is performed by  the substitution of $t_3,t_4,\ldots$ in  (\ref{equation:id-eq-T-P}) 
under the condition ${t_3}^n\equiv 0$ , where $n$ is an appropriate integer.

\begin{tabular}{lll}
Substitute $t_3=t_4=t_6=t_7=0$ under ${t_1}^2\equiv 0$ &$\Longrightarrow$&$ d_1$\\
Substitute $t_3=t_4=t_6=t_7=0$ under ${t_1}^4\equiv 0$ &$\Longrightarrow$&$ d_2,\>d_{21}$\\
Substitute $t_3=t_4=t_6=t_7=0$ under ${t_1}^7\equiv 0$ &$\Longrightarrow$&$ d_4,\>d_{7}$\\
Substitute $t_3=t_4=t_6=t_7=0$ under ${t_1}^{10}\equiv 0$ &$\Longrightarrow$&$ d_{12},\>d_{20}$\\
Substitute $t_3=t_4=t_6=t_7=0$ under ${t_1}^{13}\equiv 0$ &$\Longrightarrow$&$ d_{36}$\\
Substitute $t_3=t_4=0,t_6=t_7=t_9=1$ under ${t_1}^{2}\equiv 0$ &$\Longrightarrow$&$ d_{8}$\\
Substitute $t_3=t_6=0,t_4=t_7=t_9=1$ under ${t_1}^{2}\equiv 0$ &$\Longrightarrow$&$ d_{22},\>d_{15}$\\
Substitute $t_4=t_6=0,t_3=t_7=t_9=1$ under ${t_1}^{2}\equiv 0$ &$\Longrightarrow$&$ d_{3},\>d_{10}$\\
Substitute $t_3=t_4=t_6=t_7=t_9=1$ under ${t_1}^{2}\equiv 0$ &$\Longrightarrow$&$ d_{28}$\\
Substitute $t_3=t_4=t_6=0,t_7=t_9=1$ under ${t_1}^{3}\equiv 0$ &$\Longrightarrow$&$ d_{13}$\\
Substitute $t_3=t_6=t_7=0,t_4=t_9=1$ under ${t_1}^{3}\equiv 0$ &$\Longrightarrow$&$ d_{5}$\\
Substitute $t_4=t_6=0,t_3=t_7=t_9=1$ under ${t_1}^{3}\equiv 0$ &$\Longrightarrow$&$ d_{23}$\\
Substitute $t_3=t_7=0,t_4=t_6=t_9=1$ under ${t_1}^{3}\equiv 0$ &$\Longrightarrow$&$ d_{25}$\\
Substitute $t_6=t_7=0,t_3=t_4=t_9=1$ under ${t_1}^{3}\equiv 0$ &$\Longrightarrow$&$ d_{16},\>d_{30}$\\
Substitute $t_4=t_6=t_7=t_9=0,t_3=1$ under ${t_1}^{4}\equiv 0$ &$\Longrightarrow$&$ d_6,d_{18},\>d_{33}$\\
Substitute $t_3=t_4=t_7=0,t_6=t_9=1$ under ${t_1}^{4}\equiv 0$ &$\Longrightarrow$&$ d_{14}$\\
Substitute $t_3=t_6=t_7=0,t_4=t_9=1$ under ${t_1}^{4}\equiv 0$ &$\Longrightarrow$&$ d_{9}$\\
Substitute $t_4=t_7=t_9=0,t_3=t_6=1$ under ${t_1}^{4}\equiv 0$ &$\Longrightarrow$&$ d_{26}$\\
Substitute $t_3=t_6=t_7=0,t_4=1$ under ${t_1}^{5}\equiv 0$ &$\Longrightarrow$&$ d_{29}$\\
Substitute $t_4=t_6=0,t_3=t_7=1$ under ${t_1}^{5}\equiv 0$ &$\Longrightarrow$&$ d_{11}$\\
Substitute $t_3=t_4=t_6=t_7=t_9=1$ under ${t_1}^{5}\equiv 0$ &$\Longrightarrow$&$ d_{17},\>d_{31}$\\
Substitute $t_3=2,t_4=t_6=t_7=t_9=1$ under ${t_1}^{5}\equiv 0$ &$\Longrightarrow$&$ d_{24}$\\
Substitute $t_3=t_6=1,t_4=t_7=t_9=1$ under ${t_1}^{6}\equiv 0$ &$\Longrightarrow$&$ d_{19},d_{34}$\\
Substitute $t_3=2,t_6=1,t_4=t_7=t_9=0$ under ${t_1}^{6}\equiv 0$ &$\Longrightarrow$&$ d_{27}$\\
Substitute $t_3=t_6=t_7=t_9=0,t_4=1$ under ${t_1}^{8}\equiv 0$ &$\Longrightarrow$&$ d_{32}$\\
\end{tabular}

At last we compute both sides of  (\ref{equation:id-eq-T-P}) under the condition ${t_1}^{10}\equiv 0$.
Then the value of  $d_{35}$ is determined.

(Step A-3) The general case.

The constants $h_l$'s are still not determined.
The computation of the determination of the constants $h_l$'s are performed by the similar manner to (Step A-2).

\begin{tabular}{lll}
Substitute $t_3=t_4=t_6=t_9=0,t_7=1$ under ${t_1}\equiv 0,\>{s_5}^2\equiv 0$ &$\Longrightarrow$&$h_0,\>h_{7}$\\
Substitute $t_3=t_4=1,t_6=t_7=t_9=1$ under ${t_1}\equiv 0,\>{s_5}^2\equiv 0$ &$\Longrightarrow$&$h_3,\>h_{9}$\\
Substitute $t_3=1,t_4=t_6=t_7=t_9=1$ under ${t_1}^2\equiv 0,\>{s_5}^2\equiv 0$ &$\Longrightarrow$&$h_4,\>h_{11}$\\
Substitute $t_6=1,t_3=t_4=t_7=t_9=1$ under ${t_1}^2\equiv 0,\>{s_5}^2\equiv 0$ &$\Longrightarrow$&$h_1,\>h_{8}$\\
Substitute $t_4=1,t_3=t_6=t_7=t_9=1$ under ${t_1}^4\equiv 0,\>{s_5}^2\equiv 0$ &$\Longrightarrow$&$h_2,\>h_{10}$\\
Substitute $t_3=1,t_4=t_6=t_7=t_9=1$ under ${t_1}^5\equiv 0,\>{s_5}^2\equiv 0$ &$\Longrightarrow$&$h_5,\>h_{12}$\\
Substitute $t_3=t_4=t_6=t_7=t_9=0$ under ${t_1}^8\equiv 0,\>{s_5}^2\equiv 0$ &$\Longrightarrow$&$h_{13}$\\
Substitute $t_3=t_4=t_6=t_7=t_9=0$ under ${t_1}^3\equiv 0,\>{s_5}^5\equiv 0$ &$\Longrightarrow$&$h_{6}$\\
\end{tabular}

The transformation given in Theorem 2 (iii) is the same as the one  obtained by the argument above
under the condition $c_8=a_4$ and 
(\ref{equation:trans-e8-case}) holds if one assumes
the identity equation (\ref{equation:id-eq-T-P}).

\vspace{5mm}
{\bf Acknowledgements}

The author thanks Professor J. Michel for teaching him the conjecture of Dubrovin on the algebraic potentials
and
T. Douvropoulos for sending him the data on primitive conjugacy classes which is very helpful
to constuct algebraic potentials.
He also thanks the organizers of  ``the workshop
on hyperplane arrangements and reflection groups'' held at University of  Hannover, 2019
for giving him the chance to talk.
This work was partially supported by JSPS KAKENHI Grant Number 17K05269.


\end{document}